\documentclass[10pt]{article}

\usepackage{a4wide}
\usepackage{amssymb}
\usepackage{amsfonts}
\usepackage{amsmath}
\input xy
\xyoption{arrow} \xyoption{matrix}

\date{}

\newtheorem{proposition}{Proposition}[section]
\newtheorem{theorem}[proposition]{Theorem}
\newtheorem{lemma}[proposition]{Lemma}

\newtheorem{corollary}[proposition]{Corollary}

\def\GK{{\rm  GK}\,}
\def\Kdim{{\rm K.dim }\,}

\def\Hom{{\rm Hom}}
\def\der{\partial }

\def\nFM0{{\nu }_{F,M_0}}
\def\nFN0{{\nu }_{F,N_0}}
\def\nGN0{{\nu }_{G,N_0}}

\def\N0{ {\bf N}_0 }

\def\t{\otimes}
\def\g{\gamma}
\def\v{\varphi}
\def\ra{\rightarrow}

\def\lra{\leftrightarrow}
\def\Xpm{X^{\pm }}

\def\s{\sigma}
\def\Z{\mathbb{Z}}

\def\l1{{\lambda}_1}

\def\a{\alpha}
\def\a0{ {\alpha }_0}
\def\a1{ {\alpha }_1}

\def\l{\lambda}
\def\o{\omega}

\def\nFGM0{{\nu }_{F,G,M_0}}

\def\nFN0{{\nu}_{F,N_0}}

\def\sm{{\sigma}^m}

\def\sm1{{\sigma}^{-1}}

\def\smtp1{{\sigma}^{-t+1}}

\def\o{\omega }
\def\S1{S^{-1}}

\def\Xpm1{X^{\pm 1}_1}

\def\sPM1{{\sigma }^{\pm 1}}
\def\sMP1{{\sigma }^{\mp 1 }}

\def\b{\beta}
\def\d{\delta}

\def\di{{\rm d.ind}}

\def\M{{\cal M}}
\def\L{\Lambda}
\def\O{\Omega}

\def\G{\Gamma}

\def\OO{{\cal O}}
\def\CA{{\cal A}}

\def\CD{{\cal D}}

\def\Ytm1{Y^{t-1}}
\def\Yim1{Y^{i-1}}

\def\CM{{\cal M}}

\def\CG{{\cal G}}
\def\CH{{\cal H}}
\def\ass{{\rm ass}}

\def\i{{\bf i}}

\def\Aut{{\rm Aut}}

\def\Der{{\rm Der }}
\def\ad{{\rm ad }}
\def\dim{{\rm dim }}

\def\pad{{\rm pad }}

\def\bj{\bar{j}}

\def\ker{ {\rm ker } }

\def\CJ{ {\cal J}}

\def\D{ \Delta }

\def\SL2Z{ {\rm SL}_2({\bf Z}) }

\def\CR{ {\cal R}}

\def\th{ \theta }

\def\Gp1{ G^{1 , 1 } }
\def\P11{ P^{-1 , 1 } }
\def\Pp1{ P^{1 , 1 } }

\def\bP{ \bar{P}}

\def\th{\theta}

\def\CV{{\cal V}}

\def\nCLsr{{}^\nu\kern-2pt {\cal L}^{\sigma , \rho  }}
\def\nP{{}^\nu \kern-2pt P}
\def\nL{{}^\nu\kern-2pt L}
\def\nLL{{}^\nu\kern-2pt \Lambda}
\def\nPsr{{}^\nu\kern-2pt P^{\sigma , \rho  }}
\def\nLsr{{}^\nu\kern-2pt L^{\sigma , \rho  }}
\def\nuCL{{}^\nu\kern-2pt  {\cal L}}
\def\nCLsr{{}^\nu\kern-2pt {\cal L}^{\sigma , \rho  }}
\def\nCL1m{{}^\nu\kern-2pt {\cal L}^{-1 , 1  }}
\def\x1nu{x^\frac{1}{\nu}}
\def\xm1nu{x^{-\frac{1}{\nu}}}

\def\trdeg{{\rm tr.deg}}

\def\CR{ {\cal R}}

\def\ra{\rightarrow }

\def\CB{{\cal B}}

\def\CI{{\cal I}}

\def\CC{ {\cal C}}

\def\CH{ {\cal H}}
\def\CP{ {\cal P}}

\def\nAM0{{\nu }_{{\cal A},M_0}}
\def\nAN0{{\nu }_{{\cal A},N_0}}

\def\Kdim{ {\rm Kdim } }

\def\End{ {\rm End }}
\def\Der{ {\rm Der }}
\def\CJ{ {\cal J }}
\def\CR{ {\cal R }}
\def\CP{ {\cal P }}
\def\det{ {\rm det }}
\def\ad{ {\rm ad }}

\def\ga{\mathfrak{a}}

\def\gc{\mathfrak{c}}

\def\gm{\mathfrak{m}}
\def\gp{\mathfrak{p}}

\def\GL{{\rm GL}}
\def\SL{{\rm SL}}
\def\j{{\bf j}}

\def\Hom{{\rm Hom}}

\def\JJ{{\bf J}}

\def\di!{\frac{\der^i}{i!}}
\def\dik!{\frac{\der^k_i}{k!}}

\def\gl{\mathfrak{l}}

\def\id{{\rm id}}

\def\N{\mathbb{N}}

\def\0{\overline{0}}
\def\1{\overline{1}}

\def\Ln1{\L_{n,\overline{1}}}

\def\oa{\overline{a}}

\def\a1{a_{\overline{1}}}

\def\bs{\overline{s}}

\def\S{\Sigma}

\def\grad{{\rm grad}}

\def\CU{{\cal U}}

\def\sign{{\rm sign}}
\def\vn1{\overrightarrow{n-1}}

\def\gl{{\rm gl}}
\def\sl{{\rm sl}}

\def\bu{\overline{u}}
\def\PZ{{\rm PZ}}

\def\mJ{\mathbb{J}}
\def\mI{\mathbb{I}}

\def\ann{{\rm ann}}
\def\lann{{\rm l.ann}}
\def\rann{{\rm r.ann}}

\def\mF{\mathbb{F}}

\def\ind{{\rm ind}}

\def\K1{{\rm K}_1}

\def\hmI1{\widehat{\mI_1}}
\def\tmI1{\widetilde{\mI_1}}
\def\tmJ1{\widetilde{\mJ_1}}
\def\hB1{\widehat{B_1}}
\def\hCB1{\widehat{\CB_1}}

\def\bS{\overline{S}}

\def\Den{{\rm Den}}

\def\Den{{\rm Den}}

\def\br{\overline{r}}

\def\bs{\overline{s}}

\def\ga{\mathfrak{a}}

\def\tor{{\rm tor}}

\def\sl2{\mathfrak{sl}_2}
\def\gl2{\mathfrak{gl}_2}

\def\bCR{\overline{\CR}}

\def\b1{\overline{1}}

\def\PDer{{\rm PDer}}

\def\PIDer{{\rm PIDer}}

\def\Sym{{\rm Sym}}
\def\bCP{\overline{\CP}}
\def\OCP{\Omega_{\CP}}
\def\bCU{\overline{\CU}}

\begin{document}

\author{V. V. \  Bavula 
}

\title{The PBW Theorem and simplicity criteria for the Poisson  enveloping algebra and the  algebra of Poisson differential operators}

\maketitle

\begin{abstract}
For an arbitrary Poisson algebra $\CP$ over an arbitrary field, an (analogue of) the  Poincar\'{e}-Birkhof-Witt Theorem is proven and several presentations/constructions for its Poisson enveloping algebra $\CU (\CP )$ are given. As a result, explicit sets of generators and defining relations are given for  $\CU (\CP )$ and the algebra $P\CD (\CP)$ of  Poisson differential operators on $\CP$. 
 Simplicity criteria for the algebras $\CU (\CP )$ and 
 $P\CD (\CP )$ are given. In the case when the algebra $\CP$ is of essentially finite type, a criterion for the algebra $\CU (\CP )$ to be a domain is presented and a criterion for a natural epimorphism $\CU (\CP )\ra P\CD (\CP )$ to be  an isomorphism is given. The kernel of the epimorphism is described and for  large classes of Poisson algebras an explicit set of generators is given.  Explicit formulae for the Gelfand-Kirillov dimension of the algebras $\CU (\CP )$ and $P\CD (\CP)$ are given. In the case when  the Poisson algebra $\CP$ is a regular domain of essentially finite type
an explicit  simplecticity criterion for $\CP$  is found and  a criterion is presented for the algebra $\CU (\CP )$ to be isomorphic to the algebra $\CD (\CP )$ of differential operators on $\CP$.

$\noindent $



$${\bf Contents}$$
\begin{enumerate}
\item Introduction.
\item Generetors and defining relations of the Poisson enveloping algebra of a Poisson algebra. 
\item The PBW Theorem  for the Poisson  enveloping algebras and the module $\OCP$ of K\"{a}hler differentials of a Poisson algebra $\CP$.  
 \item Criterion for the algebra  $\CU (\CA )$ to be a domain where a Poisson algebra $\CA$  is a domain of essentially finite type. 
 \item Criterion for $\Der_K(\CA ) = \CA \CH_\CA$ where $\CA$ is a regular domain of essentially finite type. 
 \item The kernel of the epimorphism $\CU (\CP ) \ra P\CD (\CP )$ and the defining relations of the algebra $P\CD (\CP )$.
 \item Simplicity criteria for the algebras $\CU (\CP )$ and $ P\CD (\CP )$. 
\end{enumerate}
\end{abstract}


\section{Introduction}\label{P-INTR}

In this paper, module means a left module,
 $K$ is a field, algebra means a $K$-algebra (if it is not stated otherwise) and $K^\times =K\backslash \{ 0\}$.

 An associative  
 commutative algebra $\CP $ is called a {\em Poisson algebra} if it is a Lie algebra $(\CP , \{ \cdot, \cdot \})$ such that $\{ a, xy\}= \{ a, x\}y+x\{ a, y\}$ for all elements $a,x,y\in \CP $.  \\

{\bf The  Poincar\'{e}-Birkhof-Witt Theorem for Poisson algebras.} The (classical) Poincar\'{e}-Birkhof-Witt Theorem states that for each Lie algebra $\CG$ there is a natural isomorphism of   graded algebras,  $${\rm gr}\, U(\CG ) \simeq \Sym (\CG ),$$ where ${\rm gr}\, U(\CG )$ is   the associated graded algebra   of the   universal enveloping algebra  $ U(\CG )$ of the Lie algebra $\CG$ and $\Sym (\CG )$  is the symmetric algebra of $\CG$. For a {\em smooth} Poisson algebra $\CP$,  a
similar result holds \cite[Theorem 3.1]{Rinehart-1963}: 
 $${\rm gr}\, \CU(\CP ) \simeq \Sym_\CP (\O_\CP ),$$ where ${\rm gr}\, \CU(\CP )$ is   the associated graded algebra   of the   Poisson  enveloping algebra  $ \CU(\CP )$ of the Poisson algebra $\CP$ and $\Sym_\CP (\O_\CP )$  is the symmetric algebra of the $\CP$-module  $\O_\CP$ of K\"{a}hler differentials of the associative algebra $\CP$. In fact,  \cite[Theorem 3.1]{Rinehart-1963} holds in sightly more general situation, namely, for the universal enveloping algebra of the Lie-Reinhatr algebra. The pair $(\CP, \O_\CP )$ is an example of a Lie-Reinhart algebra and its universal enveloping algebra $V(\CP, \O_\CP )$ is isomorphic to the Poisson enveloping  algebra (PEA) $\CU (\CP )$, \cite[Section 2, p.197]{Rinehart-1963} (see also \cite{Huebschmann-90}).
 Recently,  it was proven  that the PBW Theorem holds for certain {\em singular} Poisson hypersurfaces,  \cite[Theorem 3.7]{Lambre-Ospel-Vanhaecke-2019}.
 One of the main results of this paper, Theorem \ref{A30Jul19}, states that the  Poincar\'{e}-Birkhof-Witt Theorem holds for all  Poisson algebras (over an arbitrary field). \\

In \cite{Oh-1999}, the Poisson enveloping algebra of a Poisson algebra was introduced as a  universal object in a certain category and an alternative  to Reinhart's proof of its existence was given. For certain classes of Poisson algebras explicit descriptions  of their Poisson enveloping algebras were presented in \cite{Oh-Park-Shin-2002, Umirbaev-2012, Yang-Yao-Ye-2013,Lu-Wang-Zhuang-2015, Lu-Oh-Wang-Yu-2018, Lambre-Ospel-Vanhaecke-2019}. \\



{\bf Simplicity criterion for the algebra $P\CD (\CP )$ of Poisson differential operators on $\CP$.} An ideal $I$ of a Poisson algebra $\CP$ is called a {\em Poisson ideal} if $\{\CP , I\}\subseteq I$. A Poisson algebra $\CP$ is called {\em Poisson simple} if the ideals $0$ and $\CP$ are the only Poisson ideals of the Poisson algebra $\CP$. Let $\Der_K(\CP )$ be the Lie algebra of $K$-derivations of the (associative) algebra $\CP$. For each element $a\in \CP$, the derivation $\pad_a :=\{a, \cdot \}\in \Der_K(\CP)$ is called called the {\em Hamiltonian vector field} associated with the element $a$.  Then 
$\CH_\CP :=\{ \pad_a \, | \, a\in \CP\}$ is a Lie subalgebra of the Lie algebra $\Der_K(\CP )$.  The subalgebra $P\CD (\CP )$ of the $K$-endomorphism algebra $\End_K(\CP )$ which generated by $\CP$ and $\CH_\CP$ is called the {\bf algebra of Poisson differential operators} of the Poisson algebra $\CP$.
 The algebra $P\CD (\CP )$ is a subalgebra of the algebra $\CD (\CP )$ of differential operators on $\CP$. In general, $P\CD (\CP )\neq \CD (\CP )$.
 Theorem \ref{X22Jul19} is a simplicity criterion for the algebra $P\CD (\CP )$ of Poisson differential operators.

\begin{theorem}\label{X22Jul19}
Let $\CP$ be a Poisson algebra over an arbitrary  field $K$.
Then the  following statements are equivalent:
\begin{enumerate}
\item The algebra $P\CD (\CP )$ is a simple algebra.
\item The Poisson algebra $\CP$ is a Poisson simple algebra. 
\end{enumerate}
\end{theorem}

{\bf Simplicity criteria for the Poisson    enveloping algebra $\CU (\CP )$.} There is a natural algebra epimorphism $$\pi_\CP :\CU (\CP) \ra P\CD (\CP ),$$ see  (\ref{UPDP}). The algebra $\CP$ is a $\CD (\CP )$-module and hence $P\CD (\CP )$- and $\CU (\CP)$-module (via $\pi_\CP$).

\begin{theorem}\label{Y22Mar19}
Let $\CP$ be a Poisson algebra over an arbitrary  field $K$. Then the  following statements are equivalent:
\begin{enumerate}
\item The algebra $\CU (\CP )$ is a simple algebra.
\item The algebra $P\CD (\CP )$ is a simple algebra and $\CU (\CP )\simeq P\CD (\CP )$. 
\item  The Poisson algebra $\CP$ is a Poisson simple algebra and $\CP$ is a faithful left $\CU (\CP )$-module. 
\end{enumerate}
If one of the equivalent conditions holds then $\CU (\CP ) \simeq P\CD (\CP )$. 
\end{theorem}

A localization of an affine commutative algebra is called an {\em algebra of essentially finite type}. 
In the case when the Poisson algebra $\CP = \CA$ is an algebra  of essentially finite type over a field of characteristic zero,  Theorem \ref{Y22Mar19} can be strengthen, see Theorem \ref{7Aug19}.  Let us introduce necessary definitions in order to formulate Theorem \ref{7Aug19}:  $P_n=K[x_1, \ldots , x_n]$ is a polynomial algebra over
$K$, $I=(f_1, \ldots , f_m)$ is a  prime but not a maximal ideal of $P_n$, $ \CA = S^{-1}(P_n/I)$ is a domain of essentially finite type and 
 $Q=Q(\CA )$ is its field of fractions,  $r= r\Big(\frac{\der  f_i}{\der x_j}\Big)$ is the rank (over $Q$) of the Jacobian matrix $\Big(\frac{\der  f_i}{\der x_j}\Big)$ of $\CA$ and $d=d_\CA = r(\CC_\CA )$ is  the rank (over $Q$) of the $n\times n$ matrix $\CC_\CA =(\{ x_i , x_j\} )\in M_n(\CA )$, $\GK$ stands for the Gelfand-Kirillov dimension.

\begin{theorem}\label{7Aug19}
Let  a Poisson algebra $\CA$ be an algebra of essentially finite type over the field $K$ of characteristic zero.
 Then the  following statements are equivalent:
\begin{enumerate}
\item The algebra $\CU (\CA )$ is a simple algebra.
\item The 
algebra $P\CD (\CA )$ is a simple algebra and one of the equivalent conditions of Theorem \ref{3Aug19} holds. 
\item  The algebra $\CA$ is Poisson simple 
  and one of the equivalent conditions of Theorem \ref{3Aug19} holds. 
\end{enumerate}
If one of the equivalent conditions holds then the algebra $\CA = S^{-1}(P_n/I)$ is a regular, Poisson simple  domain of essentially finite type over the  field $K$ of characteristic zero, the algebra epimorphism 
$\pi_\CA :  \CU (\CA) \ra  P\CD (\CA )$ is an isomorphism (see (\ref{UPDP})), $d=n-r$ where  $d= r(\CC_\CA )$ and $r= r\Big(\frac{\der  f_i}{\der x_j}\Big)$, and the algebra $\CU (\CA )$ is a simple Noetherian domain with  $$\GK \, \CU (\CA ) = \GK \, P\CD (\CA) = \GK \, {\rm gr}\, \CU (\CA ) = \GK \, \Sym_\CA (\O_\CA )= 2\GK (\CA ) = 2(n-r).$$ 
\end{theorem}

The proofs of Theorem \ref{X22Jul19}, Theorem \ref{Y22Mar19} and Theorem \ref{7Aug19} are given in Section \ref{SIMUPDP}. \\

{\bf Generators and defining relations for the Poisson enveloping algebra $\CU (\CP )$.} In Section \ref{GRREL}, for each Poisson algebra $\CP$  explicit sets of generators and defining relations  for its Poisson enveloping algebra $\CU (\CP )$ are given  (Theorem \ref{23Jun19}). Several (expected) results about the Poisson enveloping algebras are proven that are used later in the paper.  It is proven that localizations commute with the operation of taking the Poisson    enveloping algebra (Theorem \ref{30Jun19}).  The symmetric algebra $S(\CG ) =\Sym (\CG )$ of a Lie algebra $\CG$ admits the canonical Poisson structure that is determined by the Lie structure on $\CG$. Proposition \ref{C23Jun19} is an explcit description of the algebra $\CU (S (\CG ))$. It is shown that $\CU (S (\CG_1\times \CG_2))\simeq \CU (S (\CG_1) \otimes  \CU (S( \CG_2))$ (Corollary \ref{a23Jun19}) where $\CG_1\times \CG_2$ is a direct product of Lie algebras. The structure of the PEA  of the Poisson symmetric algebra of a semi-direct product of Lie algebras is described (Corollary \ref{bC23Jun19}). \\

 {\bf The Gelfand-Kirillov dimension of the algebras $\CU (\CA )$,    ${\rm gr}\, \CU (\CA )$ and $ \Sym_\CA (\O_\CA )$ where $\CA$ is a domain  of essentially finite type.} 
 
\begin{theorem}\label{AAA29Jul19}
Let a Poisson algebra $\CA = S^{-1} (P_n/I)$ be a  domain of essentially finite type over a perfect field $K$ where  $I=(f_1, \ldots , f_m)$ is a prime   ideal of $P_n$ and $r=r\Big(\frac{\der f_i}{\der x_j}\Big)$ is the rank of the Jacobian matrix $\Big(\frac{\der f_i}{\der x_j}\Big)$ over the field of fractions of the domain $P_n/I$.  Then the algebra $\CU (\CA )$ is a Noetherian algebra with $$\GK \, \CU (\CA ) = \GK \, {\rm gr}\, \CU (\CA ) = \GK \, \Sym_\CA (\O_\CA ) = 2\GK (\CA ) = 2(n-r).$$
\end{theorem}

{\bf The Gelfand-Kirillov dimension  of 
the algebra $P\CD (\CA )$ of Poisson differential operators on $\CA $.} Proposition \ref{BB26Jul19} gives the exact figure for the Gelfand-Kirillov dimension  of 
the algebra $P\CD (\CA )$.
\begin{proposition}\label{BB26Jul19}
Let a Poisson algebra $\CA$ be a domain of essentially finite type over the field $K$ of characteristic zero, $r$ is the rank of Jacobian matrix of $\CA$  and $d=r(\CC_\CA )$. Then 
 $$\GK (P\CD (\CA )) = \GK (\CA ) +d= n-r+d .$$
\end{proposition}
The proof of Proposition \ref{BB26Jul19} is given in Section \ref{KERCUCP}, see Proposition \ref{B26Jul19}.\\

{\bf Criterion for the algebra $\CU (\CA )$ to be a domain.} Theorem \ref{XBA29Jul19}  (see Theorem \ref{A29Jul19}.(3))  and Theorem  \ref{A31Jul19} are criteria for the algebra $\CU (\CA )$ to be a domain where the Poisson algebra $\CA$ is a domain of essentially finite type. The first one is given in terms of the Jacobian matrix and the Jacobian ideal of $\CA$ and the second one -- in terms of  the grades of the Jacobian ideals and prime ideals (but for certain class of Poisson algebras of essentially finite type).

 For $\i =(i_1, \ldots , i_r)$ such that $1\leq
i_1<\cdots <i_r\leq m$ and $\j =(j_1, \ldots , j_r)$ such that
$1\leq j_1<\cdots <j_r\leq n$, 

$$\D  (\i , \j ):=\det \bigg(\frac{\der f_{i_\nu}}{\der x_{ j_\mu}}\bigg), \;\; \nu , \mu
=1, \ldots, r,$$  denotes the corresponding minor of the Jacobian matrix of the algebra $\CA$  where  $r$ is the rank of the Jacobian matrix of $\CA$. The $r$-tuple $\i$ (resp., $\j $) is called {\bf
non-singular} if $\D (\i , \j')\neq 0$ (resp., $\D (\i', \j )\neq
0$) for some $\j'$ (resp., $\i'$). We denote by $\mI_r$ (resp.,
$\mJ_r$) the set of all the non-singular $r$-tuples $\i$ (resp.,
$\j $). By  \cite[Lemma 2.1]{gendifreg},  $\D (\i , \j )\neq 0$ iff $\i\in \mI_r$ and $\j\in \mJ_r$. The {\bf Jacobian ideal} $\ga_r$ of the algebra $\CA$ is an ideal of $\CA$ that is generated by all the minors $\D(\i,\j )$ of the Jacobian matrix of $\CA$.

\begin{theorem}\label{XBA29Jul19}
Let a Poisson algebra $\CA = S^{-1} (P_n/I)$ be a domain of essentially finite type over a perfect field $K$ where  $I=(f_1, \ldots , f_m)$ is a prime but not maximal  ideal of $P_n$, $r=r\Big(\frac{\der f_i}{\der x_j}\Big)$ is the rank of the Jacobian matrix $\Big(\frac{\der f_i}{\der x_j}\Big)$  and $\ga_r$ is the Jacobian ideal of $\CA$. Then 
 the following statements are equivalent:
\begin{enumerate}
\item The algebra $\CU (\CA)$ is a domain.
\item The algebra ${\rm gr} \, \CU (\CA)$ is a domain.
\item The algebra $\Sym_{\CA} (\O_\CA )$ is a domain.
\item The elements $\{ \D (\i , \j )\, | \, \i \in \mI_r, \j\in \mJ_r\}$ are regular in $\CU (\CA )$. 
\item The element $\D (\i , \j )$ is a regular element of the algebra $\CU ( \CA )$ for some $\i \in \mI_r$ and $\j \in \mJ_r$. 
\item $\lann_{\CU (A)}(\ga_r) = \rann_{\CU (A)}(\ga_r)=0$. 
\end{enumerate}
\end{theorem}

In proving the theorem below a result of Huneke is used, see Theorem \ref{HunekeT1.1}. For an $\CA$-module $M$, we denote by $v(M,\CA)$ its minimal number of generators.

\begin{theorem}\label{A31Jul19}
Let a Poisson algebra $\CA = S^{-1}(P_n/I)$ be a universally catenarian domain of essentially finite type over a perfect field $K$ satisfying Serre's condition $S_m$ and $I=(f_1, \ldots , f_m)$. The following statements are equivalent: 
\begin{enumerate}
\item The algebra $\CU (\CA )$ is a domain.
\item The algebra ${\rm gr} \, \CU (\CA)$ is a domain.
\item The algebra $\Sym_{\CA} (\O_\CA )$ is a domain.
\item ${\rm grade} (\ga_t)\geq m+2-t$ for $1\leq t \leq m$ where $\ga_t$ is the ideal of $\CA$ generated by $t\times t$ minors of the Jacobian matrix $\Big(\frac{\der f_i}{\der x_j}\Big)$. 
\item $v(\O_\gp , \CA_\gp ) \leq n-m+{\rm grade} (\gp ) -1$ for all nonzero primes $\gp$ of $\CA$.
\end{enumerate}
If one of the equivalent conditions holds then the algebra ${\rm gr} \, \CU (\CA ) \simeq \Sym_{\CA } (\O )$ is a complete intersection in the polynomial algebra $\CA [ \d_1, \ldots , \d_n]$. In particular, if the algebra $\CA$ is Cohen-Macaulay (resp., Gorenstein) then so is the algebra ${\rm gr}\, \CU (\CA ) \simeq \Sym_\CA (\O)$. 
\end{theorem}

 The proofs of Theorem \ref{XBA29Jul19}  (see Theorem \ref{A29Jul19}.(3))  and Theorem  \ref{A31Jul19} are given in Section \ref{CUADOM}. \\
 
 {\bf Simplecticity criterion for the Poisson algebra $\CA$  of essentially finite type, i.e. $\Der_K(\CA ) = \CA \CH_\CA$.}  For each Poisson algebra $\CP$, $\Der_K(\CP)\supseteq \CP \CH_\CP$. A Poisson algebra $\CP$ which is a regular (affine) domain  is called a {\em synmplectic algebra} if $\Der_K(\CP)= \CP \CH_\CP$. Theorem \ref{22Jul19} is a criterion for $\Der_K(\CA ) = \CA \CH_\CA$ where $\CA$ is a regular domain of essentially finite type. Let $\gc_{\CA , d}$ be the ideal of the algebra $\CA$ which is generated by all the $d\times d$ minors of the matrix $\CC_\CA $ where $d:=r(\CC_\CA )$ is its rank.

\begin{theorem}\label{22Jul19}
Let a Poisson algebra $\CA$ be  a regular domain of essentially finite type over a field $K$ of characteristic zero and $d:=r(\CC_\CA )$. Then the following statements are equivalent:
\begin{enumerate}
\item $\Der_K(\CA)=\CA \CH_\CA$. 
\item $d=n-r$ and $\gc_{\CA , d}=\CA$.
\item For each $\i\in \mI_r$ and $ \j\in \mJ_r$, $\D (\i , \j )^{n-r}\in \gm_\j$ where $\gm_\j$ is the ideal of $\CA$ generated by all the $(n-r)\times (n-r)$ minors of the $n\times (n-r)$ matrix $\CC_{\CA, \j}$ (see Proposition \ref{A23Jul19}).
\item For each $\i\in \mI_r$ and $ \j\in \mJ_r$, $\D (\i , \j )^k\in \gc_{\CA , n-r}$ for some $k\geq 1$. 
\end{enumerate}
\end{theorem}

 Lemma \ref{a24Jul19} and  Corollary \ref{b24Jul19} are regularity and symplecticity criteria for certain generalized Weyl Poisson algebras.\\
 
 {\bf Criteria for $\ker (\pi_\CA )=0$.}
 Recall that there  is a natural algebra epimorphism $\pi_\CP :\CU (\CP) \ra P\CD (\CP )$, see  (\ref{UPDP}).  In the case when the Poisson algebra $\CP =\CA$ is a {\em regular} domain of essentially finite type, Theorem \ref{3Aug19} is an efficient explicit criterion for $\ker (\pi_\CA )=0$, i.e. for the epimorphism $\pi_\CA : \CU (\CA ) \ra P\CD (\CA )$ to be an isomorphism.

Let $\kappa_\CA$ be an ideal of the algebra $\CU (\CA )$ which is generated by a finite set of explicit  elements  $\d_{\i , i_\nu; \j }\in \O_\CA$ where $\i , \j \in \mI_\CA (d)$, see (\ref{dinnj}). Then  $\kappa_\CA\subseteq \ker (\pi_\CA )$, see (\ref{kAkpA}). Since $\Der_K(\CA ) \simeq \Hom_\CA (\O_\CA , \CA )$, there is a {\em pairing} of left $\CA$-modules  (which is  an $\CA$-bilinear map, see (\ref{PairDer})):
$$
\Der_K(\CA ) \times \O_\CA \ra \CA , \;\; (\der , \o ) \mapsto (\der, \o ) := \der (\o ).
$$

\begin{theorem}\label{3Aug19}
Let a Poisson algebra $\CA= S^{-1}(P_n/I)$ be a regular domain of essentially finite type over the field $K$ of characteristic zero, $d=r(\CC_\CA )$ and $r=r\Big( \frac{\der f_i}{\der x_j}\Big)$. Then the following statements are equivalent (the derivations $\der_{\i ; \j , j_\nu}$ of $\CA$ are defined in Theorem \ref{9bFeb05}):
\begin{enumerate}
\item $\ker (\pi_\CA )=0$ $(\Leftrightarrow \pi_\CA : \CU (\CA ) \simeq P\CD (\CA ))$.  
\item $\kappa_\CA =0$. 
\item $d=n-r$ and $(\der_{\i ; \j , j_\nu}, \d_{\i' ; \j' , j_\mu'})=0$ for all elements $\i \in \mI_r$, $\j\in \mJ_r$, $\i'\in \mI_\CA (d)$, $\j'\in \mJ_\CA (d)$, $\nu = r+1, \ldots , n$ and $\mu = d+1, \ldots , n$ where for $\j = (j_1, \ldots , j_r)$ and $\i' = (i_1', \ldots , i'_d)$, $\{ j_{r+1}, \ldots , j_n\} =\{ 1,\ldots , n\} \backslash \{ j_1, \ldots , j_r\}$ and $\{ i'_{d+1}, \ldots , i'_n\} =\{ 1, \ldots , n\} \backslash \{ i_1', \ldots , i_d'\}$. 
\end{enumerate}
\end{theorem}

Theorem \ref{A5Aug19} is a criterion for $\ker (\pi_\CA )=0$ in the general situation, i.e. without the assumption of regularity  in the previous theorem.

\begin{theorem}\label{A5Aug19}
Let a Poisson algebra $\CA= S^{-1}(P_n/I)$ be a  domain of essentially finite type over a field of characteristic zero,  $r=r\Big( \frac{\der f_i}{\der x_j}\Big)$ and $d=r(\CC_\CA )$. Then the following statements are equivalent: 
\begin{enumerate}
\item $\ker (\pi_\CA )=0$ $(\Leftrightarrow \pi_\CA : \CU (\CA ) \simeq P\CD (\CA ))$.  
\item $\kappa_\CA =0$ and $\CU (\CA )$ is a domain (or any of the equivalent conditions of Theorem \ref{XBA29Jul19} holds). 
\item $d=n-r$, $\CU (\CA )$ is a domain (or any of the equivalent conditions of Theorem \ref{XBA29Jul19}  holds)  and $(\Der_K(\CA ), \d_{\i' ; \j' , j_\mu'})=0$ for all elements $\i'\in \mI_\CA (d)$, $\j'\in \mJ_\CA (d)$ and $\mu = d+1, \ldots , n$ where for $\i' = (i_1', \ldots , i'_d)$, $\{ j_{r+1}, \ldots , j_n\} =\{ 1,\ldots , n\} \backslash \{ j_1, \ldots , j_r\}$ and $\{ i'_{d+1}, \ldots , i'_n\} =\{ 1, \ldots , n\} \backslash \{ i_1', \ldots , i_d'\}$. 
\end{enumerate}
\end{theorem}

{\bf Criterion for the homomorphism $\pi_\CA : \CU (\CA ) \ra \CD (\CA )$ to be an isomorphism.} For a regular domain of essentially  finite type $\CA$, Theorem \ref{5Aug19} is a criterion for the  homomorphism $\pi_\CA : \CU (\CA ) \ra \CD (\CA )$ to be an isomorphism, i.e. the algebra $\CU (\CA )$ is the algebra of differential operators $ \CD (\CA )$ on $\CA$.

\begin{theorem}\label{5Aug19}
Let a Poisson algebra $\CA= S^{-1}(P_n/I)$ be a regular domain of essentially finite type over the field of characteristic zero, $r=r\Big( \frac{\der f_i}{\der x_j}\Big)$ and $d=r(\CC_\CA )$. Then the following statements are equivalent:
\begin{enumerate}
\item  The homomorphism $\pi_\CA : \CU (\CA ) \ra \CD (\CA )$ is an isomorphism.  
\item $\kappa_\CA =0$ and $\gc_{\CA , d}=\CA$.
\end{enumerate}
If one of the equivalent conditions holds then $\ker (\pi_\CA )=0$ and $\CU (\CA ) =P\CD (\CA ) = \CD (\CA )$. 
\end{theorem}

Every Poisson enveloping algebra is isomorphic to its opposite algebra (Theorem \ref{C29Jul19}.(2)). So, left and right properties of the Poisson enveloping algebra are identical.


\section{Generetors and defining relations of the Poisson enveloping algebra of a Poisson algebra}\label{GRREL}

The aim of this section is to prove Theorem \ref{23Jun19} that for each Poisson algebra $\CP$ given explicit sets of generators and defining relations  for the algebra $\CU (\CP )$. As a result, for each Poisson factor algebra $\CP'$ of a Poisson algebra $\CP$, explicit  sets of generators and defining relations are given for the algebra $\CU (\CP')$ (Theorem \ref{A23Jun19}). It is shown that $\CU (\CP_1\t \CP_2) \simeq \CU (\CP_1) \t \CU (\CP_2 )$ (Proposition \ref{B23Jun19}) where $\CP_1\t \CP_2$ is a tensor product of Poisson algebras. Proposition \ref{D23Jun19} shows that every endomorphism/automorphism of a Poisson algebra can be lifted to an endomorphism/automorphism of its Poisson    enveloping algebra. It is proven that localizations commute with the operation of taking the Poisson    enveloping algebra (Theorem \ref{30Jun19}). The symmetric algebra $S(\CG ) =\Sym (\CG )$ of a Lie algebra $\CG$ admits the canonical Poisson structure that is determined by the Lie structure on $\CG$. Proposition \ref{C23Jun19} is an explcit description of the algebra $\CU (S (\CG ))$. It is shown that $\CU (S (\CG_1\times \CG_2))\simeq \CU (S (\CG_1) \otimes  \CU (S( \CG_2))$ (Corollary \ref{a23Jun19}) where $\CG_1\times \CG_2$ is a direct product of Lie algebras. The structure of the PEA  of the Poisson symmetric algebra of a semi-direct product of Lie algebras is described in Corollary \ref{bC23Jun19}. A criterion for the PEA  of a Poisson algebra to be a {\em commutative} algebra is given in Corollary \ref{c23Jun19}. Examples of the PEAs are considered. At the beginning of the section we recall some definitions and results about Poisson algebras and their modules. \\

{\bf Poisson algebras.} An associative  commutative algebra $D$ is called a {\em Poisson algebra} if it is a Lie algebra $(D, \{ \cdot, \cdot \})$ such that $\{ a, xy\}= \{ a, x\}y+x\{ a, y\}$ for all elements $a,x,y\in D$.

 For a $K$-algebra  $D$,  let $\Der_K(D)$ be the set of its $K$-derivations. If, in addition, $(D, \{ \cdot , \cdot \} )$ is a Poisson algebra then 
 $$\PDer_K(D):=\{ \d \in \Der_K(D)\, | \, \d (\{ a,b\})= \{ \d (a),b\}+\{ a,\d (b)\} \;\; {\rm  for\; all}\;\; a,b\in D\}$$
 is the set of derivations of the Poisson algebra $D$. The vector space $\Der_K(D)$ is a Lie algebra, where $[\d , \der ]:= \d \der -\der \d$. The set of {\em inner derivations of the Poisson algebra} $D$,  
 $$\PIDer_K(D):= \{ \pad_a\, | \, a\in D\}\;\; {\rm  (where}\;\;  \pad_a (b) := \{ a,b\} )$$ 
 is an ideal of the Lie algebra $\PDer_K(D)$ (since $[\d ,\pad_a]=\pad_{\d (a)}$ for all $\d\in \PDer_K(D)$ and $a\in D$). By the very definition, the Poisson algebra $D$ is a Lie algebra with respect to the bracket $\{ \cdot , \cdot \}$. The map $ D\ra \PIDer_K(D)$, $a\mapsto \pad_a$, is an epimorphism of Lie algebras with kernel 
 $$\PZ (D):= \{ a\in D\, | \, \{ a, D\} =0\}$$ which is called the {\em centre} of the Poison algebra (or the {\em Poisson centre} of $D$). The Poisson  centre $\PZ(D)$  is invariant under the action of $\PDer_K(D)$: Let $z\in \PZ(D)$, $d\in D$ and $\der \in \PDer_K(D)$; then applying the derivation  $\der$  to the equality $\{ z,d\} =0$ we obtain the equality  $\{ \der (z) ,d\}=0 $, i.e. $\der (z) \in \PZ(D)$.\\

{\bf The dual Poisson algebra of a Poisson algebra.} Given an associative algebra $\CP$. Then its {\em dual (associative) algebra} $\CP^{op}$ coincides with $\CP$ as a vector space but the multiplication is given by the rule $a*b:= ba$. Every left $\CP$-module is a right  $\CP^{op}$-module, and vice versa.

Similarly, given a Poisson algebra $(\CP , \{ \cdot , \cdot \} )$. Its  dual associative algebra $\CP^{op}$ is a Poisson algebra $(\CP^{op} , \{ \cdot , \cdot \}^{op} )$, which is called {\em the dual Poisson algebra} of $\CP$, where $\{ a , b \}^{op}:= -\{ a , b \}$ for all $a,b\in \CP^{op}$. \\

{\bf  The Poisson structure constants matrix $C_\CP$ and the ideal $\gc_\CP$ of a  Poisson algebra $\CP$.} Let $\CP$ be a Poisson algebra and $\{ x_i\}_{i\in I}$ be  a set of algebra generators of $\CP$. 
The Poisson structure on an associative algebra $\CP$ is uniquely determined by the {\em Poisson structure constants}  $c_{ij}:= \{ x_i, x_j\}$ where $i,j\in I$.   Let $n={\rm card} (I)$ be the cardinality of the set $I$, the case $n=\infty$ is possible. The $n\times n$ matrix  
\begin{equation}\label{CPMatCoef}
C_\CP := (c_{ij})
\end{equation}
is called the {\em Poisson structure constants matrix} of the Poisson algebra $\CP$ and the ideal $\gc_\CP$ of $\CP$, which is generated by all the structure constants $c_{ij}$, is called the {\em Poisson structure constants ideal} of the Poisson algebra $\CP$.

For all elements $f,g\in \CP$, 
$$ \{ f, g\} = f'C_\CP g'^t$$
where $f':={\rm grad}(f) := (\frac{\der f}{\der x_1}, \ldots , \frac{\der f}{\der x_n})$ is the {\em gradient} of $f$ (a $1\times n$ matrix with coefficients in $\CP$) and $t$ is the transposition of  matrix. 

If $\{ x_s'\}_{j\in I'}$ is another set of algebra generators for $\CP$ and $C_\CP':=(c_{st}')$ 
 is the Poisson matrix constants that correspond to the second set of generators (where $c_{st}':=\{ x_s', x_t'\}$). Then ($t$ is the trasposition of  matrix) 
\begin{equation}\label{CPMatCoef1}
C_\CP' := \CJ^t C_\CP \CJ\;\; {\rm where}\;\; \CJ :=  \CJ (x', x):=\frac{\der x'}{\der x}
\end{equation}
the {\em Jacobian matrix} of the change of the variables from $x=\{ x_i\}_{i\in I} $ to $x'=\{ x_s'\}_{s\in I'} $, i.e. the $s^{\rm th}$ row of the
 ${\rm card} (I')\times {\rm card} (I)$ matrix $C_\CP'$ is the  gradient
$\grad (x'_s):= (\frac{\der x_s'}{\der x_i})$ of the function $x_s'=x_s'(\ldots , x_i, \ldots )$ where $i\in I$. 

\begin{lemma}\label{a22Mar19}
Let $\CP$ be a Poisson algebra. Then the  Poisson structure constants ideal $\gc_\CP$  is a Poisson ideal   of $\CP$ (i.e. $\{ \CP , \gc_\CP \} \subseteq \gc_\CP$) which  does not depend on the choice of algebra generators of $\CP$.
\end{lemma}

{\it Proof}. The lemma follows at once from (\ref{CPMatCoef1}). $\Box $\\

{\bf A module of a Poisson algebra.} Let a commutative associative algebra $(\CP , \{ \cdot, \cdot \} )$ be a  Poisson algebra and $M$ be a left $\CP$-module ($\CP \times M\ra M$, $(a, m) \mapsto am$ ). The  left $\CP$-module $M$ over the associative algebra $\CP$ is called a {\em left module over the Poisson algebra} or a {\em Poisson  left $\CP$-module} if there is a bilinear map
$$ \CP \times M\ra M, \;\; (a, m) \mapsto \d_am$$
which is called a {\em Poisson action} of $\CP$ on $M$ such that for all elements $a,b\in \CP$ and $m\in M$,\\

(PM1) $\d_{\{ a,b\} }= [\d_a, \d_b ]$, \\

(PM2)  $[\d_a, b]= \{ a,b\}$, and \\

(PM3) $\d_{a b}= a\d_b+b\d_a$.\\

Every left Poisson $\CP$-module $M$ determines the homomorphism of {\em associative}  algebras,
\begin{equation}\label{PEndM}
\CP \ra \End_K(M), \;\; a\mapsto a_M: M\ra M, \;\; m\mapsto am 
\end{equation}
 and the homomorphism of {\em Lie} algebras,
\begin{equation}\label{PEndM1}
\CP \ra \End_K(M), \;\; a\mapsto \d_a: M\ra M, \;\; m\mapsto \d_am  
\end{equation}
such that 
\begin{equation}\label{PEndM2}
[\d_a , b_M]=\{ a, b\}_M\;\; {\rm for \; all}\;\; a, b\in \CP ,  
\end{equation}

\begin{equation}\label{PEndM3}
\d_{ab} = a_M\d_b+b_M\d_a\;\; {\rm for \; all}\;\; a, b\in \CP ,  
\end{equation}

and vice versa. Indeed, (\ref{PEndM}) determines a $\CP$-module structure on $M$,    (\ref{PEndM1}) determines a Lie $\CP$-module structure on $M$, and   (\ref{PEndM2}) and    (\ref{PEndM3}) are equivalent to the properties  (PM2) and (PM3),  respectively. So, a Poisson $\CP$-module is a $\CP$-module over the associative algebra $\CP$ and the Lie algebra $\CP$ and both module structures are related by  (\ref{PEndM2}) and  (\ref{PEndM3}). \\

{\it Example.}  Every Poisson algebra $\CP$ is a left Poisson $\CP$-module
where for all $a\in \CP$, $a_\CP : \CP \ra \CP$, $  b\mapsto ab$ and $ \d_a = \{ a, \cdot \}:  \CP \ra \CP$, $b\mapsto \{ a, b\}$. \\

{\bf Right  modules of a Poisson algebra.} Given a Poisson algebra $(\CP , \{ \cdot , \cdot \} )$. A {\em right Poisson module} over $\CP$ is, by definition, a {\em left}  Poisson  module over the dual Poisson algebra 
 $(\CP^{op} , \{ \cdot , \cdot \}^{op} )$. \\
 
{\it Example.}  Every Poisson algebra $\CP$ is a right Poisson $\CP$-module  where for all $a\in \CP$, ${}_\CP a : \CP \ra \CP$, $  b\mapsto ba$ and $ \d_a' = \{  \cdot , a\}:  \CP \ra \CP$, $b\mapsto \{  b, a\}$. \\


{\bf The ring $\CD (A)$ of differential operators on an  algebra $A$}.  Let us recall the definition of the ring of differential operators $\CD (A)$ on a commutative algebra $A$ over a field $K$.   The ring of ($K$-linear)
{\bf differential operators} $\CD (A)$ on $A$ is defined as a
union of $A$-modules  $\CD (A)=\bigcup_{i=0}^\infty \,\CD_i (A)$
where $\CD_0 (A)=\End_R(A)\simeq A$, $(x\mapsto ax)\lra a$,
$$ \CD_i (A)=\{ u\in \End_K(A)\, |\, [a,u]:=au-ua\in \CD_{i-1} (A)\; {\rm for \; all \; }\; a\in A\}.$$
 The set of $A$-modules $\{ \CD_i (A)\}$ is the {\bf order filtration} for
the algebra $\CD (A)$:
$$\CD_0(A)\subseteq   \CD_1 (A)\subseteq \cdots \subseteq
\CD_i (A)\subseteq \cdots\;\; {\rm and}\;\; \CD_i (A)\CD_j
(A)\subseteq \CD_{i+j} (A) \;\; {\rm for \; all} \;\; i,j\geq 0.$$

The subalgebra $\D (A)$ of $\CD (A)$ generated by $A\equiv
\End_R(A)$ and the set ${\rm Der}_K (A)$ of all $K$-derivations of
$B$ is called the {\bf derivation ring} of $A$.

Suppose that $A$ is a  regular affine  domain of Krull dimension
$n<\infty $ over a field $K$ of characteristic zero. In geometric terms, $A$ is the coordinate ring $\OO
(X)$ of a smooth irreducible  affine algebraic variety $X$ of
dimension $n$. Then
\begin{itemize}
\item ${\rm Der}_K (A)$ {\em is a finitely generated projective}
$A$-{\em module of rank} $n$, \item  $\CD (A)=\Delta (A) $, \item
$\CD (A)$ {\em is a simple (left and right) Noetherian domain of
Gelfand-Kirillov dimension}  $\GK \, \CD (A)=2n$ ($n=\GK (A)=\Kdim
(A))$.
\end{itemize}

For the proofs of the statements above the reader is referred to
\cite[Chapter 15]{MR}.
 So, the domain $\CD (A)$ is a simple finitely generated infinite dimensional Noetherian algebra, 
\cite[Chapter 15]{MR}.\\

{\bf The algebra $P\CD (\CP )$ of Poisson differential operators.}
 Let $\CP$ be a Poisson algebra. Then 
$\CH_\CP :=\PIDer_k(\CP)=\{ \pad_a :=\{ a, \cdot \}\, | \, a\in \CP\}$ is a Lie subalgebra of the Lie algebra $\Der_K(\CP )$ (since $[\pad_a, \pad_b ] = \pad_{\{ a, b\} }$) and a $\PZ (\CP )$-submodule of $\Der_K(\CP )$.  \\

{\it Definition.} Let $\CP$ be a Poisson algebra (not necessarily commutative). The subalgebra $P\CD (\CP )$ of $\End_K(\CP )$ which generated by $L_\CP$ and $\CH_\CP$ is called the {\bf algebra of Poisson differential operators}.  Clearly, $P\CD (\CP )\subseteq \D (\CP ) \subseteq \CD (\CP)$. \\

{\bf Semidirect products of algebras.} Let $D$ be a $K$-algebra, $\CG$ be a Lie algebra and $U(\CG )$ be the    enveloping algebra of the Lie algebra $\CG$ and $\d : \CG \ra \Der_K(D)$, $a\mapsto \d_a$ be a Lie algebra homomorphism $(\d_{[a,b]}= [ \d_a, \d_b]$ for all $a,b\in \CG)$. Let $D\rtimes_\d U(\CG )$  be the semidirect product of $D$ and $U(\CG )$. It is an associative algebra that is generated by the algebras $D$ and $U(\CG )$ subject to the defininf relations: For all elements $d\in D$ and $ g\in \CG$, $gd=dg+\d_g(d)$. Let $\{ x_i\}_{i\in I}$ be a $K$-basis of the Lie algebra $\CG$. Then 
$$ D\rtimes_\d U(\CG )=\bigoplus_{\alpha \in \N^{(I)}}Dx^\alpha =\bigoplus_{\alpha \in \N^{(I)}}x^\alpha D$$ 
is a free left and right $D$-module where $\N^{(I)}$ is a direct sum of $I$ copies of the set $\N$, $x^\alpha = \prod_{i\in I} x_i^{\alpha_i}$ (in the product the order of multiples is arbitrary and all but finitely many $\alpha_i$ are equal to zero). \\

{\it Example.} Let $\CP$ be a Poisson algebra and $U(\CP )$ be its universal enveloping algebra as a Lie algebra. Then $\CP \rtimes_\pad  U (\CP )$ is a semidirect  product of $\CP$ and $U (\CP )$ where $\pad : \CP \ra \PDer_K(\CP )$, $a\mapsto \pad_a:= \{ a,\cdot \}$. \\

Given another semidorect product $D'\rtimes_{\d'} U(\CG' )$, a homomorphism $\v : D\ra D'$ and a Lie homomorphism $\psi : \CG \ra \CG'$ such that $[\psi (g) , \v (d)]= \d_{\psi (g)}'(\v (d))$ for all elements $d\in D$ and $ g\in \CG$.  Then the map 
\begin{equation}\label{DUG}
(\v , \psi ) : D\rtimes_\d U(\CG )\ra  D'\rtimes_{\d'} U(\CG' ), \;\; d\mapsto \v (d), \;\; g\mapsto \psi (g)
\end{equation}
is a homomorphism. \\

{\bf Two presentations,  $\CP = S^{-1}P_\L /I\simeq \bS^{-1}(P_\L/I')$.} Let $P_\L =K[x_i]_{i\in \L}  $ be a polynomial algebra  where $\{ x_i\}_{i\in \L}$ is a set of variables (no restriction on the cardinality of $\L$), $S$ be a multiplicative subset of $P_\l$, $S^{-1}P_\L$ is a localization of $P_\L$ at $S$ and $\CP = S^{-1}P_\L /I$ be a factor algebra of  $S^{-1}P_\L$  modulo an ideal $I=(f_s)_{s\in \G}$ where $\{ f_s\}_{s\in \G}$ is a set of its generators. The algebra $\CP$ can also be written in the form $\bS^{-1}(P_\L/I')$ where $I'$ is an ideal of the polynomial algebra $P_\L$ and $\bS$ is a multiplicative subset of regular elements of the factor algebra $P_\L/I'$ (eg, $\bS$ is the image of the set $S$ under the epimorphism $P_\L \ra P_\L / I'$).  

Conversely, given an algebra $\bS^{-1}(P_\L/I')$ where $I'$ is an ideal of the polynomial algebra $P_\L$ and $\bS$ is a multiplicative subset of regular elements of the factor algebra $P_\L/I'$. Then the algebra $\bS^{-1}(P_\L/I')$ can be written as $ S^{-1}P_\L /I$ for some multiplicative subset $S$ of $P_\L$ and some ideal $I$ of $S^{-1}P_\L$ (eg, $S=\v^{-1}(\bS)$ and $I=S^{-1}\ker (\v )$ where $\v : P_\L \ra \bS^{-1}(P_\L/I')$). So,

\begin{equation}\label{CP-2-pres}
\CP = S^{-1}P_\L /I\simeq \bS^{-1}(P_\L/I')
\end{equation}
where $I'$ is an ideal of $P_\L$ and $\bS$ is a multiplicative subset of regular elements of $P_\L/I'$.

Why we stressed this seemingly obvious fact? In theoretical arguments,  the second presentation is slightly more preferable but in dealing with examples and  computations, the first  one is (as the number of generators of an ideal under localizations is dropped, as a rule). \\

{\bf Generators and defining relations of the algebra $\CU (\CP )$.}

{\it Definition.} Let $(\CP , \{ \cdot , \cdot \})$ be a Poisson algebra over $K$. A triple $(\CU , \alpha, \beta )$ is called a {\bf Poisson  enveloping algebra} (PEA, for short) of the Poisson algebra $\CP$ if $\CU$ is an (associative) algebra, $\alpha : \CP \ra \CU$ is an algebra homomorphism,  $\beta : \CP \ra \CU$ is a Lie algebra homomorphism such that  for all elements $a,b\in \CP$, 
$$ [\beta (a) , \alpha (b) ] = \alpha (\{ a,b\})\;\; {\rm and}\;\; \beta (ab) = \beta (a)\alpha (b)+\alpha (a)  \beta (b), $$
if $(\CU', \alpha' , \beta')$ is another triple as above then there is a unique algebra homomorphism $f: \CU \ra \CU'$  such that $ \alpha'= f\alpha $  and $\beta' = f\beta$. \\

 For a Poisson algebra $\CP$ which is defined by generators and defining relations (as an associative algebra), Theorem \ref{23Jun19} gives explicit sets of generators and defining relations for the Poisson    enveloping algebra $\CU (\CP )$. 

\begin{theorem}\label{23Jun19}
Let $\CP$ be a Poisson algebra, $U (\CP )$ be its universal enveloping algebra (as a Lie algebra)  and  $\CU (\CP )$ be its Poisson    enveloping algebra. Then 
\begin{enumerate}
\item $\CU (\CP ) \simeq \CP \rtimes_\pad U (\CP )/\CI (\CP )$ where $\CI (\CP )= (\d_{ab} -a\d_b - b\d_a )_{a,b\in \CP }$ is the ideal of the algebra $\CP \rtimes_\pad U (\CP )$ generated by the set $ \{ \d_{ab} -a\d_b - b\d_a \, | \, a,b\in \CP \}$.
\item If $\CP = S^{-1}K[x_i]_{i\in \L} / (f_s)_{s\in \G}$ where $S$ is a multiplicative subset  of the polynomial algebra $K[x_i]_{i\in \L}$ ($\L$ and $\G$ are index sets). Then the algebra $\CU (\CP )$ is generated by the algebra $\CP$ and the elements $\{ \d_i:= \d_{x_i}\, | \, i\in \L \}$ subject to the defining relations (a)--(c): For all elements $ i,j\in \L$ such that $i\neq j$ and $s\in \G$,
\begin{enumerate}
\item $[\d_i , \d_j] = \sum_{k\in \L} \frac{\der \{ x_i, x_j\}}{\der x_k}\d_k$, 
\item $[\d_i, x_j]= \{ x_i, x_j\}$, and 
\item $\sum_{i\in \L} \frac{\der f_s}{\der x_i}\d_i=0$. 
\end{enumerate}
So, the algebra $\CU (\CP )$ is generated by the algebra $\CP$ and the set $\d_\CP =\{ \d_a\, | \, a\in \CP \}$ subject to the defining relations: For all elenents $a,b\in \CP$ and $\l , \mu \in K$, 
\begin{enumerate}
\item $[\d_a, \d_b ] = \d_{\{ a,b\} }$, 
\item $[\d_a , b] = \{ a,b\}$, 
\item $\d_{ab} = a\d_b + b\d_a$,  
\item $\d_{\l a+\mu b} = \l \d_a+\mu \d_b$ and $\d_1=0$. 
\end{enumerate}
\item The map $\pi_\CP : \CU (\CP )\ra \CD (\CP )$, $a\mapsto a$, $ \d_b \mapsto \pad_b= \{ b, \cdot \}$ is an algebra homomorphism where $a, b\in \CP$ and its image is the algebra $P\CD (\CP )$ of Poisson differential operators of the Poisson algebra $\CP$. 
\item The algebra $\CP$ is a subalgebra of $\CU (\CP )$. Furthermore, $\CU (\CP ) = \CP \oplus \ann_{\CU (\CP )}(1)$ is a direct sum of left $\CP$-modules where $\ann_{\CU (\CP )}(1)=\sum_{i\in \L}\CU (\CP) \d_i$ is the annihilator of the identity element of the Poisson $\CP$-module $\CP$. The Poisson $\CP$-module structure on the Poisson algebra $\CP$ is obtained from the $\CD (\CP)$-module structure on $\CP$ by restriction of scalars via $\pi_\CP$. 
\end{enumerate}
\end{theorem}

{\it Proof}. 1. Consider the triple $(\CU' , \alpha, \beta )$ where   $ \CU':= \CP \rtimes_\pad U (\CP )/ \CI (\CP )$,  $\alpha : \CP \ra \CU'$, $a\mapsto a$ (it is an algebra homomorphism),  $\beta : \CP \ra \CU'$, $b\mapsto \d_b=\ad_b :=[b, \cdot ] $ (it  is  a Lie algebra homomorphism). Then for all elements $a,b\in \CP$, 
$$ [\beta (a) , \alpha (b) ] = \alpha (\{ a,b\})\;\; {\rm and}\;\; \beta (ab) = \beta (a)\alpha (b)+\alpha (a)  \beta (b).  $$
Given an associative algebra $A$, a homomorphism $\alpha' : \CP \ra A$ and a Lie homomorphism $\beta' : \CP \ra A$ such that $[\beta' (a) , \alpha' (b) ] = \alpha' (\{ a,b\})$ and $\beta' (ab) = \beta' (a)\alpha' (b)+\alpha' (a)  \beta' (b)$  for all elements $a, b\in \CP$. Then by the universal property of the semidirect product there is a unique homomorphism 
 $$ \CU':= \CP \rtimes_\pad U (\CP )/ \CI (\CP ) \ra A, \;\; a+\CI (\CP ) \mapsto \alpha' (a), \;\; \d_b+\CI (\CP ) \mapsto \beta' (b)$$
  such that $ \alpha'= f\alpha $  and $\beta' = f\beta$. Therefore, $\CU (\CP ) \simeq \CU'$.   
  
2. Statement 2 follows at once  from statement 1 as the relations (a) -- (b) of statement 1 follow from the axioms of Poisson algebra: For all elements $ i,j\in \L$ such that $i\neq j$ and $s\in \G$,

$[\d_i , \d_j] =[\d_{x_i} , \d_{x_j}]=\d_{\{ x_i, x_j\}}= \sum_{k\in \L} \frac{\der \{ x_i, x_j\}}{\der x_k}\d_k$, 

 $[ \d_i, x_j]= [ \d_{x_i}, x_j]=\{ x_i, x_j\}$,
  
$ 0=\d_0= \d_{f_s}= \sum_{i\in \L} \frac{\der f_s}{\der x_i}\d_i=0$.

3. Statement 3 follows from statement 1.

4. By the definition of the homomorphism $\pi_\CP$, the Poisson $\CP$-module structure on the Poisson algebra $\CP$ is obtained from the $\CD (\CP)$-module structure on $\CP$ by restriction of scalars via $\pi_\CP$. Since $\CU (\CP ) = \CP +I$, where $I=\sum_{i\in \L}\CU (\CP) \d_i$, $I\subseteq \ann_{\CU (\CP)}(1)$ and $\CP\cap \ann_{\CU (\CP)}(1)$, we must have $\CU (\CP ) = \CP \oplus I$ and $I= \ann_{\CU (\CP)}(1)=0$, and statement 4 follows. $\Box $

\begin{corollary}\label{ab23Jun19}

Every homomorphism of Poisson algebras $f:\CP\ra\CP'$ can be extended to a homomorphism of their Poisson enveloping algebras $f: \CU (\CP) \ra \CU (\CP')$ by the rule $f(\d_a)=\d_{f(a)}$.
\end{corollary}

{\it Proof}. The corollary follows from Theorem \ref{23Jun19}.(1,2). The map $f$ is well-defined since $f(\d_{ab}-a\d_b-b\d_a)=\d_{f(a)f(b)}-f(a)\d_{f(b)}-f(b)\d_{f(a)}$ for all elements $a,b\in \CP$. To finish the proof of the corollary it suffices to show that the relations (a)--(d) of 
Theorem \ref{23Jun19}.(2) holds. Let us check that the relation (a) holds. The other relations can be verified in a similar way. For all elements $a,b\in \CP$,
$$ f([\d_a,\d_b])=[\d_{f(a), \d_{f(b)}}]=\{f(a),f(b)\}=f(\{ a,b\}).\;\;\Box$$

{\bf The PEA  of a trivial Poisson algebra.} Corollary \ref{a23Jun19} describes the PEA  of a trivial Poisson algebra, i.e. $\{ \cdot , \cdot \}=0$.

\begin{corollary}\label{a23Jun19}
Suppose that  the algebra $\CP = S^{-1}K[x_i]_{i\in \L} / (f_s)_{s\in \G}$ is  a trivial Poisson algebra (i.e. $\{ \cdot , \cdot \}$)  where $S$ is a multiplicative subset of regular elements of the polynomial algebra $K[x_i]_{i\in \L}$. Then $\CU (\CP ) \simeq \CP [ \d_{x_i}]_{i\in \L}/(\sum_{i\in \L}\frac{\der f_s}{\der x_i}\d_{x_i})_{s\in \G}$. In particular, $\CU (K[x_i]_{i\in \L}) =K[x_i, \d_{x_i}]_{i\in \L} \simeq K[x_i]_{i\in \L}\otimes K[x_i]_{i\in \L}$ is polynomial algebra. 
\end{corollary}

{\it Proof}. The corollary follows at once from Theorem \ref{23Jun19}. 
 $\Box $\\
 
 {\bf Criterion for the algebra $\CU (\CP )$ to be a commutative algebra.} 
Corollary \ref{c23Jun19} is such a criterion.

\begin{corollary}\label{c23Jun19}
Let $\CP$ be a Poisson algebra. Then the algebra $\CU (\CP )$ is a commutative algebra iff the Poisson algebra $\CP$ is a trivial Poisson algebra. 
\end{corollary}

{\it Proof}. $(\Rightarrow )$ Suppose that  the algebra $\CU (\CP )$ is a commutative algebra. By Theorem \ref{23Jun19}.(4), $\CP \subseteq \CU (\CP )$. Then, by Theorem \ref{23Jun19}.(2), $0=\{ \d_i, x_j\} = \{ x_i, x_j \} \in \CP$, i.e. $\CP$ is a trivial Poisson algebra. 

$(\Leftarrow )$ This implication follows from Corollary \ref{a23Jun19} or  Theorem \ref{23Jun19}.(2). $\Box$ \\

{\bf Generators  and defining relations of the PEA  of a factor algebra of a Poisson algebra.} Proposition \ref{A23Jun19} represents the PEA  of a factor algebra of a Poisson algebra as a factor algebra the PEA  of the original Poisson algebra.

\begin{proposition}\label{A23Jun19}
 Let $\CP$ be a Poisson algebra, $\ga$ be a Poisson ideal of $\CP$, $\bCP =\CP / \ga$ be the Poisson factor algebra of $\CP$ and $\CP \ra \bCP$, $a\mapsto \oa = a+\ga$.  Then $\CU (\bCP )\simeq \CU (\CP ) / (\ga , \d_\ga )$ where $\d_\ga = \{ \d_a\, | \, a\in \ga \}$. Furthermore, if $\{ a_i\}_{i\in I} $ is a set of generators of the ideal $\ga$ then $\CU (\bCP )\simeq \CU (\CP ) / (\ga , \d_{a_i} )_{i\in I}$.
\end{proposition}

{\it Proof}. By Theorem \ref{23Jun19}.(1), the map $\CU (\CP ) \ra \CU (\bCP )$, $a\mapsto \oa$, $\d_a\mapsto \d_{\oa}$ (where $a\in \CP$) is an algebra epimorphism  as the relations (a)--(d) in Theorem \ref{23Jun19} are respected by the map. Then, by Theorem \ref{23Jun19}.(1), the statement follows.   $\Box $\\

{\bf The PEA  of a tensor product of Poisson algebras.} Let $\CP_1$ and $\CP_2$ be Poisson algebras. Their tensor product $\CP = \CP_1\t \CP_2$ is a Poisson algebra
 with respect to a unique Poisson bracket that is determined by the Poisson brackets of the tensor components and the condition that $\{ \CP_1, \CP_2\} =0$: For all elements $a_1,a_2\in \CP_1$ and $b_1, b_2\in \CP_2$, 
\begin{equation}\label{TPPAs}
\{ a_1\t b_1, a_2\t b_2\} = \{ a_1, a_2\} \t b_1b_2 + a_1a_2\t \{ b_1, b_2\}.  
\end{equation}

{\it Example.} The Poisson algebra $P_{2n} = K[x_1, \ldots , x_n, y_1, \ldots , y_n]$ where $\{ x_i, x_j\}=0$, $\{ y_i, y_j\}=0$ and $\{ y_i, x_j\}=\d_{ij}$ (the Kronecker delta) is a tensor product $P_{2n}= P_2^{\t n}$ of $n$ copies of the Poisson algebra $P_2= K[x,y]$ where $\{ y,x\} =1$.\\

Proposition \ref{B23Jun19} shows that the PEA  of a tensor product of Poisson algebras is a tensor product of PEA s of the tensor components.

\begin{proposition}\label{B23Jun19}
 Let $\CP = \CP_1\t \CP_2$ be a tensor product of Poisson algebras $\CP_1$ and $\CP_2$. Then $\CU (\CP )\simeq \CU (\CP_1)\t \CU (\CP_2)$.
\end{proposition}

{\it Proof}. The result follows from  Theorem \ref{23Jun19}. $\Box$

\begin{proposition}\label{AB23Jun19}
Let $P_{2n} = K[x_1, \ldots , x_n, y_1, \ldots , y_n]$ where $\{ x_i, x_j\}=0$, $\{ y_i, y_j\}=0$ and $\{ y_i, x_j\}=\d_{ij}$ (the Kronecker delta). Then the algebra $\CU (P_{2n})$ is isomorphic to the Weyl algebra $A_{2n}$. 
\end{proposition}

{\it Proof}. Recall that $P_{2n} \simeq P_2^{\t n}$, see above. By Proposition \ref{B23Jun19}, $\CU (P_{2n})\simeq \CU (P_2)^{\t n}$. So, it suffices to prove the statement for $n=1$, i.e. $P_2= K[x,y]$ with $\{ y,x\} =1$.  By Theorem \ref{23Jun19}.(2), the algebra $\CU (P_2)$ is generated by the elements $x$, $y$, $\d_x$ and $\d_y$ subject to the defining relations: $xy=yx$, $[\d_x, \d_y]=0$, $[\d_y, x]=\{ y,x\} =1$, $[\d_y, y\} = \{ y,y\} =0$, $[\d_x, y] = \{ x,y\}=-1$ and $[\d_x, x]= \{ x,x\}=0$, i.e. $\CU (P_2) \simeq A_2$. $\Box$\\

For an associative algebra $A$,  we denote by $\End_K(A)$ (resp., $\Aut_K(A)$) the set of all algebra endomorphisms (resp., automorphisms) of $A$. Clearly, $\End_K(A)$ is a monoid  w.r.t. the composition of maps and $\Aut_K(A)$ is its group of units. For a Poisson algebra $\CP$, we denote by $\End_{\rm Pois} (\CP )$ (resp., $\Aut_{\rm Pois} (\CP ))$  the set of all Poisson  algebra endomorphisms (resp., automorphisms) of $A$. Clearly, $\End_{\rm Pois} (\CP )$ is a monoid  w.r.t. the composition of maps and $\Aut_{\rm Pois} (\CP )$ is its group of units.

Proposition \ref{D23Jun19} shows that every endomorphism/automorphism of a Poisson algebra is extended to an endomorphism/automorphism of its Poisson    enveloping algebra.

\begin{proposition}\label{D23Jun19}
Let $\CP$ is a Poisson algebra. Then 
\begin{enumerate}
\item The map $\End_{\rm Pois} (\CP ) \ra \End_K(\CU (\CP ))$, $\s\mapsto \s : a\mapsto \s (a)$, $\d_a\mapsto \d_{\s (a)}$ ($a\in \CP $) is a monoid monomorphism. 
\item The map $\Aut_{\rm Pois} (\CP ) \ra \Aut_K(\CU (\CP ))$, $\s\mapsto \s : a\mapsto \s (a)$, $\d_a\mapsto \d_{\s (a)}$ ($a\in \CP $) is a group monomorphism. 
\end{enumerate}
\end{proposition}

{\it Proof}. 1. The proposition follows from Theorem \ref{23Jun19}. $\Box$ \\

{\bf Localizations commuts with the operation of taking PEA  of a Poisson algebra.} 

\begin{theorem}\label{30Jun19}
Let $\CP$ be a Poisson algebra and $S$ be a multiplicative subset of $\CP$, $\ga = \ass (S)=\{ a\in \CP \, | \, sa =0$ for some $s\in S\}$. Then the ideal $\ga $ is a Poisson ideal, $\bCP = \CP / \ga$ is a Poisson algebra  and $\bS = \{ \bs \, | \, s\in S\}$ is a multiplicative subset of $\bCP$ that consists of regular elements of the algebra $\bCP$, $S^{-1}\CP \simeq \bS^{-1}\bCP$, $S\in \Den (\CU (\CP ), \ga )$, $\bS \in \Den ( \CU (\bCP ), 0)$, $\CU (S^{-1} \CP ) \simeq \CU ( \bS^{-1} \bCP ) \simeq S^{-1} \CU ( \CP )\simeq \CU ( \CP ) S^{-1} \simeq \bS^{-1} \CU ( \bCP )\simeq  \CU ( \bCP )\bS^{-1}$. 
\end{theorem}

{\it Proof}. (i) $\ga $ {\em is a Poisson ideal}: By the very definition, $\ga$ is an ideal of $\CP$. It remains to show that $\{ \CP , \ga \} \subseteq \ga$. Given elements $a\in \ga$ and $b\in \CP$. Then $sa=0$ for some element $s\in S$. Then $0=\{ b, 0\} = \{ b , s^2a\} = 2sa \{ b, s\} + s^2\{ b,a\} = s^2 \{ b,a\} $, and so $\{ b, a\} \in \ga$. 

By the statement (i), $\bCP$ is a Poisson algebra, the set $\bS$ consists of regular elements of $\bCP$ and $S^{-1}\CP \simeq \bS^{-1} \bCP$. 

(ii) $S\in \Den (\CU (\CP ))$: The statement (ii) follows from Theorem \ref{23Jun19} and the fact that for all elements $a\in \CP$ and $s\in S$, $\d_as^2-s^2\d_a= \{ a, s^2\} = 2s\{ a,s\}$. 

(iii) $\CU ( S^{-1} \CP ) \simeq S^{-1} \CU (\CP ) \simeq \CU ( \CP  S^{-1})\simeq \CU ( \CP )  S^{-1}$:  




By the universal property of the PEA ,
$\CU (S^{-1}\CP ) \simeq  S^{-1}\CU ( \CP )$. By the statement (ii), $S^{-1}\CU ( \CP )\simeq\CU ( \CP ) S^{-1}$, and the statement (iii) follows.

(iv) $\CU (S^{-1} \CP ) \simeq \CU ( \bS^{-1} \bCP ) \simeq  \bS^{-1} \CU ( \bCP )\simeq  \CU ( \bCP )\bS^{-1}$: The statement (iv) follows from the statement (iii). $\Box$ \\

{\bf The PEA  of the Poisson symmetric algebra $S(\CG )$ of a Lie algebra $\CG $.} Let $\CG$ be a Lie algebra over the field $K$, $\{ x_i\}_{i\in I}$ be its $K$-basis and the lie bracket  on $\CG$, 
$$ [ x_i, x_j] = \sum_{k\in I} c_{ij}^k x_k,$$
is determined by the structure constants $c_{ij}^k\in K$. The universal enveloping algebra $U = U(\CG )$ of $\CG$ admits the standard filtration $\{ U_i\}_{i\in \N }$ by the total degree of the generators   $\{ x_i\}_{i\in I}$ of $U$. The associated graded algebra $S(\CG ) := {\rm gr} (U) = \bigoplus_{i\geq 0 } U_i/ U_{i-1}$ (where $U_{-1}=0$), the so-called, {\em symmetric algebra} of $\CG$, is a polynomial algebra in the  variables  $\{ x_i\}_{i\in I}$. Furthermore, the commutative algebra $S(\CG )$ is a Poisson algebra where 
\begin{equation}\label{SGPA}
\{  x_i, x_j\} = \sum_{k\in I} c_{ij}^k x_k.
\end{equation}

\begin{proposition}\label{C23Jun19}
Let $\CG$, $U(\CG)$ and $(S(\CG ) , \{ \cdot , \cdot \})$ be as above. Then $\CU (S(\CG ))\simeq S(\CG ) \rtimes_\pad U(\CG )$ where $\pad : \CG \ra \PDer_K(S(\CG ))$, $x\mapsto \pad_x = \{ x, \cdot \}$. In more detail, the algebra $\CU (S(\CG ))$ is generated by the set $\{ x_i , \d_i\, | \, i\in I\}$ subject to the defining relations: For all elements $i,j\in I$, $x_ix_j = x_jx_i$, $[\d_{x_i} , \d_{x_j}]=\sum_{k\in I} c_{ij}^k \d_{x_k}$ and $[\d_{x_i} , x_j]=\{ x_i, x_j\} $ $(=\sum_{k\in I} c_{ij}^k x_k$). The subalgebra $K\langle \d_{x_i}\, | \, i\in I\rangle$ of $\CU (S(\CG ))$ is isomorphic to the universal enveloping algebra $U(\CG )$  of the Lie algebra $\CG$ via $\d_{x_i} \mapsto x_i$. 
\end{proposition}

{\it Proof}. We keep the notation as above. By Theorem \ref{23Jun19}, the algebra $\CU = \CU (S(\CG ))$ is generated by the elements $\{ x_i, \d_{x_i}\, | \, i\in I\}$ subject to the defining relations as in the theorem. The algebra $S(\CG ) \rtimes_\pad U(\CG )$ is generated by the same set of generators $\{ x_i, \d_{x_i}\, | \, i\in I\}$ subject to the same defining relations ($S(\CG )=K[x_i]_{i\in I}$ and $U(\CG )=K\langle \d_{x_i}\rangle_{i\in I}$).
Therefore, $\CU (S(\CG ))\simeq S(\CG ) \rtimes_\pad U(\CG )$.  $\Box $

\begin{corollary}\label{aC23Jun19}
Let $\CG = \CG_1\times \CG_2$ be a direct product of Lie algebras. Then $\CU (S(\CG ))\simeq \CU (S(\CG_1 ))\t \CU (S(\CG_2 ))$ is a direct product of algebras. 
\end{corollary}

{\it Proof}. $S(\CG ) = S(\CG_1\times \CG_2) \simeq S( \CG_1) \t S(\CG_2)$, a tensor product of Poisson algebras. Now, the corollary follows from Proposition \ref{B23Jun19}.  $\Box $\\

{\bf $\CU ( \CG_1\rtimes_\d \CG_2)$ where $\CG_1\rtimes_\d \CG_2$ be  a semidirect product of Lie algebras.}  Let  $\CG :=\CG_1\rtimes_\d \CG_2$ be  a semidirect product of Lie algebras where $\d : \CG_2\ra \Der_{\rm Lie} (\CG_1)$, $x_2\mapsto \d_{x_2}$ is a Lie homomorphism $(\d ([x_2,y_2])= [ \d_{x_2}, \d_{y_2}]$ for all $x_2,y_2\in \CG_2$),  
i.e. the Lie bracket on $\CG$ is given by the rule: For all elements $x_1, y_1\in \CG_1$ and $x_2, y_2\in \CG_2$, 
 $$[x_1+x_2, y_1+y_2]=  [x_1, y_1]_{\CG_1} + [x_2, y_2]_{\CG_2} + \d_{x_2}(y_1) - \d_{y_2}(x_1).$$
Then 
$$ U(\CG ) \simeq U(\CG_1) \rtimes_\d U(\CG_2), $$
and so $ S(\CG )\simeq S( \CG_1) \rtimes_\d  S(\CG_2)$ $(=S(\CG_1) \t S(\CG_2)$, a tensor product of algebra), the Poisson bracket on the algebra $S(\CG )$ is given by the rule: For all elements $x_1, y_1\in \CG_1$ and $x_2, y_2\in \CG_2$, 
$$\{ x_1, y_1\} = [x_1, y_1]_{\CG_1}, \;\; \{ x_2, y_2\} = [x_2, y_2]_{\CG_2} \;\; 
{\rm and}\;\; \{ x_2, x_1\} = \d_{x_2}(x_1).$$

\begin{corollary}\label{bC23Jun19}
Let $\CG = \CG_1\rtimes \CG_2$ be a semidirect product of Lie algebras. Then $\CU (S(\CG ))\simeq (S( \CG_1) \rtimes_\d  S(\CG_2))\rtimes_\ad (U( \CG_1) \rtimes_\d  U(\CG_2))$. 
\end{corollary}

{\it Proof}. By Proposition \ref{C23Jun19}, $\CU (S(\CG ))\simeq S(\CG ) \rtimes_\ad U(\CG )$, and the result follows (since  $ S(\CG )\simeq S( \CG_1) \rtimes_\d  S(\CG_2)$ and $U(\CG ) \simeq U( \CG_1) \rtimes_\d  U(\CG_2))$. $\Box $


\section{The PBW Theorem  for the Poisson  enveloping algebras and the module $\OCP$ of K\"{a}hler differentials of a Poisson algebra $\CP$}\label{MKAHLD}

The aim of this section is to prove the  PBW Theorem  for the Poisson  enveloping algebras (Theorem \ref{A30Jul19}); to give a projectivity criterion  for the algebras $\CU (\CP)$ and ${\rm }\,\CU (\CP)$ (Corollary \ref{31Jul19});  to prove that the algebra $\CU ( \CP )$ is isomorphic to its  opposite algebra $\CU ( \CP )^{op}$; a criterion for the algebra {\rm gr} $\CU (\CP ) $  to be a left/right Noetherian algebra or a finitely generated algebra is given (Proposition \ref{B31Jul19}).


The derivation algebras are introduced at the end of the section.  Theorem \ref{25Mar19} is a criterion for a derivation algebra to have PBW basis.  Theorem \ref{A25Mar19} and 
Theorem \ref{A11Jab20} show that  every polynomial Poisson algebra and every Poisson algebra with free module of K\"{a}hler differential are derivation algebras with PBW basis, respectively.\\

{\bf The algebras $F_b$, $F_{ab}$ and $\CU ( \CP )$. } Let $P_\L = K[x_i]_{i\in \L}$ be a polynomial algebra where the cardinality ${\rm card } (\L ) = |\L |$ of the set $\L $ can be arbitrary, $I=(f_s)_{s\in \G }$ be an ideal of $P_\L$, $\bP_\L :=P_\L / I$ and $\CP = S^{-1}\bP_\L$ where $S$ is a multiplicative subset in $\bP_\L$ that consists of regular  elements (i.e. non-zero-divisors) of $\bP_\L$, eg, $S= \{ 1\}$. Suppose that $\CP$ is a Poisson algebra. Let $F=\langle \CP , \d_{x_i}=\d_i\rangle_{i\in \L}$ be a $K$-algebra generated freely by the algebra $\CP$ and a set of  free generators $\{ \d_i\, | \, i\in \L\}$.   Let $\CR_a$, $\CR_{b}$ and $\CR_c$ be ideals of the algebra $F$ that are generated by the relations (a), (b)  and (c) of Theorem \ref{23Jun19}, respectively, 
\begin{eqnarray*}
 \CR_a&=& \bigg( [\d_i , \d_j] - \sum_{k\in \L} \frac{\der \{ x_i, x_j\}}{\der x_k}\d_k\bigg)_{i,j\in \L},\\
 \CR_b&=& \bigg([\d_i, x_j]- \{ x_i, x_j\}\bigg)_{i,j\in \L}, \\
\CR_c &=& \bigg(\sum_{i\in \L} \frac{\der f_s}{\der x_i}\d_i\bigg)_{s\in \G}. 
\end{eqnarray*}
Let $\CR_{ab}=(\CR_a, \CR_b)$ and $\CR = (\CR_a, \CR_b, \CR_c)$, ideals of $F$. By Theorem \ref{23Jun19}.(2), $\CU (\CP ) \simeq F/ \CR$. There are algebra epimorphisms
$$ F\ra F_b:=F/\CR_b\ra F_{ab}:= F/ \CR_{ab} \ra \CU (\CP ) \simeq F/ \CR.$$ 

\begin{proposition}\label{B29Jul19}
Let $\CP = S^{-1}(P_\L / I)$ be a Poisson algebra where $P_\L = K[x_i]_{i\in \L}$ be a polynomial algebra, $I=(f_s)_{s\in \G }$ be an ideal of $P_\L$, $\bP_\L :=P_\L / I$ and  $S$ be a multiplicative subset in $\bP_\L$ that consists of regular  elements  of $\bP_\L$. Then 
\begin{enumerate}
\item $F_b = \bigoplus_{\d \in \CM (\L )} \CP \d = \bigoplus_{\d \in \CM (\L )} \d \CP $ where $\CM (\L )$ is a free multiplicative monoid on the set $\{ \d_i\} _{i\in \L }$. 
\item $F_{ab} = \bigoplus_{\d \in \CM (\L )_c} \CP \d = \bigoplus_{\d \in \CM (\L )_c} \d \CP $ where $\CM (\L )_c$ is a free commutative  multiplicative monoid on the set $\{ \d_i\} _{i\in \L }$, i.e. $F_{ab} = \bigoplus_{\alpha \in \N^{(\L )}} \CP \d^\alpha  = \bigoplus_{\alpha \in \N^{(\L )}}  \d^\alpha \CP$ where $\d^\alpha = \prod_{i\in \L } \d_i^{\alpha_i}$ and $\alpha = (\alpha_i)_{i\in \L}$. In particular, the algebra $F_{ab}$ is a free left and right $\CP$-module.  If the algebra $\CP$ is Noetherian and $|\L |<\infty$ then the algebra $F_{ab}$ is a Noetherian algebra. 
\item The elements $\{ \d_{f_s} := \sum_{i\in \L} \frac{\der f_s}{\der x_i}\d_{x_i} \, | \, s\in \G\}\subseteq F_{ab}$ belong to the centralizer of the algebra $\CP$ in $F_{ab}$ (i.e. $[\d_{f_s}, x_j]=0$ for all $s\in \G$ and $j\in \L$); $[\d_{f_s}, \d_j]=\sum_{k\in \L} \frac{\der \{f_s,x_j\}}{\der x_k} \d_k$; $\CU (\CP ) \simeq F_{ab} / (\d_{f_s})_{s\in \G}$; and the ideal $\bCR_c=(\d_{f_s})_{s\in \G}$ of the algebra $F_{ab}$ is equal to the left and  right ideals of $F_{ab}$ that are generated by the elements $\{ \d_{f_s} \, | \, s\in \G\}$. 
\item If the algebra $\CP$ is Noetherian and $|\L |<\infty$ then the algebra $\CU (\CP )$ is a Noetherian algebra.
\end{enumerate}
\end{proposition}

{\it Proof}. 1. Statement 1 is obvious.

2. Fix a total ordering on the set $\{ \d_i:= \d_{x_i}\, | \, i\in \L \}$. It determines the deg-by-lexicographic ordering on the monoid $\CM (\L )$. In order to prove that the first equality  of statement 2  holds (about the direct sum) it suffices to show that the ambiguities $\{ \d_i\d_j a\, | \, i>j$ and $ a\in \CP \}$ and $\{ \d_i\d_j \d_k\, | \, i>j>k\}$ are resolved:
\begin{eqnarray*}
\bullet\;\;  \d_i\d_j a &=&\d_ia\d_j+ \d_i\{ x_j, a\}=a\d_i\d_j+\{ x_i, a\}\d_j+\{ x_j, a\} \d_i + \{ x_i, \{ x_j, a\}\} \\
 &=& a\d_j\d_i+a\sum_{k\in \L} \frac{\der\{ x_i, x_j\} }{\der x_k}\d_k+\{ x_i, a\}\d_j+\{ x_j, a\} \d_i + \{ x_i, \{ x_j, a\}\} \\
\bullet\;\;  \d_i\d_j a &=& \d_j\d_ia+\sum_{k\in \L} \frac{\der\{ x_i, x_j\} }{\der x_k}\d_ka=a\d_j\d_i+\{ x_j, a\} \d_i + \d_j\{ x_i, a\} + \sum_{k\in \L} \frac{\der\{ x_i, x_j\} }{\der x_k}(a\d_k+\{ x_k, a\}) \\
 &=& a\d_j\d_i+\{ x_j, a\} \d_i + \{ x_i, a\}\d_j + \{ x_j, \{ x_i, a\} \}+a\sum_{k\in \L} \frac{\der\{ x_i, x_j\} }{\der x_k}\d_k+\{ \{ x_i, x_j\} , a\}.
\end{eqnarray*}
The RHSs of both equalities are equal since $\{ x_i, \{ x_j, a\}\}= \{ x_j, \{ x_i, a\} \}+\{ \{ x_i, x_j\} , a\}$. Similarly,
\begin{eqnarray*}
\bullet\;\;  \d_i\d_j \d_k &=& \d_j\d_i\d_k+ \d_{\{ x_i, x_j\} } \d_k= \d_j\d_k\d_i+\d_j \d_{\{ x_i, x_k\} } + \d_{\{ x_i, x_j\} } \d_k\\
&=& \d_k\d_j\d_i + \d_{\{ x_j, x_k\} } \d_i +\d_{\{ x_i, x_k\} } \d_j+\d_{\{ x_j , \{ x_i, x_k\} \} }  +\d_{\{ x_i, x_j\} }\d_k.\\
\bullet\;\;  \d_i\d_j \d_k &=& \d_i\d_k \d_j+ \d_i\d_{\{ x_j , x_k\} }= \d_k \d_i \d_j + \d_{\{ x_i, x_k\} }\d_j+ \d_{\{ x_j, x_k\} } \d_i+ \d_{\{ x_i , \{ x_j , x_k\} \} } \\
&=&\d_k\d_j\d_i+ \d_{\{ x_i , x_j\} } \d_k+ \d_{ \{ x_k , \{ x_i , x_j \} \}  } + \d_{ \{ x_i , x_k \} } \d_j + \d_{\{ x_j , x_k \} } \d_i+ \d_{\{ x_i , \{ x_j , x_k\} \} }.
\end{eqnarray*}
The RHSs of both equalities are equal since $\{ x_j, \{ x_i, x_k\}\}= \{ x_k, \{ x_i, x_j\} \}+\{x_i,  \{ x_j, x_k\} \}$. So, $F_{ab} = \bigoplus_{\d \in \CM (\L )} \CP \d$. Hence, $F_{ab} = \bigoplus_{\d \in \CM (\L )}\d  \CP $.

 The associated graded algebra ${\rm gr}  (F_{ab})$ of the algebra $F_{ab}$ w.r.t. the filtration on the algebra $F_{ab}$ by the total degree of the elements $\{ \d_i \, | \, i\in \L \}$ is a polynomial algebra $\CP [\d_i]_{i\in \L }$ which  is a Noetherian algebra provided  the algebra $\CP $ is Noetherian and $|\L |<\infty$. Hence, under these conditions  the algebra $F_{ab}$ is also  Noetherian.  

3. For all $s\in \G$ and $j\in \L$, we have the following equalities in the algebra $F_{ab}$:
\begin{eqnarray*}
 [\d_{f_s} , x_j] &=& \bigg[\sum_{i\in \L } \frac{\der f_s}{\der x_i}\d_i, x_j\bigg]=\sum_{i\in \L }  \frac{\der f_s}{\der x_i}[\d_i, x_j]=\sum_{i\in \L }  \frac{\der f_s}{\der x_i}\{ x_i, x_j\}\\
 &=&\{ f_s, x_j\} =0,\\
 \d_{f_s}  \d_j- \d_j \d_{f_s}&=& \bigg[\sum_{i\in \L } \frac{\der f_s}{\der x_i}\d_i, \d_j\bigg]= \sum_{i\in \L} \bigg( \bigg[\frac{\der f_s}{\der x_i}, \d_j\bigg]\d_i+\frac{\der f_s}{\der x_i}[\d_i, \d_j]\bigg)\\
 &=& \sum_{i\in \L }\bigg( \sum_{k\in \L }  \frac{\der^2f_s}{\der x_k \der x_i} \{ x_k , x_j\} \d_i+  \frac{\der f_s}{\der x_i}  \sum_{k\in \L } \frac{\der \{ x_i, x_j \} }{\der x_k} \d_k\bigg)\\
 &=& \sum_{k\in \L } \frac{\der}{\der x_k } \bigg( \sum_{i\in \L }\frac{\der f_s}{\der x_i}\{ x_i, x_j\} \bigg) \d_k\\
 &=&\sum_{k\in \L } \frac{\der \{ f_s, x_j\} }{\der x_k}\d_k .
\end{eqnarray*}
Clearly, $\CU (\CP ) \simeq F_{ab} / (\d_{f_s})_{s\in \G}$, and the rest of the statement 3 follows. 

4.  The algebra $F_{ab}$ is a Noetherian algebra, by statement 2. Hence, so is its factor algebra $\CU ( \CP )$, by statement 3.  $\Box $\\

{\bf The Poisson    enveloping algebra  is isomorphic to its opposite algebra.} Let $A$ be an algebra. The {\em opposite algebra} $A^{op}$ of $A$ is an algebra such that $A^{op} =A$ (the equality of vector spaces) and the product in $A^{op}$ is given by the rule: For all elements $a, b\in A^{op}$, $a\cdot b = ba$.   Theorem \ref{C29Jul19} shows the Poisson    enveloping algebra  is isomorphic to its opposite algebra.

\begin{theorem}\label{C29Jul19}
We keep the notation of Proposition \ref{B29Jul19}. 
\begin{enumerate}
\item The map $F_{ab} \ra F_{ab}^{op}$, $a\mapsto a$, $\d_i\mapsto - \d_i$ is an algebra isomorphism where  $a\in \CP$ and $i\in \L$. 
\item The map $\CU (\CP ) \ra \CU (\CP )^{op}$, $a\mapsto a$, $\d_i\mapsto - \d_i$ is an algebra isomorphism where  $a\in \CP$ and $i\in \L$, i.e. for each Poisson algebra its PEA  $\CU$ is isomorphic to its opposite algebra $\CU^{op}$ and as a result properties of the algebra $\CU$ is left-right symmetric. 
\item For all elements $f\in \CP$, $\sum_{j\in \L} \{ \frac{\der f}{\der x_j}, x_j\}=0$,
\end{enumerate}
\end{theorem}

{\it Proof}. 3. Suppose that char$(K)\neq 2$.  Then  \begin{eqnarray*}
 a&=& \sum_{j\in \L} \{ \frac{\der f}{\der x_j}, x_j\}=\sum_{i,j\in \L} \frac{\der^2 f}{\der x_i\der x_j}\{ x_i, x_j\}=\sum_{i,j\in \L}  \frac{\der^2 f}{\der x_j\der x_i}\{ x_j, x_i\}\\
 &=&-\sum_{i,j\in \L}  \frac{\der^2 f}{\der x_i\der x_j}\{ x_i, x_j\}=-a, 
\end{eqnarray*}
i.e. $2a=a$, and so $a=0$ since ${\rm char} (K)\neq 2$. 

In the case  char$(K)= 2$, it suffices to show that the equality holds for monomials  that are products of distinct $x_i$'s. Let $m =x_1\cdots x_l$, $m_i$ and $m_{ij}$ ($i\neq j$) be the monomial $m$ where the elements $x_i$ and $x_ix_j$ are deleted. Then
$$\sum_{j\in \L} \{ \frac{\der m}{\der x_j}, x_j\}=\sum_{j=1}^l \{ m_j, x_j\}=\sum_{i\neq j=1}^l m_{ij}\{ x_i, x_j\}=2\sum_{i< j=1}^l m_{ij}\{ x_i, x_j\}=0.$$

1 and 2.  The map, say $\th$, in the statements 1 and 2 respects each of the defining relations $\CR_a$, $\CR_b$ and $\CR_c$ (by using statement 3):
\begin{eqnarray*}
 \CR_a&:& \th ([\d_i, \d_j]-\sum_{k\in \L} \frac{\der \{ x_i ,  x_j\} }{\der x_k}\d_k)=[-\d_j, -\d_i]-\sum_{k\in \L}(-\d_k) \frac{\der \{ x_i ,  x_j\} }{\der x_k}\\
 &=& -\bigg( [\d_i, \d_j]-\sum_{k\in \L}\frac{\der \{ x_i ,  x_j\} }{\der x_k}\d_k-\sum_{k\in \L} \{ x_k, \frac{\der \{ x_i ,  x_j\} }{\der x_k}\} \bigg)= -
 \bigg( [\d_i, \d_j]-\sum_{k\in \L}\frac{\der \{ x_i ,  x_j\} }{\der x_k}\d_k\bigg).\\
 \CR_b &:& \th ([\d_i, x_j]-\{x_i, x_j\}) = [x_j , -\d_i]-\{ x_i, x_j \} = [\d_i, x_j]-\{x_i, x_j\}.\\
 \CR_c &:& \th (\d_{f_s})= \th (\sum_{i\in \L} \frac{\der f_s}{\der x_i}\d_i)=\sum_{i\in \L}(-\d_i ) \frac{\der f_s}{\der x_i} =-\d_{f_s}-\sum_{i\in \L} \{ x_i,  \frac{\der f_s}{\der x_i}\}=-\d_{f_s}.\;\; \Box 
  \end{eqnarray*}

{\bf The opposite Poisson algebra $\CP^{op}$ of a  Poisson algebra $\CP$ and $\CU (\CP^{op})\simeq \CU (\CP) ^{op}$.} Let $\CP$ be a Poisson algebra. The {\em opposite Poisson algebra} of $\CP$, denoted by $\CP^{op}$, is a Poisson algebra that coincides with $\CP$ as an associative algebra but the Poisson bracket $\{ \cdot , \cdot \}^{op}$ in $\CP^{op}$ is given by the rule: $\{ a,b\}^{op}=\{b,a\}$. We can easily verify that $(\CP^{op}, \{ \cdot , \cdot \}^{op})$ is a Poisson algebra. Clearly, $(\CP^{op})^{op}\simeq \CP$. 

\begin{proposition}\label{D29Jul19}
 Let $\CP$ be a Poisson algebra. Then $\CU (\CP^{op})\simeq \CU (\CP )^{op}\simeq \CU (\CP)$.
\end{proposition}

{\it Proof}. In view of Theorem \ref{C29Jul19}.(2), we  have to show that only the first isomorphism of the proposition holds. Let $\CP$ be as in Theorem \ref{C29Jul19}.(2). Then the map 
\begin{equation}\label{CUopCU}
\phi : \CU (\CP^{op})\ra  \CU (\CP )^{op}, \;\ a\mapsto a, \;\; \d_i\mapsto \d_i \;\; (a\in \CP, \; i\in \L )
\end{equation}
is an algebra isomorphism. By the very definition, the map is an algebra isomorphism provided it is well-defined. In order to check this we have to show
 that the map $\phi$ respects the defining relations of both algebras, i.e. the relations (a)--(c) in Theorem \ref{23Jun19}.(2). This can be easily checked directly by using Theorem \ref{C23Jun19}.(2). Furthermore, the map $\phi$ respects {\em each type} of the relations (a)--(c).  For example, the image under the map $\phi$ of a relation of the type (b) is a relation again of the type (b):
 $$\phi ([\d_i, x_j]-\{ x_i,x_j\}^{op})=\phi ([\d_i, x_j]-\{ x_j,x_i\})= [x_j,\d_i]-\{ x_j,x_i\}=-([\d_i, x_j]-\{ x_i,x_j\})=-0=0. \;\; \Box $$
 
{\bf The module $\OCP$ of K\"{a}hler differentials of a Poisson algebra $\CP$ and $\CU ( \CP )$.} Let $\CP$ be a Poisson algebra and let $\O = \OCP$ be the {\em module of K\"{a}hler differentials} of the (associative)  algebra $\CP$. In more detail, if $\CP = S^{-1} P_\L/I$ where $I=(f_s)_{s\in \G}$ then 
\begin{equation}\label{KPdif}
\O =\O_{S^{-1} P_\L /I} \simeq \bigoplus_{i\in \L} \CP dx_i/\sum_{s\in \G} \CP df_s \;\; {\rm where} \;\; df_s=\sum_{i\in \L} \frac{\der f_s}{\der x_i}dx_i. 
\end{equation}
Notice that $\mF := \bigoplus_{i\in \L} \CP \d_i$ is a $\CP$-submodule of $F_{ab}$ (Proposition \ref{B29Jul19}.(2)) that contains the $\CP$-submodule $\sum_{s\in \G} \CP \d_{f_s}$. Then the map
\begin{equation}\label{KPdif1}
\OCP \ra \mF / \sum_{s\in \G} \CP \d_{f_s}, \;\; dx_i\mapsto \d_i \;\; (i\in \G )
\end{equation}
is a $\CP$-modules {\em isomorphism}. 

Let $\bCR_c$ be the image of the ideal $\CR_c$ of $F$ under the algebra epimorphism $F\ra F_{ab}$.  By Proposition \ref{B29Jul19}.(3),  $ \bCR_c$ is a left/right/two-sided  ideal of the algebra $F_{ab}$
 which is generated by the set $\{ \d_{f_s} \, | \, s\in \G\}$ of  of the algebra $F_{ab}$. By Proposition \ref{B29Jul19}.(3),
\begin{equation}\label{KPdif2}
\bCR_c=\sum_{s\in \G} \d_{f_s} F_{ab} = M_0F_{ab} = F_{ab}M_0\;\; {\rm where}\;\; M_0:= \sum_{s\in \G} \CP \d_{f_s}=\sum_{s\in \G}  \d_{f_s}\CP\subseteq F_{ab}. 
\end{equation}
By Proposition \ref{B29Jul19}.(2), 
$$F_{ab} =\bigoplus_{i\geq 0} F_{ab, i}=\bigoplus_{i\geq 0} F_{ab, i}' \;\; {\rm where}\;\; F_{ab,i}=\bigoplus_{|\alpha |=i} \CP \d^\alpha\;\; {\rm and} \;\; F_{ab,i}'=\bigoplus_{|\alpha |=i} \d^\alpha \CP.$$
So,  $F_{ab}$ is a direct sum of left $\CP$-modules $F_{ab, i}$ and right  $\CP$-modules $F_{ab, i}'$. By (\ref{KPdif2}) and  Proposition \ref{B29Jul19}.(3),  
\begin{equation}\label{KPdif3}
\bCR_c=\bigoplus_{i\geq 0} M_0 F_{ab, i} =\bigoplus_{i\geq 0} F_{ab, i}' M_0, \;\; M_0F_{ab,i}\subseteq F_{ab, i+1}\;\; {\rm and}\;\; F_{ab,i}'M_0\subseteq F_{ab, i+1}', \;\; i\geq 0.
\end{equation}
By (\ref{KPdif3}), the algebra 
\begin{equation}\label{KPdif4}
\CU (\CP ) = F_{ab}/\bCR_c=\bigoplus_{i\geq 0} F_{ab, i}/M_0F_{ab, i-1}= \bigoplus_{i\geq 0} F_{ab, i}'/F_{ab, i-1}'M_0
\end{equation}
is a direct sum of left $\CP$-modules $F_{ab, i}/M_0F_{ab, i-1}$ and right $\CP$-modules $F_{ab, i}'/F_{ab, i-1}'M_0$. By (\ref{KPdif1}) and (\ref{KPdif4}), 
\begin{equation}\label{KPdif5}
\CU (\CP ) = \CP \oplus \OCP \oplus \bigoplus_{i\geq 2} F_{ab, i}/M_0F_{ab, i-1}= \CP \oplus \OCP' \oplus \bigoplus_{i\geq 2} F_{ab, i}'/F_{ab, i-1}'M_0
\end{equation}
where $\OCP':=\mF'/M_0$ and $\mF':=\oplus_{i\in \L}\d_i\CP$. The first (resp., second) direct sum in (\ref{KPdif5}) is a direct sum of left (resp., right) $\CP$-modules. 
\begin{theorem}\label{30Jul19}
Let $\CP$ be a Poisson algebra and $\O = \OCP$ be the $\CP$-modules of K\"{a}hler differentials of $\CP$. Then 
\begin{enumerate}
\item  $\O$  is a Lie subalgebra of the Lie algebra $(\CU ( \CP ), [\cdot , \cdot ])$ where $[a,b]=ab-ba$. In particular, for all elements $a_1, a_2, b_1, b_2\in \CP$, 
$$[a_1da_2, b_1db_2]=a_1\{ a_2, b_1\} db_2-b_1\{ b_2, a_1\} da_2+a_1b_1d\{ a_2, b_2\}.$$
\item The left $\CP$-module homomorphism $\O \ra \Der_K(\CP )$, $adb\mapsto a\{ b, \cdot \}$ is a Lie algebra homomorphism.
\item The left $\CP$-submodule $\CP\oplus \O$ of $\CU (\CP )$ is a Lie subalgebra which is a semi-direct product $\CP \rtimes \O$ of the Lie algebra $\O$ by the abelian Lie algebra $\CP$ via the Lie homomorphism in statement 2, i.e. $[adb , p]=a\{ b,p\}$ for all elements $a,b,p\in \CP$. In particular, the Poisson structure on $\CP$ can be recovered from the Lie algebra structure on $\CP \rtimes \O$ by the rule: For all elements $a,b\in\CP$, $\{ a,b\} = [ da, b]$. 
\item $\CU (\CP ) \simeq \CP \rtimes U(\O )/R(\O )$ where $U(\O )$ is the universal enveloping algebra of the Lie $K$-algebra $\O$, $\CP \rtimes U(\O )$ is a skew product of the algebra $U(\O )$ and $\CP$ and $R(\O )= (a\cdot db-adb)_{a,b\in \CP }$ (where $a\cdot db$ is a product of the elements $a\in \CP$ and $db\in \O$  in the algebra $\CP \rtimes U(\O )$ but $adb\in \O$).
\end{enumerate}
\end{theorem}

{\it Proof}. 1-3. Statements 1--3 follow from (\ref{KPdif5}).

4. Statement 4 follows from statement 3 and Theorem \ref{23Jun19}.(2): If we represent the algebra $\CP$ as the factor algebra $P_\L/(f_s)_{s\in \G }$ then, in view of (\ref{KPdif1}), the algebra  $\CU (\CP ) \simeq \CP \rtimes U(\O )/R(\O )$ is generated by the algebra $\CP$ and the set $\{ \d_i\}_{i\in \L}$ that satisfies precisely the defining relations (a)--(c) in Theorem \ref{23Jun19}.(2).  $\Box $\\

{\bf The PBW Theorem holds for the Poisson enveloping algebra.} By Proposition \ref{B29Jul19}.(2), the algebra $F_{ab}$ admits a filtration $\{ F_{ab, \leq i}=\oplus_{|\alpha |\leq i}\CP \d^\alpha = \oplus_{|\alpha |\leq i} \d^\alpha\CP\}_{i\geq 0}$ by the total degree of the elements $\{ \d_{x_j}\, | \, j\in \L\}$. Its image under the epimorphism $F_{ab}\ra \CU (\CP )$ is the filtration $\{ \CU (\CP )_{\leq i} \}_{i\geq 0}$ where $\CU (\CP )_{\leq i}$ is the image of $F_{ab, \leq i}$. This is the filtration on the algebra $\CU (\CP )$ by the total degree of the elements $\{ \d_{x_j}\, | \, j\in \L\}$. The associated graded algebra 
\begin{equation}\label{grRab}
{\rm gr} F_{ab}=\bigoplus_{i\geq 0} F_{ab, \leq i}/F_{ab, \leq i-1}\simeq \CP [ \d_{x_j}]_{j\in \L}
\end{equation}
is a polynomial algebra in the variables $\{ \d_{x_j}\, | \, j\in \L\}$ with coefficients in the algebra $\CP$. Theorem \ref{A30Jul19} shows that the  PBW Theorem holds for the Poisson enveloping algebra.

\begin{theorem}\label{A30Jul19}
Let $\CP$ be a Poisson algebra and ${\rm gr} \, \CU (\CP ) =\bigoplus_{i\geq 0} \CU (\CP)_{\leq i} / \CU (\CP)_{\leq i-1}$ be the associated graded algebra. Then the algebra  ${\rm gr} \, \CU (\CP )$  is isomorphic to the symmetric algebra $\Sym_\CP (\O )$ of  the module   of K\"{a}hler differentials $\O = \OCP$ of the algebra $\CP$. 
\end{theorem}

{\it Proof}. By (\ref{KPdif5}), 
$$\CU (\CP)_{\leq i} / \CU (\CP)_{\leq i-1}=\begin{cases}
\CP & \text{if }i=0,\\
\O & \text{if }i=1.\\
\end{cases}
$$
 By (\ref{KPdif4}) and  (\ref{KPdif5}), for all $i\geq 2$, $\CU (\CP)_{\leq i} / \CU (\CP)_{\leq i-1}=F_{ab, i} / M_0F_{ab, i-1}$. Therefore, by (\ref{grRab}), 
 $$ {\rm gr}\, \CU (\CP ) = \CP [\d_{x_j}]_{j\in \L}/(M_0) \simeq \Sym_{\CP} (\oplus_{j\in \L} \CP \d_{x_j}/M_0)\simeq \Sym_\CP (\O ). \;\; \Box$$
 
 \begin{corollary}\label{B30Jul19}
Let $\CP$ be a Poisson algebra. Then the algebras $\CU (\CP ) \simeq  {\rm gr} \,\CU (\CP) \simeq \Sym_\CP (\O )$ are isomorphic as $\N$-filtered left/right $\CP$-modules. Furthermore, for all $i\geq 0$ the left (resp., right) $\CP$-modules $\CU (\CP )_{\leq i}$ and $\bigoplus_{j=0}^i\CU (\CP)_{\leq j} /\CU (\CP)_{\leq j-1}$ are isomorphic.
\end{corollary}

{\it Proof}. The corollary follows from (\ref{KPdif5}) and Theorem \ref{A30Jul19}. $\Box $\\

{\bf Criterion for the algebras $\CU (\CP )$ and ${\rm gr} \, \CU (\CP )$ to be projective left/right $\CP$-modules.}

\begin{corollary}\label{31Jul19}
Let $\CP$  be a  Poisson algebra. Then the following statement are equivalent:
\begin{enumerate}
\item The algebra $\CU (\CP )$ is a projective left (resp., right) $\CP$-module. 
\item The algebra ${\rm gr} \, \CU (\CP )\simeq \Sym (\O_\CP)$ is a  projective left (resp., right) $\CP$-module. 
\item The left $\CP$-module $\O_\CP$ is projective.
\end{enumerate}
\end{corollary}

{\it Proof}. The corollary follows from Corollary 
\ref{B30Jul19}. $\Box$

Corollary \ref{31Jul19} shows that the algebras $\CU (\CP )$ and ${\rm gr} \, \CU (\CP )$ are projective left and right $\CP$-modules provided the Poisson algebra is a regular algebra of essentially finite type.

\begin{corollary}\label{ax31Jul19}
Let $\CP$  be a regular  Poisson algebra of essentially finite type. Then
\begin{enumerate}
\item The algebra $\CU (\CP )$ is a projective left and right  $\CP$-module. Furthermore, for each $i\geq 0$, the left and right $\CP$-module $\CU (\CP)_{\leq i}$ is  projective and  finitely generated.
\item The algebra ${\rm gr} \, \CU (\CP )\simeq \Sym (\O_\CP )$ is a  projective left and right $\CP$-module. Furthermore, for each  $i\geq 0$, the left and right $\CP$-module $\CU (\CP)_{\leq i} / \CU (\CP)_{\leq i-1}$ is  projective and finitely generated.
\end{enumerate}
\end{corollary}

{\it Proof}. 2. The algebra $\CP$ is a regular algebra of essentially finite type. It is well-known that the $\CP$-module $\O = \OCP$ is  projective and finitely generated.  By Theorem \ref{A30Jul19}, ${\rm gr} \, \CU (\CP )\simeq \Sym_{\CP}(\OCP )$ and statement 2 follows.

1. Statement 1 follows from statement 2 and Corollary  \ref{B30Jul19}.  $\Box$ \\

{\bf The canonical Poisson structure on ${\rm gr}\, \CU (\CP )$. } Given an $\N$-filtered algebra $R=\bigcup_{i\geq 0} R_i$ such that the associated graded algebra ${\rm gr}\, R$ is a commutative algebra. Then the algebra ${\rm gr} \, R$ is a Poisson algebra where the Poisson structure is given by the rule: For all elements $\br_i=r_i+R_{i-1}\in R_i/R_{i-1}$ and $\br_j=r_j+R_{j-1}\in R_j/R_{j-1}$, 
$$\{ \br_i, \br_j \} := [r_i,r_j]+R_{i+j-2}.$$
Let $\CP$ be a Poisson algebra. The algebra $\CU (\CP )$ admits the filtration $\{ \CU (\CP )_{\leq i}\}_{i\geq 0}$ such that ${\rm gr} \, \CU (\CP) \simeq \Sym_\CP (\OCP )$ is a {\em commutative} algebra. So, the algebra ${\rm gr} \, \CU (\CP )$ is a Poisson algebra. This Poisson structure on ${\rm gr}\, \CU (\CP )$  is called  the {\em canonical Poisson structure}. 

\begin{proposition}\label{A8Aug19}
Let $\CP$ be a Poisson algebra. Then the algebra  ${\rm gr} \, \CU (\CP ) \simeq \Sym_\CP (\OCP )$ is  a Poisson algebra where for all elements $a,b\in \CP$, $\{ a,b\} =0$, $\{ da , b\} = \{ a,b\}$  and $\{ da , db\} =d \{ a,b\}$. 
\end{proposition}

By Proposition \ref{A8Aug19}, the algebra $\CP$ is a trivial Poisson subalgebra of  ${\rm gr} \, \CU (\CP )$ (i.e. $\{ \CP , \CP \} =0)$ but 
the Poisson structure on $\CP$ can be {\em recovered} from the Poisson structure on 
 ${\rm gr} \, \CU (\CP )$ by the rule: For all elements $a,b\in \CP$, $\{ a,b\} = \{ da, b\}$. 

\begin{proposition}\label{B8Aug19}
Let $\CP$ and $\CP'$ be  Poisson algebras. The following statements are equivalent:
\begin{enumerate}
\item  The Poisson algebras $\CP$ and $\CP'$  are isomorphic.
\item The $\N$-filtered algebras  $\CU (\CP ) $ and $\CU (\CP')$  are isomorphic  as $\N$-filtered algebras (where $\CU (\CP ) =\bigcup_{i\geq 0} \CU (\CP )_{\leq i}$ and $\CU (\CP') =\bigcup_{i\geq 0} \CU (\CP' )_{\leq i}$). 
\item The $\N$-graded Poisson algebras  ${\rm gr}\,\CU (\CP ) \simeq \Sym_{\CP} (\OCP )$ and ${\rm gr}\,\CU (\CP')\simeq\Sym_{\CP'} (\O_{\CP'} ) $  are isomorphic as $\N$-graded algebras. 
\end{enumerate}
\end{proposition}

{\it Proof}. $(1\Rightarrow 2 \Rightarrow 3)$ The implications are obvious.

$(3\Rightarrow 1)$ The implication follows from Proposition \ref{A8Aug19}. $\Box$ \\

{\bf Criterion for the algebra ${\rm gr} \, \CU (\CP )$ to be a Noetherian algebra.} As a corollary of Theorem \ref{A30Jul19}, we obtain a criterion for the algebra ${\rm gr} \, \CU (\CP )$ to be a Noetherian algebra or finitely generated.

\begin{proposition}\label{B31Jul19}
Let $\CP$  be a  Poisson algebra.
\begin{enumerate}
\item The algebra ${\rm gr} \, \CU (\CP )$ is a Noetherian algebra iff the algebra $\CP$ is  a Noetherian algebra and the $\CP$-module $\OCP$ is finitely generated.
\item  The algebra ${\rm gr} \, \CU (\CP )$ is a finitely generated algebra iff the algebra $\CP$ is  finitely generated and the $\CP$-module $\OCP$ is finitely generated.
\end{enumerate}
\end{proposition}

{\it Proof}. 1. By Theorem \ref{A30Jul19}, the algebra  ${\rm gr} \, \CU (\CP )\simeq \Sym_{\CP}(\OCP )$ is a positively graded algebra, and statement 1 follows. 

2. $(\Rightarrow )$ The algebra  ${\rm gr} \, \CU (\CP )\simeq \Sym_{\CP}(\OCP )$ is a commutative algebra. If it  finitely generated 
then necessarily the algebra  $\CP \simeq \Sym_{\CP}(\OCP )/(\OCP ) $ is also finitely generated, and the algebra  ${\rm gr} \, \CU (\CP )$ is also Noetherian. Then, by statement 1, the $\CP$-module $\OCP$ is finitely generated. 

$(\Leftarrow )$ The implication follows from ${\rm gr} \, \CU (\CP )\simeq \Sym_{\CP}(\OCP )$. $\Box $\\



{\bf The derivation algebras and the polynomial Poisson algebras.} 

{\it Definition.}  A {\bf polynomial Poisson algebra} is a Poisson algebra the  associative algebra of it is a polynomial algebra. \\

We aim to show that the PEA  of a polynomial Poisson algebra is a derivation algebra with PBW basis (Theorem \ref{A25Mar19}.(1)).\\

{\it Definition.} Let $D$ be a $K$-algebra,  $\d_1, \ldots , \d_n$ be $K$-derivations of $D$ and $C=\{ c_{ij}^k\in D\, | \, i,j,k=1, \ldots, n; i\neq j \}$. A {\bf derivation algebra} $A=D[t; \d , C]$ of rank $n$, where $t=(t_1, \ldots , t_n)$ and $\d = (\d_1, \ldots , \d_n)$,  is a $K$-algebra which is generated by the algebra $D$ and the elements $t_1, \ldots , t_n$ that satisfy the defining relations: For all $i,j=1, \ldots , n$ and $i\neq j$, 
\begin{equation}\label{DefRDA}
[t_i, d]=\d_i(d) \;\; {\rm and}\;\; [t_i, t_j]=\sum_{k=1}^nc_{ij}^kt_k.
\end{equation}
Using (\ref{DefRDA}), the (Jacobi) identity
$$ [t_i,[t_j,t_k]]=[[t_i,t_j],t_k]+[t_j,[t_i,t_k]]$$
can be rewritten as follows (we adopt the convention that the  summation is assumed over repeated lower and upper indices, i.e. $c_{ij}^k\d_k$ means $\sum_k c_{ij}^k\d_k$):
\begin{equation}\label{DefRDA1}
\d_i(c_{jk}^\l)t_\l+c_{jk}^\l c_{i\l}^\mu t_\mu = -\d_k(c_{ij}^\l ) t_\l +c_{ij}^\l c_{\l k}^\mu t_\mu +\d_j(c_{ik}^\l )t_\l + c_{ik}^\l c_{j\l }^\mu t_\mu .
\end{equation}

We say that the  derivation algebra $A$ has {\bf PBW basis (over $D$)} if $A= \bigoplus_{\alpha \in \N^n}Dt^\alpha$ is a free left $D$-module  where $t^\alpha = t_1^{\alpha_1}\cdots  t_n^{\alpha_n}$ and $\alpha =(\alpha_1, \ldots , \alpha_n)$. 
Clearly, if the algebra $A$ has PBW basis then the order of multiples in $t^\alpha$ can be arbitrary and $A=\bigoplus_{\alpha \in \N^n} t^\alpha D$ is also a free right $D$-module.  

When the algebra $D$ is a commutative algebra, Theorem \ref{25Mar19}  a criterion for the derivation algebra $A$ to have PBW basis. 

\begin{theorem}\label{25Mar19}
Let $D$ be a commutative algebra, $A=D[t; \d , C]$ be a derivation algebra of rank $n$, $\CG := \sum_{i=1}^n Dt_i$ and $\CG':= D+\CG'$. Then the following statements are equivalent:
\begin{enumerate}
\item The algebra $A$ has PBW basis. 
\item The direct sum of free left $D$-modules of rank 1, $\CG'=D\oplus\bigoplus_{i=1}^nDt_i$, is a Lie algebra (w.r.t. $[\cdot , \cdot ]$). In particular, $D$ is an abelian Lie ideal of the Lie algebra $\CG'$ and $[\d_i, \d_j]=\sum_kc_{ij}^k\d_k$ for all $i$ and $j$. 
\item The direct sum of free left $D$-modules of rank 1, $\CG= \bigoplus_{i=1}^nDt_i$, is a Lie algebra (w.r.t. $[\cdot , \cdot ]$),  $\CG'=D\oplus \CG$ and the map $\CG\ra \Der_K(D)$, $ dt_i\mapsto [dt_i, \cdot ]$ is a Lie algebra homomorphism.
\end{enumerate}
If one of the equivalent conditions holds then the algebra $A$ is isomorphic to the factor algebra $U(\CG')/I$ of the universal enveloping algebra $U(\CG')$ of the Lie algebra $\CG'$ modulo the ideal $I$ which is generated by the elements $\{ dt_i-d\cdot t_i \, | \, i=1, \ldots , n; d\in D\}$ (where $dt_i\in Dt_i\subseteq \CG'$ and $d\cdot t_i$ is a product of two elements in $U(\CG')$). 
\end{theorem}

{\it Proof}. $(1\Rightarrow 2)$ By statement 1 and the fact that $D$ is a commutative algebra,  $\CG'=D\oplus\bigoplus_{i=1}^nDt_i$ is a Lie subalgebra of the Lie algebra  $(A, [\cdot , \cdot ])$ since for all $d,d'\in D$ and $i,j=1, \ldots , n$: 
$$[dt_i, d']=d\d_i(d') \;\; {\rm and }\;\; [dt_i, d't_j]= d\d_i(d')t_j-d'\d_j(d)t_i+dd'\sum_{k=1}^n c_{ij}^kt_k.$$

$(2\Leftrightarrow 3)$ Since the algebra $D$ is commutative  statements 2 and 3 are equivalent; and $\CG'=D\rtimes \CG$ is a semi-direct product of the abelian Lie algebra $D$ and the Lie algebra $\CG$ that acts on $D$ by the Lie algebra homomorphism $\CG\ra \Der_K(D)$, $ dt_i\mapsto [dt_i, \cdot ]$. 

$(2\Rightarrow 1)$ Let $A':= U(\CG')/I$.   Then there is a natural/tautological algebra epimorphism  $A'\ra A$ (since the relations $I$ hold in the algebra $A$). The epimorphism is an isomorphism since the algebra $A'$ has the same generators and defining relations as the algebra $A$. 

Using the ordering $d <t_i<t_j$ for all $d\in D$ and $i<j$, in order to finish the proof it suffices to show that all the  ambiguities can be resolved (and hence the Diamond Lemma
 guarantees the result). In view of (\ref{DefRDA}),  there are  two sorts of ambiguities: $t_jt_id$ where $i<j$ and $d\in D$; and $t_kt_jt_i$ where $i<j<k$. 
 
$\bullet$ For all elements $d\in D$ and  $i<j$, 
 \begin{eqnarray*}
 t_j(t_id)&=& t_j(dt_i+ \d_i(d))= dt_jt_i+\d_j(d)t_i+\d_i(d)t_j +\d_j\d_i(d)\\
 &=& dt_it_j+dc_{ji}^kt_k+\d_j(d)t_i+\d_i(d)t_j +\d_j\d_i(d),\\
 (t_jt_i)d&=& (t_it_j+c_{ji}^kt_k)d=t_idt_j+t_i\d_j(d)+c_{ji}^kd t_k+c_{ji}^k\d_k(d)\\
 &=& dt_it_j +\d_i(d)t_j+\d_j(d)t_i+\d_i\d_j(d)+c_{ji}^kd t_k+c_{ji}^k\d_k(d).
 \end{eqnarray*} 
 By comparing the quadratic, linear and free terms (w.r.t. the variables $t_i$) of both equalities above, we obtained identities (by using the equality $[\d_i, \d_j]=c_{ij}^k\d_k$ and the fact that the algebra $D$ is commutative).

 $\bullet$ For all $i<j<k$,
 
\begin{eqnarray*}
 t_k(t_jt_i)&=& t_k t_it_j+t_kc_{ji}^\l t_\l =t_it_kt_j+c_{ki}^\mu t_\mu t_j+ c_{ji}^\l t_kt_\l +\d_k(c_{ji}^\l )t_\l \\
 &=& t_it_jt_k+ c_{kj}^\nu t_it_\nu + \d_i(c_{kj}^\nu ) t_\nu + \sum_{\mu\leq j} c_{ki}^\mu t_\mu t_j+ \sum_{\mu > j} c_{ki}^\mu t_j t_\mu +\sum_{\mu > j} c_{ki}^\mu c_{\mu j}^\alpha t_\alpha \\
 &+&  \sum_{k\leq \l} c_{ji}^\l t_kt_\l + \sum_{k> \l} c_{ji}^\l t_\l t_k + \sum_{k> \l} c_{ji}^\l c_{k\l }^\beta t_\beta +\d_k(c_{ji}^\l )t_\l , \\
 (t_kt_j)t_i&=& t_jt_kt_i+c_{kj}^\l t_\l t_i     \\
 &=& t_jt_it_k+t_jc_{ki}^\mu t_\mu +c_{kj}^\l t_\l t_i \\
 &=& t_it_jt_k+c_{ji}^\alpha t_\alpha t_k+ c_{ki}^\mu t_j  t_\mu +\d_j(c_{ki}^\mu )  t_\mu+c_{kj}^\l t_\l t_i \\
 &=& t_it_jt_k+\sum_{\alpha \leq k} c_{ji}^\alpha t_\alpha t_k+
 \sum_{\alpha > k} c_{ji}^\alpha t_k t_\alpha +\sum_{\alpha > k} c_{ji}^\alpha c_{\alpha k}^\o t_\o \\
 &+ & \sum_{j\leq \mu} c_{ki}^\mu t_j  t_\mu +\sum_{j> \mu} c_{ki}^\mu  t_\mu t_j+\sum_{j> \mu} c_{ki}^\mu c_{j\mu}^\beta t_\beta +
 \d_j(c_{ki}^\mu )  t_\mu \\
 &+& \sum_{\l\leq  i} c_{kj}^\l t_\l t_i + \sum_{\l >  i} c_{kj}^\l t_i t_\l  +\sum_{\l> i} c_{kj}^\l c_{\l i}^\g t_\g .\\
 \end{eqnarray*} 
  By comparing the cubic, quadratic and linear terms (w.r.t. the variables $t_i$) of both equalities above we obtain identities. The equality of quadratic terms is  straightforward
  when we use the equality
  $$c_{kj}^\nu t_it_\nu =\sum_{i\leq \nu} c_{kj}^\nu t_it_\nu+ \sum_{i>\nu}c_{kj}^\nu t_\nu t_i +\sum_{i>\nu}c_{kj}^\nu c_{i\nu}^\g t_\g .$$
  The equality of the linear terms follows from  (\ref{DefRDA1}). 
  $\Box$


\begin{theorem}\label{A25Mar19}
Let $\CP = P_n=K[x_1, \ldots , x_n]$ be a Poisson polynomial algebra and 
 $C_{P_n}=(c_{ij})$ where $c_{ij}:= \{ x_i, x_j\} \in P_n$. Then 
 \begin{enumerate}
 \item The Poisson    enveloping algebra $\CU (P_n)$ of the Poisson algebra $P_n$ is isomorphic to the derivation algebra $P_n[t; \d , C=\{ c_{ij}^k\} ]$ of rank $n$  that has PBW basis over $P_n$ where $t=(t_1, \ldots , t_n)$, $\d = (\d_1:= \{ x_1, \cdot \} , \ldots , \d_n:= \{ x_n, \cdot \})$ and $c_{ij}^k:= \frac{\der c_{ij}}{\der x_k} $. In particular, 
 $\CU (P_n) =\bigoplus_{\alpha \in \N^n} P_nt^\alpha = \bigoplus_{\alpha \in \N^n} t^\alpha P_n$ is a free left and right $P_n$-module.
  \item As an abstract algebra, the algebra $\CU (P_n)$  is generated over $K$ by the elements $x_1, \ldots , x_n, t_1, \ldots , t_n$  subject the defining relations:
\begin{equation}\label{DefRUPAC}
[x_i,x_j]=0, \;\; [t_i, x_j]=c_{ij}\;\; {\rm and}\;\;  [t_i, t_j]=\sum_{i=1}^n  \frac{\der c_{ij}}{\der x_k}t_k.
\end{equation}
 \item $\CG'=P_n\oplus\bigoplus_{i=1}^nP_nt_i$ is a Lie subalgebra of $(\CU (P_n), [\cdot , \cdot ])$ and $\CG'\simeq P_n\rtimes \CG$  is a semidirect product of the Lie algebras $P_n$ and $\CG = \bigoplus_{i=1}^nP_nt_i$.  
 \item The algebra $\CU (P_n)$ is isomorphic to the factor algebra $  U(\CG')/I$ of the universal enveloping algebra $U(\CG')$ of the Lie algebra $\CG'$ modulo the ideal $I$ which is generated by the elements $\{ pt_i-p\cdot t_i \, | \, i=1, \ldots , n; p\in P_n\}$ (where $pt_i\in P_nt_i\subseteq \CG'$,  and $p\cdot t_i$ is a product of two elements in $U(\CG')$). 

\end{enumerate}

\end{theorem}

{\it Proof}. 2. Statement 2 follows from Theorem \ref{23Jun19}.(2).

1. Statement 1 follows from Theorem \ref{30Jul19}.(4), Theorem \ref{A30Jul19} and the fact that the module of K\"{a}hler differentials of the polynomial algebra $P_n$,  $\O_{P_n}=\bigoplus_{i=1}^n P_n t_i$, is a free $P_n$-module where $t_i=dx_i$. 

3 and 4. Statements 3 and 4 follow from statement 1 and Theorem \ref{25Mar19}. $\Box$

\begin{theorem}\label{A11Jab20}
Suppose that  $\CP = S^{-1}K[x_i]_{i\in \L} / (f_s)_{s\in \G}$ where $S$ is a multiplicative subset  of the polynomial algebra $K[x_i]_{i\in \L}$ ($\L$ and $\G$ are index sets) and the algebra  $\CP$ is a Poisson  algebra such that the module of K\"{a}hler differentials $\O_\CP = \bigoplus_{i=1}^n \CP dy_i$ is a free left $\CP$-module.  Then 
 \begin{enumerate}
 \item The Poisson    enveloping algebra $\CU (\CP )$ of the Poisson algebra $\CP $ is isomorphic to the derivation algebra $\CP [t; \d , C=\{ c_{ij}^k\} ]$ of rank $n$  that has PBW basis over $\CP$ where $t=(dy_1, \ldots , dy_n)$, $\d = (\d_1:= \{ y_1, \cdot \} , \ldots , \d_n:= \{ y_n, \cdot \})$ and $d \{ y_i,  y_j\} =\sum_{k=1}^n c_{ij}^kdy_k $. In particular, 
 $\CU (\CP) =\bigoplus_{\alpha \in \N^n} \CP t^\alpha = \bigoplus_{\alpha \in \N^n} t^\alpha \CP$ is a free left and right $\CP$-module.
  \item As an abstract algebra, the algebra $\CU (\CP )$  is generated over $K$ by the algebra $\CP$ and the  elements $t_1, \ldots , t_n$  subject the defining relations:
\begin{equation}\label{DefRUPAC}
 [t_i, x_j]=\{ y_i, x_j\} \;\; {\rm and}\;\;  [t_i, t_j]=\sum_{i=1}^n  c_{ij}^k t_k.
\end{equation}
 \item $\CG'=\CP \oplus\bigoplus_{i=1}^n\CP t_i$ is a Lie subalgebra of $(\CU (\CP ), [\cdot , \cdot ])$ and $\CG'\simeq \CP \rtimes \CG$  is a semidirect product of the Lie algebras $\CP$ and $\CG = \bigoplus_{i=1}^n\CP t_i$.  
 \item The algebra $\CU (\CP )$ is isomorphic to the factor algebra $  U(\CG')/I$ of the universal enveloping algebra $U(\CG')$ of the Lie algebra $\CG'$ modulo the ideal $I$ which is generated by the elements $\{ pt_i-p\cdot t_i \, | \, i=1, \ldots , n; p\in \CP \}$ (where $pt_i\in \CP t_i\subseteq \CG'$,  and $p\cdot t_i$ is a product of two elements in $U(\CG')$). 

\end{enumerate}

\end{theorem}

{\it Proof}. 1. Statement 1 follows from Theorem \ref{30Jul19}.(4), Theorem \ref{A30Jul19} and the fact that the module of K\"{a}hler differentials of the polynomial algebra $P_n$,  $\O_{P_n}=\bigoplus_{i=1}^n P_n t_i$, is a free $P_n$-module where $t_i=dx_i$. 

 2. Statement 2 follows from statement 1 and Theorem \ref{23Jun19}.(2).

3 and 4. Statements 3 and 4 follow from statement 1 and Theorem \ref{25Mar19}. $\Box$\\

{\it Remark.} Theorem \ref{25Mar19}, Theorem \ref{A25Mar19} and 
Theorem \ref{A11Jab20} are true (with the same proofs) for $n$ of {\em arbitrary} cardinality (not necessarily finite).


\section{Criterion for the algebra  $\CU (\CA )$ to be a domain where a Poisson algebra $\CA$  is a domain of essentially finite type}\label{CUADOM}

The aim of this section is  to prove Theorem \ref{A29Jul19}.(3) which is a criterion  for the algebra  $\CU (\CA )$ to be a domain where a Poisson algebra $\CA$  is a domain of essentially finite type. As a corollary we obtain Theorem \ref{29Jul19} which states that the algebra $\CU (\CA )$ is a domain if, in addition, the algebra $\CA$ is regular. Using a result of Huneke (Theorem  \ref{HunekeT1.1})  we obtain another criterion for the algebra  $\CU (\CA )$ to be a domain (Theorem  \ref{A31Jul19}) provided the algebra $\CA$ satisfies Serre's Condition $S_m$. 

A localization of an affine commutative algebra is called an {\bf algebra of essentially finite type}. 
In this section the following notation will remain  fixed if
it is not stated otherwise: $P_n=K[x_1, \ldots , x_n]$ is a polynomial algebra over {\bf perfect} field 
$K$, $\der_1:=\frac{\der}{\der x_1}, \ldots ,
\der_n:=\frac{\der}{\der x_n}\in \Der_K(P_n)$, $I=(f_1, \ldots , f_m)$ is a  prime but not a maximal ideal of $P_n$, $ \CA = S^{-1}(P_n/I)$ is a domain of essentially finite type and 
 $Q=Q(\CA )$ is its field of fractions,  $r= r\Big(\frac{\der  f_i}{\der x_j}\Big)$ is the rank (over $Q$) of the {\bf Jacobian matrix} $\Big(\frac{\der  f_i}{\der x_j}\Big)$ of $\CA$, $\ga_r$ is the {\bf Jacobian ideal} of the algebra $\CA$ which is
(by definition) generated by all the $r\times r$ minors of the
Jacobian matrix  $\CA$ (the algebra $\CA$ is {\em regular} iff $\ga_r=\CA$, it is the
{\bf Jacobian criterion of regularity}, \cite[Corollary 16.20]{Eisenbook}),
$\Omega_\CA$ is the module of {\bf K\"{a}hler} differentials for the
algebra $\CA$. 

For $\i =(i_1, \ldots , i_r)$ such that $1\leq
i_1<\cdots <i_r\leq m$ and $\j =(j_1, \ldots , j_r)$ such that
$1\leq j_1<\cdots <j_r\leq n$, $\D
 (\i , \j )$ denotes the corresponding minor of the Jacobian matrix of $\CA$,  and the $\i$ (resp., $\j $) is called {\bf
non-singular} if $\D (\i , \j')\neq 0$ (resp., $\D (\i', \j )\neq
0$) for some $\j'$ (resp., $\i'$). We denote by $\mI_r$ (resp.,
$\mJ_r$) the set of all the non-singular $r$-tuples $\i$ (resp.,
$\j $).

Since $r$ is the rank of the Jacobian matrix of $\CA$, it is easy to
show that $\D (\i , \j )\neq 0$ iff $\i\in \mI_r$ and $\j\in
\mJ_r$, \cite[Lemma 2.1]{gendifreg}. We denote by $\mJ_{r+1}$ the set of all
$(r+1)$-tuples $\j =(j_1, \ldots , j_{r+1})$ such that $1\leq
j_1<\cdots <j_{r+1}\leq n$ and when deleting  some element, say
$j_\nu$, we have a non-singular $r$-tuple $(j_1, \ldots
,\widehat{j_\nu},\ldots , j_{r+1})\in \mJ_r$ where the hat over
 a symbol means that the symbol is omitted from
the list. The set $\mJ_{r+1}$ is called the {\bf critical set} and
any element of it is called a {\bf critical singular}
$(r+1)$-tuple. $\Der_K(\CA)$ is the $\CA$-module of $K$-derivations of
the algebra $\CA$.  The action of a derivation $\d$ on an element $a$
is denoted by $\d (a)$. 
 
 The next theorem gives a finite set of generators and a
finite set of  defining relations for the left $\CA$-module
$\Der_K(\CA)$ when $\CA$ is a  regular algebra.

\begin{theorem}\label{9bFeb05}
(\cite[Theorem 4.2]{gendifreg} if char$(K)=0$; \cite[Theorem 1.1]{gendifregcharp} if char$(K)=p>0$;) Let the algebra $\CA$ be a regular domain of essentially finite type over the perfect field $K$. Then the left $\CA$-module
$\Der_K(\CA)$ is generated by the derivations $\der_{\i , \j }$, $\i
\in \mI_r$, $\j \in \mJ_{r+1}$ where
\begin{eqnarray*}
  \der_{\i , \j }=\der_{i_1, \ldots , i_r; j_1, \ldots , j_{r+1}}:= \det
 \begin{pmatrix}
  \frac{\der f_{i_1}}{\der x_{j_1}} & \cdots &  \frac{\der f_{i_1}}{\der
  x_{j_{r+1}}}\\
  \vdots & \vdots & \vdots \\
  \frac{\der f_{i_r}}{\der x_{j_1}} & \cdots &  \frac{\der f_{i_r}}{\der
  x_{j_{r+1}}}\\
  \der_{j_1}& \cdots & \der_{j_{r+1}}\\
\end{pmatrix}
\end{eqnarray*}
that satisfy the following defining relations (as a left
$\CA$-module): 
\begin{equation}\label{Derel}
\D (\i , \j )\der_{\i', \j'}=\sum_{l=1}^s(-1)^{r+1+\nu_l}\D (\i';
j_1', \ldots , \widehat{j_{\nu_l}'}, \ldots , j_{r+1}')\der_{\i;\j
, j_{\nu_l}'}
\end{equation}
for all $\i, \i'\in \mI_r$, $\j=(j_1, \ldots , j_r)\in \mJ_r$, and
$\j'=(j_1', \ldots , j_{r+1}')\in \JJ_{r+1}$ where $\{ j_{\nu_1}',
\ldots , j_{\nu_s}'\}=\{ j_1', \ldots , j_{r+1}'\}\backslash \{
j_1, \ldots , j_r\}$ and $\der_j=\frac{\der}{\der x_j}$.
\end{theorem}

{\it Remark.} Suppose that a Poisson algebra $\CA= S^{-1} (P_n/I)$ is of essentially finite type. Then, if necessary, we may assume that $\{x_i, x_j\} \in P_n/I$ for all $i,j=1, \ldots , n$ and the multiplicative set $S$ consist of {\em regular} elements of the algebra $P_n/I$. The second statement follows from the fact that $S^{-1}(P_n/I)\simeq \bS^{-1}(P_n/J)$ where $\ass_{P_n/I}(S)=J/I$ for some ideal $J$ of $P_n/I$ and $\bS = \{ s+J/I\, | \, s\in S\}$. Let $s\in S$ be the product of all the denominators, say $s_{ij}$,  in $\{ x_i, x_j\}= s_{ij}^{-1}a_{ij}$ where $a_{ij}\in P_n/I$. Then $P_n/I\simeq P_n[x_{n+1}]/(I, sx_{n+1}-1)$ and we are done since $s_{ij}^{-1} = s^{-1}\cdot \frac{s}{s_{ij}}$.\\

{\bf The Gelfand-Kirillov dimension  of the algebra $\CU (\CA )$.} Proposition \ref{GKA29Jul19} gives a lower bound for the Gelfand-Kirillov dimension of the algebras $\CU (\CA )$, ${\rm gr} \, \CU (\CA )$ and $\Sym_{\CA} ( \CA_\CA )$. In some important cases the lower bound is the Gelfand-Kirillov of these algebras (Theorem \ref{A29Jul19}.(4) and Theorem \ref{29Jul19}).  

\begin{proposition}\label{GKA29Jul19}
Let a Poisson algebra $\CA = S^{-1} (P_n/I)$ be a domain of essentially finite type where  $I=(f_1, \ldots , f_m)$ is a prime but not maximal  ideal of $P_n$ and $r=r\Big(\frac{\der f_i}{\der x_j}\Big)$ is the rank of the Jacobian matrix $\Big(\frac{\der f_i}{\der x_j}\Big)$ over the field of fractions of the domain $P_n/I$. Then $\GK \, \CU (\CA ) \geq \GK \, {\rm gr} \, \CU (\CA ) = \GK \, \Sym_{\CA} ( \O_\CA )\geq 2\GK (\CA ) = 2(n-r)$. 
\end{proposition}

{\it Proof}. By \cite[Lemma 8.3.20]{MR}, $\GK \CU (\CA ) \geq \GK \, {\rm gr} \, \CU (\CA )$. By Theorem \ref{A30Jul19}, ${\rm gr} \, \CU (\CA )\simeq \Sym_{\CA} ( \O_\CA )$, and so $\GK {\rm gr} \, \CU (\CA ) = \GK \, \Sym_{\CA} ( \O_\CA )$. By Theorem \ref{A29Jul19}.(1), $\GK \, \Sym_{\CA} ( \O_\CA )\geq \GK (\CA_{\D (\i , \j )})+n-r=\GK (\CA ) +n-r=2(n-r)= 2\GK (\CA )$.  $\Box $\\

{\bf Criterion for the algebra $\CU (\CA )$ to be a domain.} Theorem \ref{A29Jul19}.(3) is a criterion for the algebra $\CU (\CA )$ to be a domain where the Poisson algebra $\CA$ is a domain of essentially finite type.

\begin{theorem}\label{A29Jul19}
Let a Poisson algebra $\CA = S^{-1} (P_n/I)$ be a domain of essentially finite type where  $I=(f_1, \ldots , f_m)$ is a prime but not maximal  ideal of $P_n$ and $r=r\Big(\frac{\der f_i}{\der x_j}\Big)$ is the rank of the Jacobian matrix $\Big(\frac{\der f_i}{\der x_j}\Big)$ over the field of fractions of the domain $P_n/I$. Then 
\begin{enumerate}
\item For each $\i \in \mI_r$ and $\j \in \mJ_r$, the algebra $\CU (\CA )_{\D (\i , \j )}\simeq \CU (\CA_{\D (\i , \j )} )\simeq \bigoplus_{\alpha \in \N^{n-r}}\CA_{\D (\i , \j )}\d^{\alpha_1}_{x_{j_{r+1}}}\cdots \d^{\alpha_{n-r}}_{x_{j_n}}$ is a Noetherian domain with $\GK \, \CU (\CA_{\D (\i , \j )} )=2\GK (\CA ) = 2(n-r)$ where $\j = (j_1, \ldots , j_r)$ and $\{ j_{r+1}, \ldots , j_n) = \{ 1, \ldots , n\} \backslash \{ j_1, \ldots , j_r\}$. 
\item The kernel of the algebra homomorphism 
$$\th : \CU (\CA ) \ra \prod_{\i \in \mI_r , \j\in \mJ_r} \CU (\CA )_{\D (\i , \j )}, \;\; u\mapsto (\ldots , \frac{u}{1}, \ldots )$$
is a finitely generated left and right ideal of the algebra $\CU (\CA )$ which is  $S_{\D (\i , \j )}$-torsion for all $\i\in \mI_r$ and $\j \in \mJ_r$. In particular, $\D (\i , \j )^k\ker (\th ) =0$ and $\ker (\th )\D (\i , \j )^l =0$ for some natural numbers $k,l\geq 0$. Furthermore, $\ga_r^s\ker (\th ) =0$ and 
$\ker (\th )\ga_r^t =0$ for some natural numbers $s,t\geq 0$. 
\item The following statements are equivalent:
\begin{enumerate}
\item The algebra $\CU (\CA)$ is a domain.
\item The algebra ${\rm gr} \, \CU (\CA)$ is a domain.
\item The algebra $\Sym_{\CA} (\O_\CA )$ is a domain.
\item The elements $\{ \D (\i , \j )\, | \, \i \in \mI_r, \j\in \mJ_r\}$ are regular in $\CU (\CA )$. 
\item The element $\D (\i , \j )$ is a regular element of the algebra $\CU ( \CA )$ for some $\i \in \mI_r$ and $\j \in \mJ_r$. 
\item $\lann_{\CU (A)}(\ga_r) = \rann_{\CU (A)}(\ga_r)=0$. 
\end{enumerate}
\item Suppose that the algebra $\CU (\CA )$ is a domain. Then $\GK \, \CU (\CA )= \GK \, {\rm gr} \, \CU (\CA )=\GK \, \Sym_{\CA} ( \O_\CA )=2\GK (\CA )=2(n-r)$.
\end{enumerate}
\end{theorem}

{\it Proof}. 1. Statement 1 follows from Theorem \ref{A11Jab20} since $\O_{\D (\i , \j )}=\bigoplus_{\nu = r+1}^n\CA_{\D (\i , \j )}\d_{x_{i_\nu}}$.


2. Since the algebra $\CU (\CA )$ is a Noetherian algebra statement 2 follows from statement 1 and the fact that the sets $\mI_r$ and $\mJ_r$ are finite sets. 

3. $(b\Rightarrow a\Rightarrow f)$ The implications are obvious. 

$(b\Leftrightarrow c)$ By Theorem \ref{30Jul19}, ${\rm  gr} \, \CU (\CA ) \simeq \Sym_{\CA} (\O_\CA )$, and the equivalence follows.

$(f\Rightarrow e)$ The implication follows from the fact that the ideal $\ga_r$ is generated by the elements $\D ( \i , \j )$ where $ \i \in \mI_r$ and $\j \in \mJ_r$. 

$(e\Rightarrow d)$ The implication is obvious (see statement 2).

$(d\Rightarrow b)$ By Corollary \ref{B30Jul19}, the left and right $\CA$-modules $\CU (\CA )$ and ${\rm gr} \, \CU (\CA )$ are isomorphic. Therefore, the algebra ${\rm gr} \, \CU (\CA )$ is $S_\D$-torsion  free where $\D = \D (\i , \j )$. In particular, the algebra ${\rm gr} \, \CU (\CA )$ is a subalgebra of its localization at the powers of the element $\D$, 
$$ {\rm gr} \, \CU (\CA )\subseteq ({\rm gr} \, \CU (\CA ))_\D \simeq {\rm gr} \, \CU (\CA )_\D\simeq {\rm gr} \, \CU (\CA_\D )\simeq \CA_\D [ \d_{x_{j_{r+1}}}, \ldots , \d_{x_{j_n}}],$$
a polynomial algebra over $\CA_\D$ in the variables $\d_{x_{j_{r+1}}}, \ldots , \d_{x_{j_n}}$ which is obviously a domain, and so is its subalgebra ${\rm gr} \, \CU (\CA )$.

4. Since the algebra $\CU (\CA )$ is a domain, $\CU (\CA) \subseteq \CU (\CA)_\D$
 where $\D = \D (\i , \j )$, $\i \in \mI_r$ and $\j \in \mJ_r$. By Theorem \ref{30Jul19}, $\CU (\CA )_\D\simeq \CU (\CA_\D )$. Now, 
 by statement 1, $\GK \, \CU (\CA ) \leq \GK \, \CU (\CA_\D )=2\GK (\CA )$. Now, statement 4 follows from the equalities $\GK (\CA )=\GK (\CA_\D )=n-r$.  $\Box $\\
 
 {\bf The algebra $\CU (\CA )$ is a domain when $\CA$ is a regular domain of essentially finite type.} As a corollary of Theorem \ref{A29Jul19}.(3), we obtain Theorem \ref{29Jul19} that states that the algebra $\CU (\CA )$ is a domain provided the algebra $\CA$ is a  regular domain  of essentially finite type. 

\begin{theorem}\label{29Jul19}
Let a Poisson algebra $\CA = S^{-1} (P_n/I)$ be a regular domain of essentially finite type where  $I=(f_1, \ldots , f_m)$ is a prime but not maximal  ideal of $P_n$ and $r=r\Big(\frac{\der f_i}{\der x_j}\Big)$ is the rank of the Jacobian matrix $\Big(\frac{\der f_i}{\der x_j}\Big)$ over the field of fractions of the domain $P_n/I$.  Then the algebra $\CU (\CA )$ is a Noetherian domain with $\GK \, \CU (\CA ) = \GK \, {\rm gr}\, \CU (\CA ) = \GK \, \Sym_\CA (\O_\CA ) = 2\GK (\CA ) = 2(n-r)$. 
\end{theorem}

{\it Proof}. By Proposition \ref{B29Jul19}.(4), the algebra $\CU (\CA )$ is Noetherian. The algebra $\CA$ is regular, and so $\ga_r = \CA$. By Proposition \ref{A29Jul19}.(3), the algebra $\CU (\CA )$ is a domain. Now, the theorem follows from Theorem \ref{A29Jul19}.(4).  $\Box $\\

{\bf The $\CA$-torsion submodules of the algebras $\CU (\CP )$, ${\rm gr}\, \CU (\CA )$ and $\Sym (\O_\CA )$.} Theorem \ref{BA2Aug19} describes the $\CA$-torsion submodules of the three algebras above. 

\begin{theorem}\label{BA2Aug19}
Let a Poisson algebra $\CA = S^{-1}(P_n/I)$ be a domain of  essentially finite type over the field $K$  and $\CV :={\rm gr}\, \CU (\CA ) \simeq \Sym (\O_\CA )$.   Then 
\begin{enumerate}
\item $\tor_{\CA\backslash \{ 0\} }(\CU (\CA ))=\tor_{S_{\D (\i , \j )}}(\CU (\CA ))=\lann_{\CU (\CA )}(\ga_r^i)=\rann_{\CU (\CA )}(\ga_r^i)=\lann_{\CU (\CA )}(\D (\i , \j )^i)=\rann_{\CU (\CA )}(\D (\i , \j )^i)$ for all $i\gg 1$, $\i  \in \mI_r$ and $\j \in \mJ_r$.
\item $\tor_{\CA\backslash \{ 0\} }(\CV )=\tor_{S_{\D (\i , \j )}}(\CV )=\lann_{\CV }(\ga_r^i)=\rann_{\CV }(\ga_r^i)=\lann_{\CV }(\D (\i , \j )^i)=\rann_{\CV}(\D (\i , \j )^i)$ for all $i\gg 1$, $\i  \in \mI_r$ and $\j \in \mJ_r$.
\end{enumerate}
\end{theorem}

{\it Proof}. By Proposition \ref{B30Jul19}, the left and right $\CA$-modules $\CU (\CA)$ and $\CV$ are isomorphic. So, statements 1 and 2 are equivalent. So, it suffices to prove, say, statement 1.
 The algebras $\CU (\CP )$, ${\rm gr}\, \CU (\CA )$ and $\Sym (\O_\CA )$ are Noetherian finitely generated algebras. The multiplicative sets $\CA\backslash \{ 0\}$ and $S_\D :=\{ \D^i\, | \, i\geq 0\}$ are denominator sets of them where $\D = \D (\i , \j )$, $\i \in \mI_r$ and $\j \in \mJ_r$. By Theorem \ref{A29Jul19}.(1), the algebra $\CU (\CA_\D )\simeq \CU (\CA )_\D$ is a Noetherian domain. Hence, $\tor_{\CA\backslash \{ 0\} }(\CU (\CA ))=\tor_{S_{\D (\i , \j )}}(\CU (\CA ))$ for all $\i\in \mI_r$ and $\j\in \mJ_r$. The other equalities in statement 1 follow at once from the facts that the algebra $\CU (\CA )$ is Noetherian, the Jacobian ideal $\ga_r$ is generated by the finite set $\{ \D(\i , \j )\, | \, \i \in \mI_r, \j\in \mJ_r\}$ and the sets $\CA\backslash \{ 0\}$ and $S_\D$ are (left and right) denominator sets of $\CU (\CA )$. $\Box$ \\

 {\bf The Gelfand-Kirillov dimension of the algebras $\CU (\CA )$,    ${\rm gr}\, \CU (\CA )$ and $ \Sym_\CA (\O_\CA )$ where $\CA$ is a domain  of essentially finite type.}  Theorem \ref{AAA29Jul19} shows that the regularity condition in Theorem \ref{29Jul19} can be dropped but the result about the Gelfand-Kirillov dimension holds.\\

{\bf Proof of Theorem \ref{AAA29Jul19}}. It suffices to show that the theorem holds  for a finitely generated algebra $\CA $. 

We prove the theorem by induction on the Gelfand-Kirillov dimension of the algebra $\CA$. If $\GK (\CA )=0$, i.e. the algebra $\CA$ is a finite field extension of the perfect field $K$. In particular, the field $\CA$ is perfect. Hence, $\O_\CA =0$, and the statement is obvious (Theorem \ref{30Jul19}.(4)).

Suppose that $N=\GK (\CA )\geq 1$ and the statement is true for all 
algebras $\CA'$ with $\GK (\CA')<N$. By Proposition \ref{GKA29Jul19}, it suffices to show that $\GK (\CU (\CA ))\leq 2N$.

 (i) {\em The algebra $\CU (\CA )$ is an almost (left and right) centralizing extension of $\CA$}: The statement (i) follows from Theorem \ref{23Jun19}.(2).
 
 (ii) {\em The algebra $\CU (\CA )$ is a somewhat commutative algebra}: The statement (ii) follows from the statement (i) and 
 \cite[Proposition 8.6.9]{MR}.
 
(iii)  {\em The Gelfand-Kirillov dimension over $\CU (\CA)$ is an exact function}: The statement (iii) follows from the statement (ii) and \cite[Corollary 8.4.9(i)]{MR}.

Fix $\D = \D (\i, \j )$. The Poisson algebra $\CA_\D$ is a regular finitely generated domain. Let $\CU'$ and $J$ be the image and the kernel of the algebra $\CU (\CA )$ under the localization map $\CU (\CA ) \ra \CU (\CA_\D )\simeq \CU (\CA )_\D$, $u\mapsto \frac{u}{1}$. 

(iv) $\GK (\CU')\leq 2N$: By Theorem \ref{29Jul19},  
$\GK (\CU')\leq \GK (\CU (\CA_\D ))=2\GK (\CA_\D ) =2N$. 

(v) $\GK (J)<2N$: The algebra $\CU (\CA )$ is a finitely generated Noetherian algebra. Let $\min (\ga_r)$ be the set of minimal primes over the Jacobian ideal $\ga_r$ of the algebra $\CA$. We assume that $\ga_r \neq \CA$,  otherwise the result follows from Theorem 
 \ref{29Jul19}. 
 By Theorem  \ref{BA2Aug19} and the statement (iii), there is a natural number $i\geq 1$ such that 
\begin{eqnarray*}
 \GK (J)&\leq & \GK (\CU (\CA)/(\ga_r^i))=\max \{ \GK (\CU (\CA)/(\gp ))\, | \, \gp\in \min (\ga_r) \}\\
 &=&\max \{ \GK (\CU (\CA /\gp ))\, | \, \gp\in \min (\ga_r) \}\\
 &=&\max \{ 2\GK (\CA)/(\gp ))\, | \, \gp\in \min (\ga_r) \}<2N.
\end{eqnarray*}
(vi) $\GK (\CU (\CA))\leq 2N$: By the statement (iii), 
$$\GK (\CU ( \CA )) =\max \{ \GK (\CU'), \GK (J)\},$$
 and the statement (vi) follows from the statements (iv) and (v). 
 $\Box $\\
  
Let $B$ be a matrix with coefficients in a commutative ring $R$. By $I_s(B)$ we denote the ideal generated by the $s\times s$ minors of $B$. For an $R$-module $M$, let $v(M, R)$ be the least number of its generators. Recall that a commutative ring $R$ satisfies {\em Serre's property} $S_k$ if
$$ {\rm depth} (R_\gp )\geq \min \{ \dim (R_\gp ), k\}$$
for all prime ideals $\gp$ of $R$. If $I$ is an ideal of $R$, we denote by ${\rm grade} (I)$ the length of a maximal $R$-sequence in $I$.

\begin{theorem}\label{HunekeT1.1}
(\cite[Theorem 1.1]{Huneke81}) Let $R$ be a universally catenarian Noetherian ring satisfying Serre's condition $S_m$ and let $M$ be an $R$-module having a finite free resolution, 
$$ 0\ra R^m\stackrel{B}{\longrightarrow}R^n\ra M\ra 0, \;\; B= (a_{ij}).$$
The following statements are equivalent:
\begin{enumerate}
\item $\Sym_R(M)$ is a domain.
\item ${\rm grade} (I_t(B))\geq m+2-t$ for $1\leq t \leq m$. 
\item $v(M_\gp , R_\gp ) \leq n-m+{\rm grade} (\gp ) -1$ for all nonzero primes $\gp$ of $R$.
\end{enumerate}
If any (and hence all) of the above conditions hold then $\Sym_R(M)$ is a complete intersection in $R[T_1, \ldots , T_n]$. In particular, if $R$ is Cohen-Macauley (resp., Gorenstein) then so is $\Sym_R(M)$. 
\end{theorem}

{\bf Proof of Theorem \ref{A31Jul19}}. $(1\Leftrightarrow 2\Leftrightarrow 3)$ Theorem \ref{A29Jul19}.(3).

$(3\Leftrightarrow 4\Leftrightarrow 5)$ Theorem \ref{HunekeT1.1}. $\Box$


\section{Criterion for $\Der_K(\CA ) = \CA \CH_\CA$ where $\CA$ is a regular domain of essentially finite type}\label{DERPHP}
In this section, $K$ is a field of  {\bf characteristic zero} (if it is not stated otherwise) and  we keep the assumptions and the notations of Section \ref{CUADOM}.  
 The aim of this section is to give a criterion for $\Der_K(\CA ) = \CA \CH_\CA$ where the Poisson algebra $\CA$ is a regular domain of essentially finite type (Theorem \ref{22Jul19}). When $\CA = P_n$ is a polynomial algebra in $n$ variables the criterion looks particularly nice (Corollary \ref{p22Jul19}). Examples are considered.  Lemma \ref{a24Jul19} and  Corollary \ref{b24Jul19} are regularity and symplecticity criteria for certain generalized Weyl Poisson algebras.

\begin{proposition}\label{A23Jul19}
Let a Poisson algebra $\CA$ be a regular domain of essentially finite type, $Q= Q(\CA )$ be its field of fractions, $\CC_\CA = (c_{ij})\in \M_n (\CA )$ where $c_{ij} = \{ x_i, x_j\}$ and $r(\CC_\CA)$ is the rank of the matrix $\CC_\CA\in M_n(Q)$, $r=r(\Big(\frac{\der f_s}{\der x_i}\Big))$ be the rank of the Jacobian matrix $\Big(\frac{\der f_s}{\der x_i}\Big)$ of $\CA$. Then 
\begin{enumerate}
\item $r(\CC_\CA )=\dim_Q(Q\CH_\CA )$. 
\item For each $\j =(j_1, \ldots , j_r)\in \mJ_r$, $\dim_Q(Q\CH_\CA ) = r(\CC_{\CA , \j })$ where $\CC_{\CA , \j }$ is an $n\times (n-r)$ matrix which is obtained from the matrix $\CC_\CA$ by deleting columns with indices $j_1, \ldots , j_r$, i.e. $\CC_{\CA , \j }=(c_{ij})$ where $i=1, \ldots ,n$  and $j\in \{ j_{r+1}, \ldots , j_n\} = \{ 1, \ldots , n\} \backslash \{ j_1, \ldots , j_r\}$. In particular, $r(\CC_\CA ) \leq \dim_Q(\Der_K (Q))=n-r$. 
\end{enumerate}
\end{proposition}

{\it Proof}.  By  \cite[Corollary 2.15]{gendifreg}, each derivation $\d \in \Der_K(\CA )$ is uniquely determined by its action on the elements $x_{j_{r+1}}, \ldots , x_{j_n}$. Since for each $i=1, \ldots , n$, $\pad_{x_i} = \{ x_i, \cdot \}= \sum_{j=1}^n\{ x_i, x_j\} \der_j$, the derivation $\pad_{x_i}$ is uniquely determined by the elements $\{x_i,x_{j_{r+1}}\}, \ldots , \{ x_i, x_{j_n}\}$ of $\CA$, and statements 1 and 2 follow.  $\Box $\\

{\bf The ideals $\gc_{\CA , s}$, $s=1,\ldots , r(\CC_\CA )$, and the Poisson ideal $\gc_{\CA , 1}$.} For each integer $s=1, \ldots ,  r(\CC_\CA )$, let 
  $\gc_{\CA , s}$ be the ideal of the algebra $\CA$ generated by all the $s\times s$ minors of the matrix $\CC_\CA$. By  (\ref{CPMatCoef1}), the ideals  $\gc_{\CA , s}$ do not depend on the set of essential generators $x_1, \ldots , x_n$ of the algebra $\CA$. Clearly,
  $$ \gc_{\CA , 1} \supseteq \gc_{\CA , 2} \supseteq\cdots \supseteq\gc_{\CA ,  r(\CC_\CA )}.$$

\begin{lemma}\label{a23Jul19}

\begin{enumerate}
\item The ideal $\gc_{\CA , 1}$ is a Poisson ideal of $\CA$ which is $P\Der_K(\CA )$-invariant. 
\item For all $s=1,\ldots , r(\CC_\CA )$ and $\d \in \PDer_K(\CA )$, $\d (\gc_{\CA , i})\subseteq \gc_{\CA , 1}\gc_{\CA , i-1}\subseteq \gc_{\CA , i-1}$. In particular, $\{\CA , \gc_{\CA , i}\} \subseteq \gc_{\CA , i-1}$.
\end{enumerate}
\end{lemma}

{\it Proof}. 1. For all $a\in \CA$ and $i=1, \ldots , n$, $\{ x_i, a\} = \sum_{j=1}^n \frac{\der a}{\der x_j}\{ x_i, x_j\}$, and statement 1 follows.

2. Statement 2 is obvious. $\Box $\\

{\bf Criterion for $\Der_K(\CA ) = \CA \CH_\CA$.}  For each Poisson algebra $\CP$, $\Der_K(\CP)\supseteq \CP \CH_\CP$. A Poisson algebra $\CP$ which is a regular affine domain  is called a {\em synmplectic algebra} if $\Der_K(\CP)= \CP \CH_\CP$. Theorem \ref{22Jul19} is a criterion for $\Der_K(\CA ) = \CA \CH_\CA$ where $\CA$ is a regular domain of essentially finite type.\\

{\bf  Proof of Theorem \ref{22Jul19}}.   
 $(1\Rightarrow 3)$ Given elements $\i \in \mI_r$ and $ \j = (j_1, \ldots , j_r) \in \mJ_r$. Let $\D = \D ( \i , \j )$  and $ \{ j_{r+1}, \ldots , j_n\} = \{ 1, \ldots , n\} \backslash \{ j_1, \ldots , j_r\}$. We assume that $ j_{r+1}< \cdots < j_n$.  By  \cite[Corollary 3.5]{gendifreg}, the derivations
$$\D^{-1} \der_{\i ; \j , j_{r+1}}, \ldots , \D^{-1} \der_{\i ; \j , j_n}$$ are the partial derivatives $\der_{j_{r+1}}=\frac{\der }{\der x_{j_{r+1}}}, \ldots , \der_{j_n}=\frac{\der }{\der x_{j_n}}$  of the localization $\CA_\D$ of the algebra $\CA$ at the powers of the element $\D$, respectively. Notice that $\der_{\i ; \j , j_{r+1}}, \ldots ,  \der_{\i ; \j , j_n}\in \Der_K(\CA )$ and $\Der_K(\CA ) = \CA \CH_\CA$ (by the assumption). Therefore, 
$$\begin{pmatrix}
\der_{\i ; \j , j_{r+1}}\\  \vdots \\
\der_{\i ; \j , j_n} \end{pmatrix}=\L \begin{pmatrix}
\pad_{x_1}\\  \vdots \\
\pad_{x_n} \end{pmatrix}= \L\CC_\CA \begin{pmatrix}
\der_1\\  \vdots \\
\der_n \end{pmatrix}$$
for some matrix $\L \in M_{n-r, n} (\CA )$ where $\der_i = \frac{\der}{\der x_i}$. By evaluating this equality of derivations at the elements $x_{j_{r+1}}, \ldots , x_{j_n}$ we obtain the equality of $(n-r)\times (n-r)$ matrices with coefficients in the algebra $\CA$, $$\D E_{n-r}=\L \CC_{\CA , \bj }$$ where $E_{n-r}$ in the $(n-r)\times (n-r)$ identity matrix and $\bj =(x_{j_{r+1}}, \ldots , x_{j_n})$. Let $e_1=(1,0, \ldots, 0), \ldots , 
e_{n-r}=(0, \ldots, 0, 1)$ be the standard basis for the free $\CA$-module $\CA^{n-r}=\bigoplus_{i=1}^{n-r}\CA e_i$. Let $R_1, \ldots , R_{n-r}\in \CA^{n-r}$ be the rows of the matrix $\D E_{n-r}=\L \CC_{\CA , \bj }$. In the $\CA$-module $\wedge^{n-r} \CA^{n-r}=\CA e_1\wedge\cdots \wedge e_{n-r}$ (the wedge product),  we have the inclusion  $$ \D^{n-r} e_1\wedge\cdots \wedge e_{n-r}=  R_1\wedge\cdots \wedge R_{n-r}\subseteq  \gm_\j   e_1\wedge\cdots \wedge e_{n-r}, $$
and so $\D^{n-r}\in \gm_\j$. 

$(3\Rightarrow 4)$  The implication follows from the inclusions $\gm_\j\subseteq \gc_{\CA , n-r}$ for all $\j \in \mJ_r$.

$(4\Rightarrow 2)$ By statement 4, the ideal $\gc_{\CA , n-r}$ contains all the elements $\D (\i , \j )^k$ where $\i \in \mI_r$,  $\j \in \mJ_r$ and  $k =k(\i , \j )$. The elements $\D(\i , \j )$ are generators of the Jacobian ideal $\ga_r$ of the algebra $\CA$. So, $$\ga_r^s\subseteq \gc_{\CA , n-r}$$ for some $s\geq 1$. The algebra $\CA$ is a regular algebra, hence $\ga_r = \CA$.  Now, $\gc_{\CA , n-r} =\CA$, and so $d\geq n-r$. By Proposition \ref{A23Jul19}.(2), $d\leq n-r$, and so $d=n-r$.

$(2\Rightarrow 1)$ Let $\CM_{n-r}$ be the set of nonzero $(n-r)\times (n-r)$ minors of the matrix $\CC_\CA$. Then $\CM_{n-r}\neq \emptyset$ since $d=n-r$ and  $\gc_{\CA , d}=\CA$. 

(i) {\em The algebraic extension $\CA \ra \prod_{\mu \in \CM_{n-r}}\CA_\mu$ is faithfully flat}: Since $\CA = \gc_{\CA , n-r}= (\mu)_{\mu \in \CM_{n-r}}$, the statement (i) follows.

(ii) {\em For all $\mu \in \CM_{n-r}$, $\Der_K(\CA_\mu ) = \CA_\mu \CH_{\CA_\mu }$ (where $\CA_\mu$ is the localization of the algebra $\CA$ at the powers of the element $\mu$)}: Fix a minor $\mu \in \CM_{n-r}$, it is deterimed by the rows with indices $i_{r+1}, \ldots , i_n$ and columns with indices $j_{r+1}, \ldots , j_n$. Let $\j =(j_1, \ldots , j_r)$ where $\{ j_1< \cdots < j_r\}=\{ 1, \ldots , n\} \backslash \{ j_{r+1}, \ldots , j_n\}$. 
Then 
$$\begin{pmatrix}
\pad_{x_{i_{r+1}}}\\  \vdots \\
\pad_{x_{i_n}} \end{pmatrix}= \L 
\begin{pmatrix}
\der_1\\  \vdots \\
\der_n \end{pmatrix}
$$
for some matrix $\L \in M_{n-r, n}(\CA )$. Let us consider the $(n-r)\times (n-r)$ matrix $M = (\{ x_{i_s}, x_{i_t}\})$ where $s , t = r+1, \ldots , n$. In particular, $\mu := \det (M)$ and $M^{-1} = \mu^{-1} \widetilde{M}$ where $\widetilde{M}$ is the adjoint matrix of the matrix $M$. Then 
$$\begin{pmatrix}
\g_{r+1}\\  \vdots \\
\g_n \end{pmatrix}:=
M^{-1}\begin{pmatrix}
\pad_{x_{i_{r+1}}}\\  \vdots \\
\pad_{x_{i_n}} \end{pmatrix}=\mu^{-1}\widetilde{M}
\begin{pmatrix}
\pad_{x_{i_{r+1}}}\\  \vdots \\
\pad_{x_{i_n}} \end{pmatrix}
\in \Bigg(\Der_K(\CA_\mu ) \cap \CA_\mu \CH_{\CA_\mu}\bigg)^{n-r}$$
and $\g_s(x_{j_t})=\d_{st}$ (the Kronecker delta) where $s , t = r+1, \ldots ,  n$. Therefore, $$\Der_K(\CA)_\mu = \bigoplus_{s=r+1}^n\CA_\mu \g_s=\CA_\mu \CH_{\CA_\mu}$$ (since $\dim_Q(Q\Der_k(\CA )) = n-r=\trdeg_K (Q)$ and the restriction of the derivations $\g_{r+1}, \ldots , \g_n$ to the subfield 
$K(x_{j_{r+1}}, \ldots , x_{j_n})$ of rational functions in the variables $x_{j_{r+1}}, \ldots , x_{j_n}$ are equal to the partial derivatives $\frac{\der}{\der x_{j_{r+1}}}, \ldots , \frac{\der}{\der x_{j_n}}$, respectively).

(iii) $\Der_K(\CA ) = \CA\CH_{\CA}$:  Let us consider the left $\CA$-module $V=\Der_K(\CA ) / \CA\CH_{\CA}$. By the statement (ii), $V_\mu=0$ for all elements $\mu \in \CM_{n-r}$. Hence $V=0$, by the statement (i). 
$\Box $

\begin{corollary}\label{p22Jul19}
Let $\CP = P_n$ be a polynomial  Poisson algebra. Then the  following statements are equivalent:
\begin{enumerate}
\item $\Der_K(P_n) = P_n\CH_{P_n}$. 
\item $\CC_\CP \in \GL_n(P_n)$.
\item $\det (\CC_\CP ) \in K^\times$. 
\end{enumerate}
\end{corollary}

{\it Proof}. $(1\Leftrightarrow 3)$  Notice that $\Der_K(P_n)=\bigoplus_{i=1}^n P_n \der_i$ where $\der_i= \frac{\der}{\der x_i}$. By Theorem \ref{22Jul19}.(1,2),  $\Der_K(P_n) = P_n\CH_{P_n}$ iff $d=n$ and $\gc_{\CP , n} = P_n$ iff $\det (\CC_\CP ) \in K^\times$ since $\gc_{\CP , n} =(\det (\CC_\CP ) )$. 

$(2\Leftrightarrow 3)$ The equivalence is obvious.
 $\Box $\\

{\it Example.} Let $\CP = P_2=K[x_1, x_2]$ be a Poisson algebra and $ a= \{ x_1, x_2\}$, an arbitrary element of $P_2$. Then $\Der_K(P_2) = P_2\CH_{P_2}$ iff $a\in K^\times$, by Corollary \ref{p22Jul19}.\\

{\bf The  generalized Weyl Poisson algebra $D[X,Y; a, \der \}$.} We apply the above results for certain generalized Weyl Poisson algebras. In particular,  Lemma \ref{a24Jul19} and  Corollary \ref{b24Jul19} are regularity and symplecticity criteria for certain GWPAs.

 {\it Definition, \cite{Pois-GWA}.} Let $D$ be a Poisson algebra, $\der = (\der_1, \ldots , \der_n)\in \PDer_K(D)^n$ be an $n$-tuple of commuting derivations of the Poisson algebra $D$, $a=(a_1, \ldots , a_n)\in \PZ (D)^n$ 
  be such that $\der_i (a_j)=0$ for all $i\neq j$. The {\em generalized Weyl algebra} $$A= D[ X, Y; (\id_D, \ldots , \id_D), a]=
 D[X_1, \ldots , X_n , Y_1, \ldots , Y_n]/(X_1Y_1-a_1, \ldots , X_nY_n-a_n)$$ admits a Poisson structure which is an extension of the Poisson structure on $D$ and is given by the rule: For all $i,j=1, \ldots , n$ and $d\in D$,

\begin{equation}\label{PGWAR1}
\{ Y_i, d\}=\der_i(d)Y_i, \;\; \{ X_i, d\}=-\der_i(d)X_i \;\; {\rm and}\;\; \{ Y_i, X_i\} = \der_i (a_i),
\end{equation}
\begin{equation}\label{PGWAR2}
\{ X_i, X_j\}=\{ X_i, Y_j\}=\{ Y_i, Y_j\} =0\;\; {\rm for\; all}\;\; i\neq j.
\end{equation}
The Poisson algebra is denoted by $A =D[ X, Y; a, \der \}$ and is called the {\bf generalized Weyl Poisson algebra} of rank $n$ (or GWPA, for short)  where $X=(X_1, \ldots , X_n)$ and $Y= (Y_1, \ldots , Y_n)$. Lemma \ref{a24Jul19} and  Corollary \ref{b24Jul19} are regularity and symplecticity criteria for certain GWPAs.

\begin{lemma}\label{a24Jul19}
Let $A= K[H][X,Y; a, \der=b\frac{d}{dH}\}$ be a GWPA of rank 1 where $a, b\in K[H]$ and $a\neq 0$. Then 
\begin{enumerate}
\item The domain $A= K[H,X,Y]/(XY-a)$ is regular iff  $(a, a')=K[H]$ where 
  $a'=\frac{da}{dH}$, i.e. the polynomials $a$ and $a'$ are co-prime. 
\item Suppose that the algebra $A$ is a regular algebra, i.e. $(a,a')=K[H]$. Then 
\begin{enumerate}
\item $\Der_K( A)=A\CH_A$ iff $b\in K^\times$.
\item $\CC_A=\begin{pmatrix}
\{X , X\}  &\{X , Y\} & \{X , H\}\\ \{ Y, X\} &\{ Y, Y\} &\{Y ,H \}\\
\{H , X\} &\{ H, Y\} & \{H , H\} \end{pmatrix} = \begin{pmatrix}
0 &-ba' & -bX\\ ba' &0 & bY\\
bX &-bY & 0 \end{pmatrix}$ and $r(\CC_A)=\begin{cases}
2& \text{if }b\neq 0,\\
0& \text{if }b=0.
\end{cases}$
\end{enumerate}
\end{enumerate}
\end{lemma}

{\it Proof}. 1. The Jacobian matrix of $A$ is equal to $(Y,X, -a')$. Then the Jacobian ideal  $\ga_1$  of $A$ is equal to  $(Y,X, -a')$, the ideal of $A$ which is generated by the elements $Y$, $X$ and  $-a'$. Now, $\ga_1 = A$ iff 
$K[H] = K[H]\cap (Y,X, -a')= (a,a')$ since $XY=a$ (use the $\Z$-grading of the GWA $A$). 

2. (b) The statement (b) is obvious.

(a) By Theorem \ref{22Jul19}.(1,2), $\Der_K(A) = A\CH_A$ iff $r(\CC_\CP )=3-1=2$ (i.e. $b\neq 0$, by the statement (b)) and $A= \gc_{A, 2}= b^2(a'^2, a' X, a'Y, XY,  X^2, Y^2)$ iff $b\in K^\times $ and $$A= (a'^2, a' X, a'Y, a, X^2, Y^2)=(a'^2, a, X, Y)$$ since $(a,a')=1$  iff $b\in K^\times $ and $K[H]=K[H]\cap (a'^2, a, X, Y)=(a'^2, a)$ iff  $b\in K^\times $ since $K[H]=(a,a')$ implies $K[H]=(a,a'^2)$. $\Box $

\begin{corollary}\label{b24Jul19}
Let $A= K[H_1, \ldots , H_n][X,Y; a=(a_1, \ldots , a_n), \der=(b_1\frac{d}{dH_1}, \ldots ,b_n\frac{d}{dH_n} )\}$ be a GWPA of rank n where $a_i, b_i\in K[H_i]$ and $a_i\neq 0$. Then 
\begin{enumerate}
\item The domain $A= K[H_1, \ldots , H_n, X_1, \ldots , X_n, Y_1, \ldots , Y_n]/(X_1Y_1-a_1, \ldots , X_nY_n-a_n)$ is regular iff  $(a_i, a_i')=K[H_i]$ for $i=1, \ldots , n$  where 
  $a_i'=\frac{da_i}{dH_i}$, i.e. the polynomials $a_i$ and $a_i'$ are co-prime in $K[H_i]$. 
\item Suppose that the algebra $A$ is a regular algebra. Then 
\begin{enumerate}
\item $\Der_K( A)=A\CH_A$ iff $b_i\in K^\times$ for $i=1, \ldots, n$.
\item  $r(\CC_A)=2n$ iff  $b_i\neq 0$ for $i=1, \ldots, n$.
\end{enumerate}
\end{enumerate}
\end{corollary}

{\it Proof}.  Notice that $A=\bigotimes_{i=1}^n A_i$, a tensor product of Poisson algebras where $A_i= K[H_i][X_i,Y_i; a_i, \der_i=b_i\frac{d}{dH_i}\}$ is a GWPA as in Lemma \ref{a24Jul19} for $i=1, \ldots, n$. Now, the corollary follows from  Lemma \ref{a24Jul19}. 
$\Box $


\section{The kernel of the epimorphism $\CU (\CP ) \ra P\CD (\CP )$ and the defining relations of the algebra $P\CD (\CP )$}\label{KERCUCP}

In this section we keep the assumptions and the notations of Section \ref{CUADOM}. In most cases the field $K$ has characteristic zero.

Let $\CP$ be a Poisson algebra. The Poisson algebra $\CP$ is a Poisson $\CP$-module, i.e. a left $\CU (\CP )$-module,  and the image of the corresponding algebra homomorphism $\CU (\CP ) \ra \End_K(\CP ) $ is the algebra of Poisson differential operators $P\CD (\CP )$ on $\CP$. So, the kernel of the algebra epimorphism
\begin{equation}\label{UPDP}
\pi_\CP : \CU (\CP )\ra P\CD (\CP ) , \;\; p\mapsto p, \;\; \d_q\mapsto \pad_q:=\{ q, \cdot \} \;\; (p,q\in \CP )
\end{equation}
is the annihilator of the left $\CU ( \CP )$-module $\CP $, 
\begin{equation}\label{ker=anP}
\ker (\pi_\CP ) = \ann_{\CU (\CP )}(\CP ). 
\end{equation}
For a  domain of essentially  finite type $\CA$ over a field of characteristic zero, Proposition \ref{B26Jul19} determines the exact value for the Gelfand-Kirillov dimension of the algebra $P\CD (\CA )$. An explicit ideal $\kappa_\CA$ of the algebra $\CU (\CA )$ is introduced such that $\kappa_\CA \subseteq \ker (\pi_\CA )$, see (\ref{dinnj}) and 
 (\ref{dinnj1}) for the definition of its explicit  generators. Corollary \ref{a2Aug19} is a criterion  for 
 $\kappa_\CA =\ker (\pi_\CA )$. Proposition \ref{26Jul19} is about properties of the ideals $\kappa_\CA $ and $\ker (\pi_\CA )$.
 Proposition \ref{26Jul19}.(3) is a sufficient condition for 
 $\kappa_\CA =\ker (\pi_\CA )$. Proposition \ref{A2Aug19} describes the torsion submodule of the module $\O_\CA$ of K\"{a}hler differential of the algebra $\CA$. 
 Theorem \ref{2Aug19} gives explicit descriptions of the kernels $\ker (\pi_\CA )$ and $\ker (\overline{\pi}_\CA )$.  Theorem \ref{3Aug19} and Theorem \ref{A5Aug19} are criteria for 
 $\ker (\pi_\CA )=0$ and their proofs are given in this section.  Theorem \ref{5Aug19} is a criterion for the  homomorphism $\pi_\CA : \CU (\CA ) \ra \CD (\CA )$ to be an isomorphism, its proof is also presented in this section. 

In this section,  the algebra $\CP = \CA = S^{-1}(P_n/I)$ is a domain of essentially finite type and we keep the assumptions and the notations of Section \ref{CUADOM}. Let  $Q=Q(\CA )$ be its field of fractions, $d=d_\CA = r(\CC_\CA )$ be the rank of the $n\times n$ matrix $\CC_\CA =(\{ x_i , x_j\} )\in M_n(\CA )$ over
 the field $Q$. For each $l=1, \ldots , n$, let $$\ind_n(l) =\{ \i =(i_1, \ldots , i_l)\, | \, 1\leq i_1 <\cdots <i_l\leq n\}.$$ For  elements $\i = (i_1, \ldots , i_l)$ and $\j = (j_1, \ldots , j_l)$ of $\ind_n(l)$, let $\CC_\CA (\i , \j )=(\{ x_{i_\nu} , x_{j_\mu } \} )$ be the $l\times l$ submatrix of the matrix $\CC_\CA$. So, the rows (resp., the columns) of the matrix  $\CC_\CA (\i , \j )$ are indexed by $i_1, \ldots , i_l$ (resp., $j_1, \ldots , j_l$). The $(i_\nu , j_\mu )^{th}$ element of the matrix  $\CC_\CA (\i , \j )$ is $\{ x_{i_\nu} , x_{j_\mu }\}$. Let $M_{\CA , l} =\{ \CC_\CA (\i , \j ) \, | \, \i, \j \in \ind_n(l)\}$ be the set of all $l\times l$ submatrices of $\CC_\CA$ and 
 $$\CC_{\CA , l} =\{ \mu (\i , \j ) :=\det (\CC_\CA (\i , \j )) \, | \, \i, \j \in \ind_n(l)\}$$ be the set of all $l\times l$ minors of $\CC_\CA$. Let 
 \begin{eqnarray*}
 \mI (l)&=& \mI_\CA (l) = \{ \i \in \ind_n(l)\, | \, \mu (\i , \j ) \neq 0\;\; {\rm for \; some}\;\; \j \in \ind_n(l)\} \\
 \mJ (l)& = &\mJ_\CA (l) = \{ \j \in \ind_n(l)\, | \, \mu (\i , \j ) \neq 0\;\;  {\rm for \; some}\;\;\i \in \ind_n(l)\}
\end{eqnarray*}
The symmetric group $S_n$ acts on the set of indices $\{ 1, \ldots , n\}$ by permutation. In particular, the symmetric group $S_l$ acts on the set $\{ i_1, \ldots , i_l\}$ by permutation. For each $\s \in S_l$, $\sign (\s )$ is the {\em sign} of $\s$. 

\begin{lemma}\label{a26Jul19}

\begin{enumerate}
\item For all $\i , \j \in \ind_n(l)$, $\CC_\CA(\i , \j )^t=-\CC_\CA (\j , \i )$ where `t' stands for the trasposition of a matrix.
\item For all $\i , \j \in \ind_n(l)$, $\mu (\i , \j ) = (-1)^l \mu (\j , \i )$.
\item $\mI_\CA (l) = \mJ_\CA (l)$ for all $l=1, \ldots, n$.
\item For all permutations $\s , \tau \in S_l$, $\mu (\s (\i ), \tau (\j ))= \sign (\s \tau ) \mu (\i , \j )$. 
\end{enumerate}
\end{lemma}

{\it Proof}. 1. The $(j_t, i_s)^{th}$ element of the matrix $\CC_\CA (\i , \j )^t$ is $\{ x_{i_s} , x_{j_t}\}= -\{ x_{j_t}, x_{i_s}\}$, and statement 1 follows.

2. $\mu (\i , \j ) = \det (\CC_\CA (\i , \j ))= \det (\CC_\CA (\i , \j )^t)\stackrel{\rm st.1}{=} 
\det (-\CC_\CA (\j , \i )) = (-1)^l \det (\CC_\CA (\j , \i ))=(-1)^l \mu (\j , \i )$. 

3. Statement 3 follows from statement 2.

4. Statement 4 is obvious.  $\Box $\\

In the case $l=d$, we can say more about the elements $\mu (\i , \j )$.

\begin{proposition}\label{A26Jul19}
Let $K$ be a field of characteristic zero. 
\begin{enumerate}
\item Given $\i = (i_1, \ldots , i_d) \in \ind_n(d)$ where $d=r(\CC_\CA )$. Then 
$\i\in \mI (d)$ iff $Q\CH _\CA  =\bigoplus_{s=1}^d Q\d_{x_{i_s}}'$ where $\d_{x_{i_s}}'=\d_{i_s}'=\sum_{j=1}^n\{ x_{i_s}, x_j\} \der_j$. 
\item For all $\i \in \mI_\CA (d)$ and $\j \in \mJ_\CA (d)$, $\mu (\i , \j ) \neq 0$.
\item Given $\i , \j  \in \ind_n(d)$. Then $\i , \j \in \mI_\CA (d)=\mJ_\CA (d)$ iff $\mu (\i , \j ) \neq 0$. 
\end{enumerate}
\end{proposition}

{\it Proof}. 1. Statement 1 follows from the fact that $\dim_Q(Q\CH_\CA ) = d$ and the elements $\d_{x_{i_s}}'=\sum_{j=1}^n\{ x_{i_s} , x_j\} \der_j$ $(s=1, \ldots , d)$ of the vector space $Q\CH_\CA$ over the field $Q$ are $Q$-linearly independent iff $\i\in \mI (d)$. 

2. Since $\i \in \mI_\CA (d)$ and $\j \in \mJ_\CA (d)$, $\mu (\i , \j')\neq 0$ and $\mu (\i', \j )\neq 0$ for some $\i', \j'\in \mI_\CA (d) = \mJ_\CA (d)$, Lemma \ref{a26Jul19}.(3). By statement 1, $$Q\CH_\CA = \bigoplus_{s=1}^dQ\d_{i_s}'= \bigoplus_{s=1}^d Q\d_{i_s'}'.$$ Let $\d_\i'=(\d_{i_1}', \ldots , \d_{i_d}')^t$ and $\d_{\i'}'=(\d_{i_1'}', \ldots , \d_{i_d'}')^t$. Then $\d_{\i'}'=\CC_\CA (\i' , \j') \CC_\CA (\i , \j' )^{-1}\d_{\i}'$ and 
 $\d_{\i}'=\CC_\CA (\i , \j) \CC_\CA (\i' , \j )^{-1}\d_{\i'}'$. Therefore, 
 $$ \d_{\i}'= \CC_\CA (\i , \j ) \CC_\CA (\i' , \j )^{-1}\CC_\CA (\i' , \j') \CC_\CA (\i , \j' )^{-1} \d_{\i}'$$
and so the product of four matrices is equal to the identity matrix. Hence, their determinants are nonzero. In particular, $\mu ( \i , \j ) =\det(\CC_\CA (\i , \j )) \neq 0$, as required. 

3. Statement 3 follows at once from statement 2. $\Box $\\


{\bf The ideal $\kappa_\CA$ of the algebra $\CU (\CA )$ and its generators $\d_{\i , i_\nu ; \j }$.} For each pair of elements $\i = (i_1, \ldots , i_d)$ and $\j = (j_1, \ldots , j_d)$ of $\mI_\CA (d)=\mJ_\CA (d)$ and each element $i_\nu\in \{ i_{d+1}, \ldots , i_n\}=\{ 1, \ldots , n\} \backslash \{ i_1, \ldots , i_d\}$, let us  consider the following elements of the algebra $\CU (\CA )$ (see Theorem \ref{23Jun19}.(2) for the definition of the elements $\d_i$), 
\begin{eqnarray}\label{dinnj}
 \d_{\i , i_\nu; \j }:= \det
 \begin{pmatrix}
 \{ x_{i_1}, x_{j_1} \} & \ldots &   \{ x_{i_1}, x_{j_d}\} & \d_{i_1}\\
  \vdots & \vdots & \vdots  & \vdots\\
   \{ x_{i_d}, x_{j_1} \} & \ldots &   \{ x_{i_d}, x_{j_d} \}& \d_{i_d}\\
  \{ x_{i_\nu}, x_{j_1} \} & \ldots &   \{ x_{i_\nu}, x_{j_d} \}& \d_{i_\nu}\\
\end{pmatrix}= \mu (\i , \j ) \d_{i_\nu} + \sum_{s=1}^d(-1)^{s+d+1}\mu (i_1, \ldots ,\widehat{i_s}, \ldots , i_d, i_\nu;\j )\d_{i_s}. 
\end{eqnarray}
 
Let $\{ j_{d+1}, \ldots , j_n\} := \{ 1, \ldots , n\} \backslash \{ j_1, \ldots , j_d\}$ and for each $\nu \in \{ d+1, \ldots , n\}$, let us  consider the following elements of $\CU (\CA )$, 

\begin{eqnarray}\label{dinnj1}
 \d_{\i ;  \j ,  j_\nu }':= \det
 \begin{pmatrix}
 \{ x_{i_1}, x_{j_1} \} & \ldots &   \{ x_{i_1}, x_{j_d}\} & \{ x_{i_1}, x_{j_\nu}\} \\
  \vdots & \vdots & \vdots  & \vdots\\
   \{ x_{i_d}, x_{j_1} \} & \ldots &   \{ x_{i_d}, x_{j_d} \}& \{ x_{i_d}, x_{j_\nu} \} \\
  -\d_{j_1} & \ldots &  -\d_{j_d}& -\d_{j_\nu}\\
\end{pmatrix}= -\mu (\i , \j ) \d_{j_\nu} - \sum_{s=1}^d(-1)^{s+d+1}\mu (\i ; j_1,\ldots , \widehat{j_s}, \ldots , j_d, j_\nu)\d_{j_s}. 
\end{eqnarray}
The elements $\{ \d_{\i ; \j , j_\nu}'\}$ are `dual' of the elements 
 $\{ \d_{\i , i_\nu; \j }\}$. In fact, they are the same up to sign, 
\begin{equation}\label{dinnj2}
\d_{\i ;  \j ,  j_\nu }'= (-1)^{d+1} \d_{\j ,  j_\nu ;\i  }.
\end{equation}
{\it Proof.} Let $M$ and $M'$ be the matrices in the definitions of the elements $\d_{\j ,  j_\nu ;\i  }$  and $\d_{\i ;  \j ,  j_\nu }'$, respectively, (see (\ref{dinnj}) and (\ref{dinnj1})). Then 
$$ \d_{\i ;  \j ,  j_\nu }'=\det (M') = \det ((M')^t) = \det (-M)=(-1)^{d+1}\det (M) = (-1)^{d+1} \d_{\j ,  j_\nu ;\i }.\;\; \Box$$

{\it Definition.} Let $\kappa_\CA$ be the ideal of the algebra $\CU (\CA )$ which is generated by all the elements $\d_{\i , i_\nu; \j }$ above where $\i , \j \in \mI_\CA (d)$, see (\ref{dinnj}).\\

 By (\ref{dinnj2}), the ideal $\kappa_\CA$ is also generated by all the elements $\d_{\i ;  \j ,  j_\nu }'$, see (\ref{dinnj1}). Since $d$ is the rank of the matrix $\CC_\CA$ and $\pi_\CA (\d_i) =\{ x_i, \cdot \}$ for all $i=1, \ldots , n$, we have that  $$\pi_\CA (\d_{\i, i_\nu ; \j} )(x_{j_\mu })\stackrel{(\ref{dinnj})}{=}\mu (\i , i_\nu; \j , j_\mu )=0\;\; {\rm  for \; all}\;\; i_\nu, j_\mu = 1, \ldots , n$$ since all the  $(d+1)\times (d+1)$ minors of the matrix $\CC_\CA$ of rank $d$ are equal to zero. Therefore, 
\begin{equation}\label{kAkpA}
\kappa_\CA \subseteq \ker (\pi_\CA ).
\end{equation}

{\bf The Gelfand-Kirillov dimension  of 
the algebra $P\CD (\CA )$.} Proposition \ref{B26Jul19} gives the exact figure for the Gelfand-Kirillov dimension  of 
the algebra $P\CD (\CA )$.
\begin{proposition}\label{B26Jul19}
Let a Poisson algebra $\CA$ be a domain of essentially finite type over the field $K$ of characteristic zero, $r$ is the rank of Jacobian matrix of $\CA$  and $d=r(\CC_\CA )$. Then 
\begin{enumerate}
\item $\GK (P\CD (\CA )) = \GK (\CA ) +d= n-r+d $.
\item Let $\mu = \mu (i_1, \ldots , i_d; j_1, \ldots , j_d)= \det(\{ x_{i_s}, x_{j_t}\} )$ be a nonzero  minor of the matrix $\CC_\CA$ where $s,t=1, \ldots , d$. Then $P\CD (\CA )_\mu \simeq P\CD (\CA_\mu ) = \bigoplus_{\alpha \in \N^d} \CA_\mu \pad_{x_{i_1}}^{\alpha_1}\cdots \pad_{x_{i_d}}^{\alpha_d}= \bigoplus_{\alpha \in \N^d} \CA_\mu \pad_{x_{j_1}}^{\alpha_1}$ $\cdots \pad_{x_{j_d}}^{\alpha_d}$, and $\GK (P\CD (\CA)_\mu ) = \GK (\CA ) +d$. 
\end{enumerate}
\end{proposition}

{\it Proof}. 2. By (\ref{dinnj}) and (\ref{kAkpA}), the localization of the algebra  $P\CD (\CA)$ at the powers of the element $\mu$ is equal to $P\CD (\CA)_\mu = \sum_{\alpha \in \N^d} \CA_\mu \pad_{x_{i_1}}^{\alpha_1}\cdots \pad_{x_{i_d}}^{\alpha_d}$. Let 
$$\begin{pmatrix}
\g_{i_1}\\  \vdots \\
\g_{i_d} \end{pmatrix}:=\CC_\CA (\i , \j )^{-1}\begin{pmatrix}
\pad_{x_{i_1}}\\  \vdots \\
\pad_{x_{i_d}} \end{pmatrix}.$$
Then $P\CD (\CA )_\mu =\sum_{\alpha\in \N^d}\CA_\mu \g_{i_1}^{\alpha_1}\cdots \g_{i_d}^{\alpha_d}$. The sum is a direct sum since
$$ \g_{i_\nu} (x_{j_\mu}) = \d_{\nu \mu}\; \; {\rm for\; all}\;\; \nu , \mu =1, \ldots , d$$
 and the field $K$ has characteristic zero where $\d_{\nu \mu}$ is the Kronecker delta. Therefore $P\CD (\CA)_\mu = \bigoplus_{\alpha \in \N^d} \CA_\mu \pad_{x_{i_1}}^{\alpha_1}\cdots \pad_{x_{i_d}}^{\alpha_d}$,  and so $\GK (P\CD (\CA )_\mu )=\GK (\CA ) +d=n$.

1. By statement 2, $\GK (P\CD (\CA ))\leq \GK (P\CD (\CA)_\mu ) = \GK (\CA ) +d$. By statement 2, $$\bigoplus_{\alpha \in \N^d} \CA \pad_{x_{i_1}}^{\alpha_1}\cdots \pad_{x_{i_d}}^{\alpha_d}\subseteq P\CD (\CA ).$$ Hence, 
$\GK (P\CD (\CA ))\geq \GK (\CA )+d$, and statement 1 follows since $\GK (\CA ) = n-r$.  $\Box $ \\

Recall that {\em every left (resp., right) Ore set of a left (resp., right) Noetherian ring is a left (resp., right) denominator set}, see \cite[Corollary 4.24]{Good-WarIntrNCNR}. 

The ring $\CU (\CA )$ is a Noetherian ring and for each nonzero element $s\in \CA$ the set $S_s=\{ s^i \, | \, i\geq 0\}$ is an Ore set of $\CU (\CA )$. Hence, $S_s$ is a denominator set of $\CU (\CA )$, and so $$\CU (\CA)_s:=S_s^{-1}\CU (\CA ) \simeq \CU (\CA ) S_s^{-1}$$ is the localization of the ring $\CU (\CA )$ at the powers of the element $s$. By (\ref{kAkpA}) and the inclusion $\CA \subseteq P\CD (\CA )$, $$\CA \cap \kappa_\CA =0.$$ The ring $\CU (\CA ) / \kappa_\CA $ is a Noetherian ring. Hence, the set $S_s$ is a denominator set of the ring $\CU (\CA ) / \kappa_\CA $. 

Recall that the algebra $P\CD (\CA )$ is equipped with two filtrations: the order filtration $\{ P\CD (\CA )_i= P\CD (\CA ) \cap \CD (\CA )_i\}_{i\geq 0} $ and the filtration $\{ P\CD (\CA )_{\leq i}\}_{i\geq 0}$ that is determined by the total degree of the elements $\d_{x_1}, \ldots , \d_{x_n}$. Clearly,
$$P\CD (\CA )_{\leq i}\subseteq P\CD (\CA)_i \;\; {\rm  for \; all}\;\;i\geq 0.$$
Theorem \ref{26Jul19} shows that if $\gc_{\CA , d} =\CA$ then $\ker (\pi_\CA ) =\kappa_\CA$ and $P\CD (\CA )_{\leq i}= P\CD (\CA)_i$ for   all $i\geq 0$.

\begin{theorem}\label{26Jul19}
Let  a Poisson algebra $\CA$ be a domain of essentially finite type over the field $K$ of characteristic zero and $d=r(\CC_\CA )$. Then 
\begin{enumerate}
\item $\gc_{\CA , d}^s\ker (\pi_\CA )\subseteq \kappa_\CA$ and $\ker (\pi_\CA )\gc_{\CA , d}^t\subseteq \kappa_\CA$ for some natural numbers $s$ and $t$.
\item For each $i\geq 0$, the left $\CA$-module $P\CD (\CA )_i$ is finitely generated and there are natural numbers $s_i$ and $t_i$ such that $\gc_{\CA , d}^{s_i}P\CD (\CA )_i\subseteq P\CD (\CA )_{\leq i}$ and $P\CD (\CA )_i\subseteq P\CD (\CA )_{\leq i}\gc_{\CA , d}^{t_i}$.
\item Suppose that $\gc_{\CA , d}=\CA$. Then 
\begin{enumerate}
\item $\ker (\pi_\CA ) = \kappa_\CA $. 
\item For all $i\geq 0$, $P\CD (\CA )_i=P\CD (\CA )_{\leq i}$. 
\end{enumerate}
\end{enumerate}
\end{theorem}

{\it Proof}. 1. (i) {\em For all} $\i , \j \in \mI_\CA (d)$, $\Big (\CU (\CA ) / \kappa_\CA \Big)_{\mu (\i , \j )} = \bigoplus_{\alpha \in \N^d} \CA_\mu 
\d_{x_{i_1}}^{\alpha_1}\cdots \d_{x_{i_d}}^{\alpha_d}$. 
 The statement (i) follows at once from (\ref{dinnj}) and Proposition \ref{B26Jul19}.(2). 
 
 (ii) {\em For all } $\i , \j \in \mI_\CA (d)$, $\Big( \ker (\pi_\CA )/\kappa_\CA \Big)_{\mu (\i , \j )}=0$: When we localize the short exact sequence of left $\CU (\CP )$-modules,
 $$ 0\ra \ker (\pi_\CA ) / \kappa_\CA \ra \CU (\CA ) / \kappa_\CA \ra P\CD (\CA )\ra 0, $$
 at the powers at the element $\mu = \mu (\i , \j )$ 
we obtain the short exact sequence of left $\CU (\CP )_{\mu}$-modules,
 $$ 0\ra \ker (\pi_\CA )_{\mu} / (\kappa_\CA )_{\mu}\ra \CU (\CA )_{\mu} / (\kappa_\CA)_{\mu} \ra P\CD (\CA )_{\mu}\ra 0. $$
By the statement (i) and Proposition \ref{B26Jul19}.(2), the second map of the short exact sequence above is an isomorphism, and the statement (ii) follows. 

The algebra $\CU (\CA )$ is a Noetherian algebra. Therefore, the ideal $\kappa_\CA$ is a finitely generated left and right ideal. The localization at the powers of the element $\mu$ is a left and right localization. Hence, $$\mu^{s(\i , \j )}\ker (\pi_\CA ) \subseteq \kappa_\CA\;\; {\rm  and}\;\;\ker (\pi_\CA ) \mu^{t(\i , \j )}\subseteq \kappa_\CA$$ for some natural numbers $s(\i , \j )$ and $t( \i , \j )$. The ideal $\gc_{\CA , d}$ is generated by the finitely many elements $\{ \mu (\i , \j ) \, | \, \}_{\i , \j \in \mI_\CA (d)}$. Now, statement 1 is obvious. 

2. By \cite[Proposition 5.3.(2)]{gendifreg}, for each $i\geq 0$, the left $\CA$-module $\CD(\CA )_i$ is finitely generated, hence Noetherian. The proof of \cite[Proposition 5.3.(2)]{gendifreg} works for right $\CA$-modules, i.e. for each $i\geq 0$, the right $\CA$-modules $\CD (\CA )_i$ is a finitely generated, hence Noetherian. Therefore, the left and right $\CA$-submodule $P\CD (\CA )_i$ of $\CD (\CA )_i$ is finitely generated, and so is its quotient $$V(i) = P\CD (\CA )_i/ P\CD (\CA )_{\leq i}.$$ By Proposition \ref{B26Jul19}.(2), $V(i)_{\mu (\i , \j )}=0$ for all $\i , \j \in \mI_\CA (d)$. Therefore, $\gc_{\CA , d}^{s_i}V(i)=0$ and $V(i)\gc_{\CA , d}^{t_i}=0$
 for some natural numbers $s_i$ and $t_i$, and statement 2 follows. 
 
 3. Since $\gc_{\CA , d}=\CA $, the statements (a) and (b) follow from statements 1 and 2, respectively.  $\Box $\\
 
 {\bf The torsion $\CA$-submodule $T_\CA$ of $\O_\CA$.} Let a Poisson algebra $\CA = S^{-1}(P_n/I)$ be a domain of essentially finite type. Then 
 $$T=T_\CA = \{ \o \in \O_\CA \, | \, a\o =0\;\; {\rm  for \; some}\;\; a\in \CA \backslash \{ 0\} \}$$ 
 is called the {\em torsion} $\CA$-submodule of $\O_\CA$. The set $\CA \backslash \{ 0\}$ of regular elements of $\CA$ is a denominator set of $\CA$ and $T$ is the $\CA \backslash \{ 0\}$-torsion submodule of $\O_\CA$. For each $\i \in \mI_r$ and $\j \in \mJ_r$, $\D = \D (\i , \j )\neq 0$, and the set $S_\D =\{ \D^i\, | \, i\geq 0\}$ is a denominator set of the algebra $\CA$. Let $$\tor_{S_\D}(\O_\CA ) := \{ \o \in \O\, | \,\D^i\o =0\;\; {\rm  for \; some} \;\;i\geq 0\},$$ the $S_\D$-{\em torsion submodule} of the $\CA$-module $\O$.\\

Proposition \ref{A2Aug19} is an explicit description of the torsion submodule $T$ of $\O$. 

\begin{proposition}\label{A2Aug19}
Let a Poisson algebra $\CA = S^{-1}(P_n/I)$ be a domain of  essentially finite type.  Then 
\begin{enumerate}
\item $T_\O = \tor_{S_{\D (\i , \j )}}(\O_\CA )$ for all $\i  \in \mI_r$ and $\j \in \mJ_r$.
\item  $T_\O =\ann_{\O}(\ga_r^t)$  for all $t\gg 0$ (where $r$ is the rank of the Jacobian matrix $\Big(\frac{\der f_i}{\der x_j}\Big)$ and $\ga_r$ is the Jacobian ideal).
\item If, in addition, the algebra $\CA$ is a regular algebra then $T_\O =0$. 
\end{enumerate}
\end{proposition}

{\it Proof}. 1. Let $\D = \D (\i , \j )$ where $\i \in \mI_r$ and $\j \in \mJ_r$. Then $S_\D \subseteq \CA \backslash \{ 0\}$, and so $\tor_{S_\D } (\O ) \subseteq \tor_{\CA \backslash \{ 0\} }(\O) =T$. The opposite inclusion follows from the fact that the $\CA$-module $$S_\D^{-1} \O = \oplus_{\nu = r+1}^n \CA_\D dx_{j_\nu}$$ is torsion free where $\{ j_{r+1}, \ldots ,  j_n\}=\{ 1, \ldots , n\}\backslash \{ j_1, \ldots , j_r\} $. 

2. The $\CA$-module $\O$ is finitely generated, hence Noehterian. Then the $\CA$-submodule $T$ of $\O$ is finitely generated. By statement 1, $\D(\i , \j )^kT=0$ for all $\i \in \mI_r$, $\j \in \mJ_r$ and some natural number $k= k(\i , \j )\geq 0$. Therefore, $$\ga_r^sT=0\;\; {\rm  for \; some}\;\; s\geq 0$$ since the ideal $\ga_r$ is generated by the finite set $\{ \D (\i , \j ) \, | \, \i \in \mI_r, \j \in \mJ_r\}$. Therefore, $T\subseteq \ann_\O (\ga_r^s)$. Since $\D (\i , \j )^l\in \ga_r^s$ for some $l\geq 1$,  $$\ann_\O (\ga_r^s) \subseteq \ann_\O (\D (\i , \j )^l)\subseteq \tor_{S_{\D (\i , \j )}}(\O ) =T,$$ by statement 1, and statement 2 follows. 

3. The algebra $\CA$ is regular. So, $\ga_r=\CA$, and statement 3 follows from statement 2. $\Box $\\

{\bf The $\CA$-torsion submodule $T_{\bCU (\CA }$ of $\bCU (\CA )$.} Let us consider the factor algebra 
\begin{equation}\label{CAkA2}
\bCU := \bCU (\CA ) :=\CU (\CA ) / \kappa_\CA .
\end{equation}
 There are  obvious algebra epimorphisms (since $\kappa_\CA \subseteq \ker (\pi_\CA )$):
\begin{equation}\label{CAkA1}
\pi_\CA : \CU (\CA ) \stackrel{\phi}{\longrightarrow}\bCU (\CA ) \stackrel{\overline{\phi}_\CA}{\longrightarrow} P\CD (\CA)\;\; {\rm and }\;\; \pi_\CA = \overline{\phi}_\CA\phi.
\end{equation}
Since $\CA \subseteq P\CD (\CA)$, 
\begin{equation}\label{CAkA}
\CA \cap \kappa_\CA =0.
\end{equation}
The set $\CA \backslash \{ 0\}$ is an Ore set of the Noetherian algebra $\CU (\CA )$. Hence the set $\CA \backslash \{ 0\}$ is a denominator set of $\CU (\CA )$. By (\ref{CAkA}), the set   $\CA \backslash \{ 0\}$ is an Ore set of the Noetherian algebra $\bCU $, hence it is a denominator set of $\bCU $.  Let $$T_{\bCU } = T_{\bCU (\CA )}= \tor_{\CA \backslash \{ 0\} }(\bCU (\CA ))$$ be the $\CA \backslash \{ 0\}$-torsion submodule of the algebra $\bCU$. The set $T_{\bCU}$ is an ideal of the algebra $\bCU$. 
 Abusing the language we often say `$\CA$-torsion' meaning `$\CA\backslash \{ 0\}$-torsion'. 

\begin{proposition}\label{B2Aug19}
Let a Poisson algebra $\CA= S^{-1}(P_n/I)$ be a domain of essentially finite type over the field $K$  of characteristic zero and $d=r(\CC_\CA )$. Then 
\begin{enumerate}
\item $T_{\bCU (\CA )} = \tor_{S_{\mu (\i , \j )}}(\bCU (\CA ))$ for all $\i \in \mI_\CA (d)$ and $\j \in \mJ_\CA (d)$ where $S_{\mu (\i , \j )}=\{ \mu (\i , \j )^k\, | \, k\geq 0\}$ is a denominator set of the algebra $\bCU (\CA )$. 
\item $T_{\bCU (\CA )} = \rann_{\bCU (\CA )}(\gc_{\CA , d}^t)$ and $T_{\bCU (\CA )} = \lann_{\bCU (\CA )}(\gc_{\CA , d}^t)$ for all $t\gg 0$. 
\end{enumerate}
\end{proposition}

{\it Proof}. 1. Let $\mu = \mu (\i , \j )$ where $\i \in  \mI_\CA (d)$ and $\j \in \mJ_\CA (d)$. Then the set $S_\mu = \{ \mu^i\, | \, i\geq 0\}$ is a denominator set of the Noetherian algebra $\CU (\CA )$ such that $S_\mu \in \CA \backslash \{ 0\}$. By (\ref{CAkA}), $S_\mu  \cap \kappa_\CA =\emptyset $. Therefore, the set $S_\mu$ is a denominator set of the Noetherian algebra $\bCU = \bCU (\CA )$. Then the inclusion $S_\mu \subseteq \CA \backslash \{ 0\}$ implies the inclusion $$\tor_{S_\mu} (\bCU ) \subseteq T_{\bCU }.$$ The inverse inclusion follows from the fact that (see the statement (i) of the proof of Theorem \ref{26Jul19}.(1))
$$ S_\mu^{-1}\bCU = \bigoplus_{\alpha \in \N^d} \CA_\mu \d_{x_{i_1}}^{\alpha_1}\cdots \d_{x_{i_d}}^{\alpha_d}= \bigoplus_{\alpha \in \N^d}  \d_{x_{i_1}}^{\alpha_1}\cdots \d_{x_{i_d}}^{\alpha_d}\CA_\mu=\bCU S^{-1}_\mu.$$

2. The algebra $\CU = \CU (\CA )$ is a Noetherian algebra, hence so is its factor algebra $\bCU$. Then the ideal $T_{\bCU }$ of $\bCU$ is a finitely generated left and right $\bCU$-module. By statement 1, $\mu (\i , \j )^kT_{\bCU } = T_{\bCU }\mu (\i , \j )^k =0$ for all $\i \in \mI_\CA (d)$, $\j \in \mJ_\CA (d)$ and some natural number $k= k(\i , \j )\geq 0$. Therefore, $\gc_{\CA , d}^sT_{\bCU }=0$ and $T_{\bCU }\gc_{\CA , d}^s=0$ for some natural number $s\geq 0$ since the ideal $\gc_{\CA , d}$ is generated by the finite set $\{ \mu (\i  , \j ) \, | \, \i \in \mI_\CA (d) , \j \in \mJ_\CA (d)\}$.  Therefore, $$T_{\bCU}\subseteq \rann_{\bCU}(\gc_{\CA , d}^s)\;\; {\rm  and}\;\; T_{\bCU}\subseteq \lann_{\bCU}(\gc_{\CA , d}^s).$$ Since $\mu^s\in \gc_{\CA , d}^s$  (where $\mu = \mu (\i , \j )$), $$\rann_{\bCU}(\gc_{\CA , d}^s)\subseteq \rann_{\bCU}(\mu^s)\subseteq \tor_{S_\mu} (\bCU )= T_{\bCU},$$ by statement 1. Similarly, $\lann_{\bCU}(\gc_{\CA , d}^s)\subseteq T_{\bCU}$
 (by replacing `r' by `l' in the above chain of inclusions), and statement 2 follows.  $\Box $\\
 
 {\bf Explicit descriptions of $\ker (\pi_\CA )$ and $\ker (\overline{\pi}_\CA )$;  defining relations of the algebra $P\CD (\CA )$.} Theorem \ref{2Aug19} gives explicit descriptions of the kernels $\ker (\pi_\CA )$ and $\ker (\overline{\pi}_\CA )$. It gives the set of defining relations of the algebra $P\CD (\CA )$ together with Theorem \ref{23Jun19}.(2).

\begin{theorem}\label{2Aug19}
Let a Poisson algebra $\CA= S^{-1}(P_n/I)$ be a domain of essentially finite type over the field $K$ of characteristic zero, $d=r(\CC_\CA )$, $\phi : \CU (\CA ) \ra \bCU (\CA )=\CU (\CA ) / \kappa_\CA$, $u\mapsto u+\kappa_\CA$ and $\overline{\pi}_\CA : \bCU (\CA ) \ra P\CD (\CA )$, $a+\kappa_\CA \mapsto a$, $\d_i+\kappa_\CA \mapsto \{ x_i, \cdot \}$ where $a\in \CA$. Then 
\begin{enumerate}
\item $\ker (\pi_\CA ) = \phi^{-1}(T_{\bCU (\CA )})=\phi^{-1}(\tor_{S_{\mu (\i , \j )}}(\bCU (\CA )))=\phi^{-1}(\rann_{\bCU}(\gc_{\CA , d}^t))=\phi^{-1}(\lann_{\bCU}(\gc_{\CA , d}^t))$ for all $\i \in \mI_\CA (d)$, $\j \in \mJ_\CA (d)$ and $t\gg 0$. 
\item  $\ker (\overline{\pi}_\CA ) = T_{\bCU (\CA )}= 
\tor_{S_{\mu (\i , \j )}}(\bCU (\CA ))= \rann_{\bCU (\CA )}(\gc_{\CA , d}^t) = \lann_{\bCU (\CA )}(\gc_{\CA , d}^t)$ for all $\i \in \mI_\CA (d)$, $\j \in \mJ_\CA (d)$ and $t\gg 0$. 
\end{enumerate}
\end{theorem}

{\it Proof}. 1. Since $\kappa_\CA \subseteq \ker (\pi_\CA )$ (see (\ref{kAkpA})), $\ker (\pi_\CA )=  \phi^{-1}(\ker (\overline{\pi}_\CA ) )$ and statement 1 follows from statement 2. 

2. Notice that $\phi (\ker (\pi_\CA ))=  \ker (\overline{\pi}_\CA ) $ and $\mu = \mu (\i , \j )\in \gc_{\CA , d}$ for all $\i \in \mI_\CA (d)$ and  $\j \in \mJ_\CA (d)$. Now, by Theorem \ref{26Jul19}.(1),  $$\ker (\overline{\pi}_\CA )\subseteq \tor_{S_\mu } (\bCU ).$$
The epimorphism $\overline{\pi}_\CA$ is an $\CA$-module homomorphism $(\overline{\pi}_\CA (a\bu ) = a\overline{\pi}_\CA (\bu )$ for all elements $a\in \CA$ and $\bu \in \bCU$) and the element $\mu$ is a regular element of the algebra $P\CD (\CA )$. Therefore,  $$\ker (\overline{\pi}_\CA )\supseteq \tor_{S_\mu } (\bCU ),$$ i.e.  $\ker (\pi_\CA )= \tor_{S_\mu } (\bCU )$ and statement 2 follows from Proposition \ref{B2Aug19}. $\Box $\\

{\bf Criterion for $\ker (\pi_\CA )= \kappa_\CA$ $(\Leftrightarrow \ker (\overline{\pi}_\CA)=0$). } 
The next corollary follows at once from (\ref{CAkA1}) and Theorem, \ref{2Aug19}.(2) and is  a criterion for $\ker (\pi_\CA )= \kappa_\CA$ $(\Leftrightarrow \ker (\overline{\pi}_\CA)=0$). 
\begin{corollary}\label{a2Aug19}
We keep the assumption  of Theorem \ref{2Aug19}. Then  $\ker (\pi_\CA )= \kappa_\CA$  iff $\ker (\overline{\pi}_\CA)=0$, i.e. $\bCU (\CA ) = P\CD (\CA )$, iff the nonzero elements of the algebra $\CA$ are (left or right) regular in $\bCU (\CA )$ iff all/some of the elements of the set $\{ \mu (\i , \j ) \, | \, \i\in \mI_\CA (d) , \j \in \mJ_\CA (d)\}$ are (left or right) regular in $\bCU (\CA )$ iff $\rann_{\bCU (\CA )}(\gc_{\CA , d}^t)=0$ for all $t\gg 0$ iff $\lann_{\bCU (\CA )}(\gc_{\CA , d}^t)=0$ for all $t\gg 0$. 
\end{corollary}

{\bf The right kernel $\O_\CA'$ of the pairing $\Der_K(\CA ) \times \O_\CA \ra \CA$.} 
Let a Poisson algebra $\CA= S^{-1}(P_n/I)$ be an algebra  of essentially finite type. Recall $\O = \O_\CA\simeq \bigoplus_{i=1}^n\CA \d_i/G_\CA$ where $G_\CA := \sum_{j=1}^m \CA \d_{f_j}$ and $\d_{f_j}=\sum_{i=1}^n\frac{\der f_j}{\der x_i}\d_i$, and $\Der_K(\CA ) \simeq \Hom_\CA (\O , \CA )$, an isomorphism of left $\CA$-modules. So, there is a {\em pairing} of left $\CA$-modules  (which is  an $\CA$-bilinear map):
\begin{equation}\label{PairDer}
\Der_K(\CA ) \times \O \ra \CA , \;\; (\der , \o ) \mapsto (\der, \o ) := \der (\o ).
\end{equation}
If $\der=\sum_{i=1}^na_i\der_i\in \Der_k(\CA )$ (where $\der_i=\frac{\der}{\der x_i}$) and $\o =\sum_{i=1}^n b_i\d_i$ where $a_i, b_i\in \CA$ then 
\begin{equation}\label{PairDer1}
 (\der, \o ) = \sum_{i=1}^n a_ib_i.
\end{equation}
In particular, $(\der, \d_i)=a_i:={\rm coef}_i(\der)$, the $i^{th}$ {\em coefficient} of $\der$. The $\CA$-submodule of $\O$ 
$$ \O' = \O_\CA':=\{ \o \in \O \, | \, (\Der_K(\CA ) , \o )=0\}$$
is called the {\em right kernel} of the pairing. In a similar way, the left kernel is defined. Since $\Der_K(\CA ) \simeq \Hom_\CA (\O , \CA )$, the left kernel of the pairing is equal to zero. 

\begin{lemma}\label{x3Aug19}
Let a Poisson algebra $\CA= S^{-1}(P_n/I)$ be a domain  of essentially finite type  and let $T$ be a multiplicative subset of $\CA$. Then $\O'_{T^{-1}\CA} \simeq T^{-1} \O'_{\CA }$. 
\end{lemma}

{\it Proof}. Since $\Der_K(T^{-1} \CA ) \simeq \Hom_{T^{-1} \CA }(\O_{T^{-1} \CA }, T^{-1} \CA )\simeq \Hom_{T^{-1} \CA }(T^{-1}\O_ \CA , T^{-1} \CA )\simeq \;\;\;$ $T^{-1}\Hom_{ \CA }(\O_ \CA ,  \CA )\simeq T^{-1}\Der_K(\CA )$, we must have $\O'_{T^{-1}\CA} \simeq T^{-1} \O'_{\CA }$.  $\Box $

\begin{proposition}\label{A3Aug19}
Let a Poisson algebra $\CA= S^{-1}(P_n/I)$ be a domain of essentially finite type. Then 
\begin{enumerate}
\item $(\O'_\CA)_{\D (\i , \j )}=0$ for all $\i \in \mI_r$ and $\j \in \mJ_r$. 
\item $\ga_r^s\O'_\CA =0$ for some natural number $s\geq1 $. 
\item If, in addition, the algebra $\CA$ is a regular algebra than $\O_\CA'=0$. 
\end{enumerate}
\end{proposition}

{\it Proof}. 1. Let $\D = \D (\i , \j )$, $\j = (j_1, \ldots , j_r)$ and $\{ j_{r+1}, \ldots , j_n\} = \{ 1, \ldots , n\} \backslash \{ j_1, \ldots , j_r\}$.  Statement 1 follows at once from the fact that $ (\O_\CA)_\D =\O_{\CA_\D} =\bigoplus_{\nu = r+1}^n \CA_\D \d_{j_\nu }$ and $(\O_\CA')_\D \simeq \O'_{\CA_\D}$.

2. The left $\CA$-module $\O'_\CA$ is finitely generated. Then, by statement 1, $\D(\i , \j )^k\O'_\CA =0$ for some $k=k(\i , \j )$. Hence, $\ga_r^s\O'_\CA =0$ since the ideal $\ga_r$ is generated by the finite set $\{ \D (\i , \j )\}$. 

3. Statement 3 follows from statement 2 and the fact that $\ga_r=\CA$ (since the algebra $\CA$ is a regular algebra).  $\Box $\\

{\bf Criteria for $\ker (\pi_\CA )=0$.} In the case when the Poisson algebra $\CA$ is a regular domain of essentially finite type, Theorem \ref{3Aug19} is an efficient explicit criterion for $\ker (\pi_\CA )=0$, i.e. for the epimorphism $\pi_\CA : \CU (\CA ) \ra P\CD (\CA )$ to be an isomorphism.\\

{\bf Proof of Theorem \ref{3Aug19}}. $(1\Rightarrow 2)$  This implication is obvious since $\kappa_\CA \subseteq \ker (\pi_\CA )$. 

$(2\Rightarrow 1)$  Since $\kappa_\CA =0$,  $\bCU (\CA ) = \CU (\CA )$, and  the implication follows from Theorem \ref{2Aug19}.(1) and the fact that the algebra $\CU (\CA )$ is a domain (Theorem \ref{29Jul19}): $\ker (\pi_\CA ) =\ker (\overline{\pi}_\CA ) = T_{\bCU (\CA )}=0$. 

$(2\Rightarrow 3)$  If $\kappa_\CA =0$ then $\ker (\pi_\CA ) =0$ (since 
 $(1\Leftrightarrow 2)$), i.e. $\CU (\CA ) \simeq P\CD (\CA )$. By Theorem \ref{29Jul19} and Proposition \ref{B26Jul19}.(2), 
 $$ 2\GK (\CA ) = \GK  (\CU (\CA )) = \GK \, P\CD (\CA ) = \GK (\CA ) +d,$$
we have that $d=\GK (\CA) = n-r$. Now, the implication follows from Proposition \ref{A3Aug19}.(3).

$(2\Leftarrow 3)$ Suppose that $d=n-r$. Then $\kappa_\CA =0$ iff all the elements $\d_{\i', i'_\mu; \j'}$ (in statement 3) belong to the $\CA$-module $G_\CA$, by Corollary \ref{a2Aug19} and  (\ref{KPdif5}), iff $(\Der_K(\CA ) , \d_{\i', i'_\mu; \j'})=0$ for all $\d_{\i', i'_\mu; \j'}$ since $\O'_\CA=0$, by Proposition \ref{A3Aug19}.(3),  iff $(\der_{\i ; \j , j_\nu}, \d_{\i' ; \j' , j_\nu'})=0$ for all elements $
\der_{\i ; \j , j_\nu}$ and $ \d_{\i' ; \j' , j_\nu'}$ as in statement 3 since the elements $\der_{\i ; \j , j_\nu}$ are $\CA$-module generators of $\Der_K(\CA )$, \cite[Theorem 1.1]{gendifreg}. Therefore, $(2\Leftarrow 3)$.  $\Box $ \\

Theorem \ref{A5Aug19} is another criterion for $\ker (\pi_\CA )=0$ where we do not assume that the algebra $\CA$ is regular. We need the following lemma. 

\begin{lemma}\label{a7Aug19}
Let $R$ be a commutative domain of essential finite type over a field of characteristic zero and $\{ \CD (R)_i\}_{i\geq 0}$ be the order filtration on $\CD (R)$. Then 
\begin{enumerate}
\item Multiplication by a  nonzero element of $R$ preserves the order of  differential operator. 
\item The ring $\CD (R)$ is a domain.
\end{enumerate}
\end{lemma}

{\it Proof}. 1. We use an induction on the order of differential operators. Given $\d \in \CD (R)_i\backslash  \CD (R)_{i-1}$. We have to show that $r\d , \d r\in \CD (R)_i\backslash  \CD (R)_{i-1}$ for all $r\in R\backslash \{ 0\}$. This is obvious for $i=0$ as $\CD (R)_0=R$ is a domain. Suppose that $i>0$ the result is true for all $i'<i$. Suppose that $r\d \in \CD (R)_{i-1}$ (resp., $ \d r \in \CD (R)_{i-1}$) for some $r\neq 0$, we seek a contradiction. Fix $r'$ such that $\d' = [r', \d ] \in \CD (R)_{i-1}\backslash \CD (R)_{i-2}$. Then by induction
$$ \CD (R)_{i-1}\backslash \CD (R)_{i-2}\ni r\d' =  [r', r\d ]\in \CD (R)_{i-2}$$ (resp., $ \CD (R)_{i-1}\backslash \CD (R)_{i-2}\ni \d'r=[\d' r, r' ]\in \CD (R)_{i-2}$), a contradiction. 

2. Take a nonzero element $s$ of the Jacobian ideal $\ga_r$ of $\CA$. By statement 1, $\CD (\CA ) \subseteq \CD (\CA )_s\simeq \CD (\CA_s)$. The algebra $\CA_s$ is a regular domain of essentially finite type, hence the algebra $\CD (\CA_s)$ is a domain and so its subalgebra $\CD (\CA )$.  $\Box $\\

{\bf Proof of Theorem \ref{A5Aug19}}. $(1\Rightarrow 2)$ If statement 1 holds then $\CU (\CA ) = P\CD (\CA )$.  Then the implication follows from the inclusion  $\kappa_\CA \subseteq \ker (\pi_\CA )$ and the fact that $\CU (\CA ) = P\CD (\CA )$ is a domain (Lemma \ref{a7Aug19}.(2)). 

$(2\Rightarrow 1)$  Since $\kappa_\CA =0$, we have that  $\bCU (\CA ) = \CU (\CA )$.
 Now, the implication follows from Theorem \ref{2Aug19}.(2) and the assumption that the algebra $\CU (\CA )=\bCU (\CA )$ is a domain: $\ker (\pi_\CA ) = \ker (\overline{\pi}_\CA ) = T_{\bCU (\CA )}=0$. 
 
$(2\Leftrightarrow 3)$  If $\kappa_\CA =0$ and $\CU (\CA )$ is a domain then   $\ker (\pi_\CA ) =0$ (since 
 $(1\Leftrightarrow 2)$), i.e. $\CU (\CA ) \simeq P\CD (\CA )$. By Theorem \ref{29Jul19}.(4) and Proposition \ref{B26Jul19}.(1), 
 $$ 2\GK (\CA ) = \GK  (\CU (\CA )) = \GK \, P\CD (\CA ) = \GK (\CA ) +d,$$
we have that $d=\GK (\CA) = n-r$. 

Suppose that $d=n-r$ and $\CU (\CA )$ is a domain. Then $\kappa_\CA =0$ iff all the elements $\d_{\i', i'_\mu; \j'}$ (in statement 3) belong to the $\CA$-module $G_\CA$, by (\ref{KPdif5}), iff $(\Der_K(\CA ) , \d_{\i', i'_\mu; \j'})=0$ for all $\d_{\i', i'_\mu; \j'}$, by Proposition \ref{A3Aug19}.(3) (since $\O'_\CA=0$).   $\Box $\\

{\bf Criterion for the homomorphism $\pi_\CA : \CU (\CA ) \ra \CD (\CA )$ to be an isomorphism.} For a regular domain of essentially  finite type $\CA$, Theorem \ref{5Aug19} is a criterion for the  homomorphism $\pi_\CA : \CU (\CA ) \ra \CD (\CA )$ to be an isomorphism.\\

{\bf  Proof of Theorem \ref{5Aug19}}. $(1\Rightarrow 2)$ Suppose that the  homomorphism $\pi_\CA : \CU (\CA ) \simeq \CD (\CA ))$ is an isomorphism.  Then $P\CD (\CA ) = \CD (\CA)$ and the epimorphism $\pi_\CA : \CU (\CA ) \ra P\CD (\CA )$ is an isomorphism. By Theorem \ref{22Jul19}, it suffices to show that $\Der_K(\CA ) =\CA \CH_\CA$. The equality follows from the Claim.\\

{\sc Claim.} {\em For all} $i\geq 1$, $P\CD (\CA )_{\leq i} = P\CD (\CA )_i$. \\ 

Indeed, for $i=1$,  $P\CD (\CA )_{\leq 1} = \CA \oplus \CA \CH _\CA$ and $P\CD (\CA )_i = \CA \oplus P\CD (\CA ) \cap \Der_K(\CA )=\CA \oplus \CD (\CA ) \cap \Der_K (\CA ) = \CA \oplus \Der_K(\CA )$, so the equality in the Claim for $i=1$ yields the equality $\CA \CH_\CA = \Der_K(\CA )$.

{\it Proof of the Claim.}  Recall that 
 $P\CD (\CA )_{\leq i} \subseteq  P\CD (\CA )_i$ for all $i\geq 1$.

Since the algebra $\CA$ is a regular domain,  the algebra $\CD (\CA )=\CU (\CA )=P\CD (\CA )$ is a simple Noetherian domain. 
 By Theorem \ref{XBA29Jul19}, the algebra ${\rm gr}\, \CU (\CA)=\bigoplus_{i\geq 0} P\CD (\CA )_{\leq i}/P\CD (\CA )_{\leq i-1 }$ is a domain. By Theorem \ref{26Jul19}.(2), there are natural numbers $s_i$  such that $$\gc_{\CA , d}^{s_i}P\CD (\CA )_i\subseteq P\CD (\CA )_{\leq i}\;\; {\rm for \; all}\;\; i\geq 0.$$
Since $\CD (\CA )$ is a domain, the degree of a differential operators is preserved by the multiplication by nonzero element of $\CA$, and the Claim follows  (since $\CD (\CA )_i\backslash \CD (\CA )_{i-1}\subseteq P\CD (\CA )_{\leq i} \backslash P\CD (\CA )_{\leq i-1}$, by the inclusions above and the fact that the algebra ${\rm gr}\, \CU (\CA)$ is a domain). 

$(2\Rightarrow 1)$ By Theorem \ref{26Jul19}.(3), $\ker (\pi_\CA ) = \kappa_\CA$, i.e. $\CU (\CA ) = \bCU (\CA )$. Since $d=n-r$ and $\gc_{\CA , d} = \CA$, $\ker (\overline{\pi}_\CA )=0$, by Theorem \ref{2Aug19}, and so $$\CU (\CA ) \simeq P\CD (\CA ).$$ Since $d=n-r$ and $\gc_{\CA , d}=\CA$, $\Der_K(\CA ) = \CA \CH_\CA$, by Theorem \ref{22Jul19}. The algebra $\CA$ is a regular algebra, hence the algebra $\CD (\CA )$ is generated by $\CA$ and $\Der_K(\CA )$, and so $\CD (\CA ) = P\CD (\CA )$. $\Box$


\section{Simplicity criteria for the algebras $\CU (\CP )$ and $ P\CD (\CP )$ }\label{SIMUPDP}

The aim of this section is to give/prove several simplicity criteria for the algebras  $ P\CD (\CP )$ (Theorem \ref{X22Jul19})  and $\CU (\CP )$ (Theorem \ref{Y22Mar19} and  Theorem \ref{7Aug19}). \\

{\bf Subalgebras $\L$ of $\CD (R)$ such that $R\subset \L $.} Let $R$ commutative $K$-algebra over an {\bf arbitrary} field $K$. Let $\L$ be a subalgebra of the algebra $\CD (R)$ of differential operators on $R$ such that properly contains $R$ (i.e. $R\subseteq \L $ and $R\neq \L$). Then $\L =\bigcup_{i\geq 0} \L_i$ is a filtered algebra where $\{ \L_i:= \L \cap \CD (R)_i\}_{i\geq 0}$ is the filtration that is induced by the order filtration of the algebra $\CD (R)$.  We call this filtration the {\em order filtration} on $\L$. The algebra $\L$ is an $R$-bimodule.

\begin{lemma}\label{a21Jul19}
Let $\L$ and $\CD (R)$ be as above. Then $\L_1= R\oplus \G$ where $\G := \L \cap \Der_K(R)$ is a left $R$-submodule  and a Lie subalgebra of $\Der_K(R)$. 
\end{lemma}

{\it Proof}. Since $R\neq \L$, we must have $\L_1\backslash R\neq \emptyset$. Since $\CD (R)_1=R\oplus \Der_K(R)$ and $R\subseteq \L$, 
$$\L_1= \L \cap \CD (R)_1= \L \cap (R\oplus \Der_K(R))=R\oplus \L \cap \Der_K(R)=R\oplus \G $$
and $\G$ is a nonzero left $R$-submodule and  a Lie subalgebra of $\Der_K(R)$.  $\Box $\\

{\it Definition.} The subalgebra $\D (\L ) =\langle R, \G \rangle$ of $\L$  is called  the {\em derivation subalgebra} of $\L$ where $\G = \L \cap \Der_K(R)$. 

\begin{theorem}\label{A21Jul19}
Let $R$ be a commutative $K$-algebra, $K$ be an arbitrary field, $\L$ be a subalgebra of the algebra $\CD (R)$ of differential operators on $R$ that properly contains $R$. If $I$ is a nonzero ideal of $\L$ then the intersection $I\cap R$ is a nonzero $\G$-stable ideal of $R$ where $\G = \L\cap \Der_K(R)$ (i.e. $\G (I\cap R) \subseteq I\cap R$, the elements of $\G$ act as derivations) such that $\L (I\cap R)\L \cap R=I\cap R$. 
\end{theorem}

{\it Proof}.  Recall that $\G$ is a nonzero left $R$-module and a  Lie subalgebra of $\Der_K(R)$ (Lemma \ref{a21Jul19}).

(i)  $I\cap  R\neq 0$: Notice that $I=\bigcup_{i\geq 0} I_i$ where $I_i=I\cap \L_i= I\cap \CD (R)_i$. Let $s=\min \{ t\geq 0\, | \, I_t\neq 0\}$. Choose a nonzero element, say $v$, of $I_s$. Then, for all elements $r\in R$, $[a, v]\in I_{s-1}=0$, i.e. $v\in R$, as required.

(ii) $I\cap R$ {\em is a $\G$-stable ideal of the algebra $R$}: For all elements $u\in I$ and $ \g \in \G$, $I\ni [ \g , u] = \g (u)$, and the statement (ii) follows.

(iii) $\L (I\cap R)\L \cap R=I\cap R$: $I\cap R\subseteq \L (I\cap R)\L \cap R \subseteq I \cap R$ and the statement (iii) follows.
 $\Box$ \\

{\it Definition.} Let $\G'$ be a set of $K$-linear maps from $R$ to $R$. We say that the algebra $R$ is $\G'$-{\em simple} if $\{ 0\}$ and $R$ are the only $\G'$-stable ideals of $R$. \\



{\bf Simplicity criterion for the algebra $P\CD (\CP )$ of Poisson differential operators on $\CP$.}\\


{\bf Proof of Theorem \ref{X22Jul19}}. 
 Without loss of generality we may assume that  $\CH_\CP\neq 0$ (since otherwise the theorem is obvious as $P\CD (\CP )= \CP$).

(i) {\em The algebra $P\CD (\CP )$ is not simple $\Rightarrow$ the Poisson algebra $\CP$ is not simple}: The algebra $P\CD (\CP )$ is a subalgebra of $\CD (\CP )$ that properly contains the algebra $\CP$ (since $\CH_\CP\neq 0$). Suppose that $I$ is a proper ideal of $P\CD (\CP )$. Then, by Theorem \ref{A21Jul19}, $I\cap \CP$ is a proper $\G$-stable ideal of $\CP$ where $\G = P\CD (\CP ) \cap \Der_K(\CP )\supseteq \CH_\CP =\{ \pad_a=\{ a\cdot \} \, | \, a\in \CP \}$. So, the intersection $I\cap \CP$ is a proper $\CH_\CP$-stable ideal of $\CP$, i.e. the Poisson algebra $\CP$ is not simple. 

(ii) {\em The Poisson algebra $\CP$ is not simple $\Rightarrow$ the  algebra $P\CD (\CP )$ is not a  simple algebra}: If $J$ is a proper Poisson ideal of $\CP$ then $J\cdot P\CD (\CP )$ is a proper ideal of $P\CD (\CP )$ (since $\CH_\CP (J ) \subseteq  J$). So, statements 1  and 2 are equivalent.  $\Box$

\begin{corollary}\label{xX22Jul19}
Let $\CP$ be a Noetherian Poisson algebra over  a field $K$ of characteristic zero.  If the algebra $P\CD (\CP )$ is a simple algebra ($\Leftrightarrow$ the Poisson algebra $\CP$ is a Poisson simple algebra, Theorem \ref{X22Jul19}) then  the algebra $\CP$ is a domain.
\end{corollary}

{\it Proof}. Since the minimal primes of a Noetherian  algebra over a field of characteristic zero are derivation-stable and the Poisson algebra $\CP$ is Poisson simple, the Poisson algebra $\CP$ must be a domain.  $\Box $\\

{\bf Simplicity criteria for the Poisson    enveloping algebra $\CU (\CP )$.}


{\bf Proof of Theorem \ref{Y22Mar19}}. $(1\Rightarrow 2 )$ The Poisson algebra $\CP$ is a left $\CU (\CP )$-module. Let $\ga$ be the kernel of the algebra epimorphism $\CU  (\CP )\ra P\CD (\CP )$, $p\mapsto p$, $\d_q\mapsto \pad_q$ for all $p,q\in \CP$.  Now, the implication is obvious.

$(2\Rightarrow 1 )$  The implication is obvious.

$(2\Leftrightarrow 3 )$ The equivalence follows from Theorem \ref{X22Jul19}. $\Box $\\

In the case when the Poisson algebra $\CP = \CA$ is an algebra  of essentially finite type over a field of characteristic zero,  Theorem \ref{Y22Mar19} can be strengthen. \\


{\bf Proof of Theorem \ref{7Aug19}}. $(1\Rightarrow 2)$ By Theorem \ref{Y22Mar19},  the algebra $P\CD (\CA )$ is a simple algebra and $\CU (\CA )\simeq P\CD (\CA )$. 
   The algebra $\CA$ is of essentially finite type (hence Noetherian)  of characteristic zero.   Hence, the Jacobian ideal $\ga_r$ and all the minimal primes of the algebra $\CA$ are $\Der_K(\CA )$-stable ideals 
(\cite[Theorem 5]{Matsumura-1982} and
  \cite[Theorem 1]{Seidenberg-1967}, respectively).  Since 
  the algebra $\CU (\CA )\simeq P\CD (\CA )$ is simple,  the algebra $\CA$ is a regular domain. Theorem \ref{3Aug19} is a criterion for $\CU (\CA ) \simeq P\CD (\CA )$ when the algebra $\CA$ is a regular domain, and statement 2 follows.
  

$(1\Leftarrow 2\Leftrightarrow 3)$  The implications  follow from Theorem \ref{Y22Mar19} and Theorem \ref{3Aug19}. 
  $\Box$\\

In particular, in the case when the algebra $\CA$ is a regular domain Theorem \ref{7Aug19} states that\\

$\bullet$  {\em the algebra $\CU (\CA )$ is a simple algebra iff the Poisson algebra $\CA $ is a Poisson simple algebra, $d=n-r$ and $(\der_{\i ; \j , j_{\nu }}, \d_{\i', i_\mu'; \j'})=0$ for all elements $\i \in \mI_r$, $\j \in \mJ_r$, $\i'\in \mI_\CA (d)$, $\j' \in \mJ_\CA (d)$, $\nu = r+1, \ldots , n$ and $\mu = d+1, \ldots , n$ where 
for $\j = \{ j_1,\ldots , j_r\}$ and $\i'=\{ i_1', \ldots , i_d'\}$we have that $\{ j_{r+1}, \ldots , j_n\}:=\{ 1,\ldots , n\} \backslash \{ j_1,\ldots , j_r\}$ and} $\{ i_{d+1}', \ldots , i_n'\} := \{ 1, \ldots , n\} \backslash
\{ i_1', \ldots , i_d'\}$. \\

So, Theorem \ref{7Aug19} is an efficient tool in proving or disproving simplicity of the algebra $\CU (\CA )$. \\




\small{

Department of Pure Mathematics

University of Sheffield

Hicks Building

Sheffield S3 7RH

UK

email: v.bavula@sheffield.ac.uk}

\end{document}